\documentclass{article}
\usepackage{iclr2027_conference,times}

\usepackage{amsmath,amsfonts,bm}

\def\eqref#1{equation~\ref{#1}}

\def\1{\bm{1}}

\DeclareMathAlphabet{\mathsfit}{\encodingdefault}{\sfdefault}{m}{sl}
\SetMathAlphabet{\mathsfit}{bold}{\encodingdefault}{\sfdefault}{bx}{n}

\usepackage{hyperref}
\hypersetup{colorlinks=true,linkcolor=[rgb]{0,0,0.6},citecolor=[rgb]{0,0.4,0},urlcolor=[rgb]{0,0,0.6},pdftitle={Learning Warm Starts for the Simplex Method from Measured Solver Cost},pdfauthor={Felix Parker, Danielle A. Ripsman, Emmett Springer, Kimia Ghobadi}}
\usepackage{url}
\usepackage{xurl}
\usepackage{needspace}
\usepackage{amsmath,amssymb,amsthm}
\usepackage{booktabs}
\usepackage{multirow}
\usepackage{graphicx}
\usepackage{xcolor}
\usepackage{tikz}
\usetikzlibrary{arrows.meta,positioning,calc,shapes.misc,fit,backgrounds}
\usepackage{caption}
\usepackage{subcaption}
\usepackage{enumitem}
\usepackage{microtype}
\usepackage{siunitx}
\usepackage{xspace}
\usepackage{algorithm}
\usepackage{algpseudocode}
\usepackage{booktabs}
\usepackage{siunitx}
\usepackage{array}
\usepackage{longtable}
\usepackage{url}

\setlist{nosep,leftmargin=1.4em}

\theoremstyle{definition}
\newtheorem{theorem}{Theorem}
\newtheorem{corollary}{Corollary}

\newcommand{\ours}{\textsc{CostStart}\xspace}
\newcommand{\oursS}{\textsc{CostStart-Shift}\xspace}
\newcommand{\oursA}{\textsc{CostStart-Alg}\xspace}
\newcommand{\oursF}{\textsc{CostStart-Full}\xspace}

\newcommand{\Trs}{\textsc{Train-Small}\xspace}
\newcommand{\Trtw}{\textsc{Train-20k}\xspace}
\newcommand{\Dtw}{\textsc{Dev-20k}\xspace}
\newcommand{\Dfif}{\textsc{Dev-50k}\xspace}
\newcommand{\Ttw}{\textsc{Test-20k}\xspace}
\newcommand{\Tfif}{\textsc{Test-50k}\xspace}
\newcommand{\Tei}{\textsc{Test-80k}\xspace}
\newcommand{\Tnc}{\textsc{Test-NewClass}\xspace}
\newcommand{\Tad}{\textsc{Test-Adapt}\xspace}
\newcommand{\Tva}{\textsc{Test-Variant}\xspace}
\newcommand{\ci}[2]{{\scriptsize[#1,\,#2]}}

\DeclareMathOperator{\nnz}{nnz}


\title{Learning Warm Starts for the Simplex Method\\from Measured Solver Cost}

\author{Felix Parker\thanks{Equal contribution.}\hphantom{$^{*}$}, Danielle A. Ripsman\footnotemark[1]\hphantom{$^{*}$}, Emmett Springer \& Kimia Ghobadi\\
Center for Systems Science and Engineering, Johns Hopkins University}

\iclrfinalcopy %

\begin{document}

\maketitle
\lhead{Preprint} %

\begin{abstract}
Learned warm starts for the simplex method build starting bases using predictions of which variables will be basic at optimality. Prior work has shown that starting bases predicted by graph neural networks can sharply reduce iteration counts, but the time saved may fall well short of the iterations saved. On 25 classes of linear programs (LPs), we find that a strong supervised basis predictor removes two thirds of the dual simplex's iterations but only a third of its time. Two effects explain the gap: the iterations removed by a good start tend to be cheap, and closeness to an optimal basis does not rank starts by cost. We prove that two starts with equal overlap with an optimal basis and equal pivot counts can differ in linear-algebra work by a factor that grows with the square root of the number of rows. Empirically, across the starts of one LP, the pivot count explains 71\% of the variance of log solve time, while the CPU instruction count explains 98.5\%. We introduce \ours, which uses supervised prediction for the basis and trains a small policy on measured instruction counts to choose, per LP, how many predicted columns to use, which algorithm and pricing rule to run, and whether to warm-start at all. On unseen LPs of the training classes, including LPs four times the training size, \ours roughly halves the solve time relative to the solver's default start and is 1.4 to 2.2 times as fast as two published methods retrained on the same data.

\end{abstract}

\section{Introduction}

Despite decades of competing algorithmic development since its 1947 debut, the simplex method remains a highly practical algorithm for solving linear programs (LPs) today.
Simplex pivots iteratively from a starting basis to neighboring bases until it reaches optimality. A near-optimal starting basis can significantly shorten this path, which solvers have long exploited by building starting bases with crash heuristics \citep{bixby1992implementing,maros2003computational}. Recent work has had success learning these starting bases, using a graph neural network (GNN) to score each column for membership in an optimal basis and a decoder to turn the scores into a valid basis. These methods substantially reduce the average number of simplex iterations \citep{fan2023smart,lu2025learning}. Learned starts are most useful for a stream of new LPs from known problem classes that must be solved by simplex, for example because a basic solution is needed at the root of branch-and-bound, for cut generation or for sensitivity analysis, or because runs must be reproducible. In that setting no earlier basis of the same LP is available, the training cost is paid once, and what matters is the solve time saved per LP.

The number of iterations (pivots) is a natural measure of starting-basis quality: it is hardware- and implementation-independent, and, intuitively, a start closer to an optimal basis should need fewer pivots.
In practice, however, users pay for solve time, and the work of each pivot depends on factors like the sparsity of both the basis and the vectors it computes, and the numerical stability of the basis factorizations. As a result, pivot costs may vary by orders of magnitude within a single solve \citep{hall2005hyper}. Published results show the difference: in the HiGHS experiments of \citet{fan2023smart}, the learned start takes $0.25\times$ the iterations of the default start but $0.41\times$ its time; \citet{lu2025learning} report $0.60\times$ and $0.73\times$ and note that ``a closer initial basis does not always result in greater acceleration'' (Table~\ref{tab:published}). We therefore minimize ``charged'' time, the overall solve time with every overhead of the learned start included, and treat the iteration count as only one of its factors.

\begin{figure}[t]
\centering
\includegraphics[width=\linewidth]{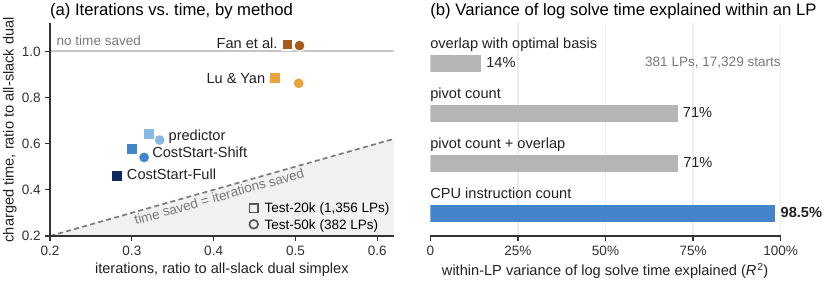}
\caption{Iteration savings overstate time savings. (a) Iterations and charged time of learned warm starts on unseen test LPs, relative to the dual simplex from a cold (all-slack) start (geometric means; \S\ref{sec:setup}). Every method lies well above the diagonal, where time saved would equal iterations saved. (b) How much of the variance in log solve time each quantity explains across different starts of the same LP. The CPU instruction count, which our method, \ours (\S\ref{sec:method}), trains on, explains nearly all of it.}
\label{fig:gap}
\end{figure}

On our data, the gap is large (Figure~\ref{fig:gap}a). On unseen LPs from 25 problem classes, the strongest supervised predictor we could train removes about two thirds of the dual simplex's iterations but only about a third of its time, and two published methods save little or no time. The gap has two causes. First, the iterations that a good start removes are often the cheap ones. Second, closeness to an optimal basis, which predictors are trained for, does not rank the candidate starts of an LP by cost. We prove that two starts that contain the same fraction of an optimal basis (their overlap) and need the same number of pivots can differ in linear-algebra work by a factor that grows with the square root of the number of rows. On our generated LPs, across different starts of the same LP, overlap explains just 14\% of the variance in log solve time, and even the pivot count explains only 71\% (Figure~\ref{fig:gap}b).

These findings call for a method that optimizes solve time, not just overlap with an optimal basis. Our method, \ours, has two parts. A supervised predictor scores each column, and a triangular crash decodes the scores into a basis that is nonsingular by construction; the resulting start alone is already faster than two published methods retrained on our data. A small decision policy then chooses, per LP, how far to shift the predicted scores (and so how many predicted columns enter the start), which simplex algorithm and pricing rule to run, or whether to skip the warm start entirely. The policy is trained on the CPU instruction count of every action on every training LP, a less noisy proxy for solve time. Its three versions, \oursS, \oursA and \oursF, have progressively larger action sets.

Whereas the published methods train one model per dataset, we train a single model on LPs from 26 problem classes with 134 variants, half of which are the dual forms of the other half, and test it on the 25 of these classes that have test instances. We then test how far it generalizes: to new instances of the training classes, to LPs up to four times the training size, to held-out variants, and to unseen classes. Guardrails limit losses on unfamiliar and badly scaled LPs.

\paragraph{Contributions.}
\begin{itemize}
\item \textbf{Analysis.} We explain why iteration savings overstate time savings for learned simplex warm starts: the removed iterations are cheap, and neither overlap nor pivot count ranks the starts of an LP by cost (\S\ref{sec:gap}).
\item \textbf{Method.} \ours pairs a supervised basis predictor with a decision policy trained on measured instruction counts, whose actions change the cost of each iteration as well as their number; guardrails limit losses on unfamiliar LPs (\S\ref{sec:method}).
\item \textbf{Evidence.} On fresh test sets, \ours roughly halves the solver's charged time, including on LPs up to four times the training size, and is 1.4 to 2.2 times as fast as two published methods retrained on our data. On new problem classes, adapting with 25 labeled LPs or an online bandit recovers part of the lost speed (\S\ref{sec:results}).
\item \textbf{Data and code.} We release our LP generator library (\url{https://github.com/flixpar/SyntheticLPs.jl}) and our code and trained models (\url{https://github.com/flixpar/learned-simplex-warm-start}); Appendix~\ref{app:repro} describes how to reproduce our results.
\end{itemize}

\section{Related work}
\label{sec:related}

\textbf{Learned initial bases for the simplex method.}
The canonical work on predicting initial simplex bases \citep{fan2023smart} trains a GNN on the constraint-variable bipartite graph of an LP \citep{gasse2019exact} to predict each structural and slack variable's status at optimality: basic, at lower bound, or at upper bound. The GNN is trained with inverse-frequency-weighted cross-entropy, and the predicted basis consists of the $m$ highest-scoring columns, repaired by a factorization-based procedure \citep{bixby1994recovering} when the resulting matrix is singular. \citet{lu2025learning} build on this approach with a duality-inspired tripartite graph, which carries objective coefficients and bounds as edges to a global node, and with richer supervision: an added loss term on the mistakes of the top-$m$ selection, multi-level labels sampled from feasible bases along the solver's own path, and label preprocessing for equality rows. Both works are our reference points; we reimplement and retrain both on identical data and report the reproduction.

\textbf{Crash bases and measures for warm starts.}
Crash bases built from the structure of the constraint matrix are classical \citep{bixby1992implementing, gould1989new, maros2003computational}, as is the observation that they can cut iterations without cutting time, because a denser basis slows each iteration and some pricing rules, such as steepest edge, are costly to initialize \citep{koberstein2005dual}. Both \citet{fan2023smart} and \citet{lu2025learning} report iterations and solve time relative to the solver's default start, with \citet{lu2025learning} observing that a closer initial basis does not always give greater acceleration. Evaluations of warm starts have recently come under scrutiny for interior-point solvers. On AC optimal power flow problems, \citet{suri2026warp} find that reported 30--46\% iteration reductions rest on benchmarking against a flat start rather than the solver's actual default, and that against the corrected baseline no primal-only method reduces iterations. \citet{taheri2026notall} independently report that most partial warm starts for the same problems increase solve time or reduce reliability. For the simplex method, we find that iteration savings similarly overstate the real gain, and we explain why.

\textbf{Algorithms with predictions.}
The learning-augmented algorithms literature bounds running time by the quality of a prediction. \citet{sakaue2022discrete} give a discrete-convex-analysis framework whose running times depend on the $\ell_\infty$ distance from the prediction to an optimum, covering weighted matroid intersection among other problems, later extended to M-convex minimization \citep{oki2023faster}. For the simplex method, we prove that closeness gives no such guarantee: two starts with equal overlap and equal pivot counts can differ in work by a factor that grows with problem size.

\section{Iteration savings overstate time savings}
\label{sec:gap}

We solve LPs $\min\{c^\top x : Ax = b,\ \ell \le x \le u\}$ with $m$ rows and $n$ structural columns with HiGHS 1.13.1 \citep{huangfu2018parallelizing}, on one thread with presolve off. A basis is a set of $m$ columns of $[A\ I]$ with a nonsingular submatrix; the columns of $I$ are slacks, and the all-slack basis is $I$ itself. A warm start installs a basis $B_0$ before the solve, and its overlap is $|B_0 \cap B^\star|/m$ for an optimal basis $B^\star$. Our main comparator, the native solver, receives no basis and uses HiGHS' ``choose'' strategy: primal simplex if the all-slack basis is primal feasible, dual simplex otherwise (HiGHS' default, always dual, is slower on our data).

Learned warm starts make $B_0$ resemble $B^\star$ and are evaluated through the number of iterations $N$. The user, however, also pays for producing and installing the start and for the work of each iteration. To compare a method $a$ with a reference $r$ on one LP, we therefore split the log time ratio into three terms: iterations, time per iteration and fixed overhead,
\begin{equation}
\label{eq:decomp}
\log\frac{T_a}{T_r} \;=\; \log\frac{N_a}{N_r} \;+\; \log\frac{\bar w_a}{\bar w_r} \;+\; \Big(\log\frac{T_a}{T^{\text{run}}_a} - \log\frac{T_r}{T^{\text{run}}_r}\Big),
\end{equation}
where $T$ is the charged time, $T^{\text{run}}$ the time of the solver's own run and $\bar w = T^{\text{run}}/N$ its mean time per iteration.

\subsection{The gap in published results}
\label{sec:published}

The gap between iteration savings and time savings is already visible in the published papers' own tables. Table~\ref{tab:published} lists, for every dataset and solver setting that \citet{fan2023smart} and \citet{lu2025learning} report, the ratio of the learned start's iterations and wall time to those of the solver's default start, computed from their tables. In 18 of the 21 settings the time ratio is larger than the iteration ratio, often by a wide margin; on STOCH with OptVerse, for example, the learned start takes $0.42\times$ the iterations but $0.96\times$ the time. These ratios come from different data, decoders and charging rules than ours, and in part from a different solver, so we do not compare them with our numbers. They show only that the gap is not specific to our setting.

\begin{table}[t]
\centering\footnotesize
\caption{Iteration and time ratios (learned start / default start) computed from the published tables: \citet{fan2023smart}, Tables 1 (OptVerse) and 2 (HiGHS), dual simplex; \citet{lu2025learning}, Tables 3 to 7 (HiGHS dual simplex, their tripartite model). Geometric means over the datasets of each column. $^\dagger$ time falls more than iterations.}
\label{tab:published}
\setlength{\tabcolsep}{5pt}
\begin{tabular}{@{}lcc@{\hspace{14pt}}lcc@{}}
\toprule
\multicolumn{3}{@{}l}{\citet{fan2023smart}: iterations / time} & \multicolumn{3}{l}{\citet{lu2025learning}: iterations / time}\\
\cmidrule(r){1-3}\cmidrule(l){4-6}
Dataset & OptVerse & HiGHS & Dataset & unpresolved & presolved\\
\midrule
LIBSVM & 0.61 / 0.66 & 0.58 / 0.62 & Mirp1 & 0.42 / 0.52 & 0.32 / 0.42\\
MIRP & 0.64 / 0.70 & 0.61 / 0.81 & Mirp2 & 0.93 / 1.05 & 0.86 / 0.83$^\dagger$\\
STOCH & 0.42 / 0.96 & 0.48 / 0.35$^\dagger$ & Anonymous & 0.59 / 0.72 & 0.74 / 1.19\\
GEN & 0.23 / 0.38 & 0.04 / 0.15 & Load balancing & 0.46 / 0.71 & 0.63 / 0.70\\
SC-1 & 0.10 / 0.29 & 0.25 / 0.52 & & &\\
SC-2 & 0.14 / 0.25 & 0.18 / 0.33 & 5 large MIRP LPs & \multicolumn{2}{c}{0.93 / 0.86$^\dagger$}\\
\midrule
Geometric mean & 0.28 / 0.48 & 0.25 / 0.41 & Geometric mean & 0.57 / 0.73 & 0.60 / 0.73\\
\bottomrule
\end{tabular}
\end{table}

\subsection{Overlap and pivot count do not determine work}
\label{sec:theory}

An optimal basis would make the solve nearly free, but for other starts, overlap and pivot count do not determine the work, even for a fixed pivot rule. Take classical dual steepest edge (DSE), which picks the leaving row $p$ maximizing $\delta_p^2/\gamma_p$, where $\delta_p$ is the row's primal infeasibility and $\gamma_p = \|B^{-\top}e_p\|_2^2$. Each iteration solves for $\rho_p = B^{-\top}e_p$, $\alpha_q = B^{-1}a_q$ and, for DSE, $\tau = B^{-1}\rho_p$. We call the total number of nonzeros of these vectors along a path its solve-vector work $K$. $K$ depends only on the bases on the path and bounds the work of any implementation from below; on the LP below, output-sensitive solves and heap-based row selection achieve $\Theta(M + K)$ up to logarithmic terms (Appendix~\ref{app:proofs}).

\begin{theorem}[Equal overlap and pivot count, unbounded work ratio]
\label{thm:chain}
For every even $t \ge 4$ let $M = t^4$ and $r = t^3$, and consider the minimum-cost-flow LP on a path with return arcs $\min \sum_{i=1}^M u_i$ subject to $u_j - u_{j+1} - v_j = b_j$ ($1 \le j < M$), $u, v \ge 0$, with an explicit rational $b \in [0.97, 1.03]^{M-1}$. It has a unique nondegenerate optimal basis, and three dual-feasible starts satisfy, under exact DSE,
\begin{center}\small
\begin{tabular}{lccc}\toprule
Start & Overlap with $B^\star$ & Pivots & Work $K$\\ \midrule
Cold (all-slack: every $v_j$ basic) & $0$ & $M-1$ & $\sim \tfrac32 t^6$\\
Warm A (last $r$ chain arcs missing: singletons) & $1 - r/(M-1)$ & $r$ & $\sim 2t^5$\\
Warm B (every second of the last $2r$ arcs missing: pairs) & $1 - r/(M-1)$ & $r$ & $\sim t^7$\\ \bottomrule
\end{tabular}
\end{center}
Every choice on the three paths is unique, every basis factorizes without fill, and the two warm bases have equal nonzero counts after every pivot. Hence $K_B/K_A \sim t^2/2 = \Theta(\sqrt{M})$, and $K_B/K_{\text{cold}} = \Theta(M^{1/4})$.
\end{theorem}

\begin{corollary}[The proxy-optimal start can be the costlier one]
\label{cor:proxy}
With $r-1$ missing arcs, warm B has strictly higher overlap and one fewer pivot than warm A, and still $\sim t^2/2$ times its work.
\end{corollary}

\begin{figure}[t]
\centering
\includegraphics[width=\linewidth]{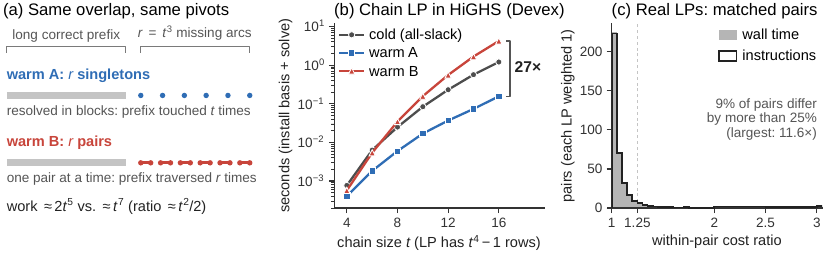}
\caption{Equal overlap and pivot count do not imply equal work. (a) The warm starts of Theorem~\ref{thm:chain}: DSE resolves the singletons of warm A in blocks and the pairs of warm B one at a time. (b) The chain LP in HiGHS (presolve and scaling off, Devex): warm B takes 27 times as long as warm A. (c) Our generated LPs: pairs of starts of one LP with overlap within 0.005 and pivots within 2\% differ in time by more than 25\% in 9\% of cases (class-balanced; the histogram gives each LP equal weight).}
\label{fig:theory}
\end{figure}

The two warm starts differ only in the arrangement of the missing columns (Figure~\ref{fig:theory}a). Most of warm B's work is the solve with the entering column. Each time one of warm B's isolated pairs is attached, this solve spans the long correct prefix of the path, whereas warm A's singletons are first assembled into blocks that touch the prefix rarely. Every dual simplex iteration forms this vector, so the separation is a property of the path and not of a particular implementation. Work on algorithms with predictions bounds the running time of warm-started algorithms by the distance of the prediction from an optimal solution \citep{dinitz2021faster,sakaue2022discrete}; for the simplex method, the theorem rules out such a bound in terms of overlap, even together with the pivot count.

The theorem is an existence result for solve-vector work under one pivot rule. The ratio is unbounded only in a cost model without per-iteration terms linear in the number of rows; an implementation with such a term, or with refactorization after a fixed number of updates, bounds the ratio of total work by a constant (Appendix~\ref{app:discussion}). The theorem motivates our design, and the empirical case rests on the measurements of \S\ref{sec:empirical-gap}.

The separation also appears in a real solver. HiGHS, run with Devex pricing and with scaling and cost perturbation off, takes the same numbers of pivots as in the theorem (Figure~\ref{fig:theory}b). At 65,535 rows it takes 1.2\,s from the cold start, 0.16\,s from warm A and 4.2\,s from warm B, although warm B needs only a sixteenth of the cold start's pivots. Under HiGHS' default edge weights the difference is hidden: a start that is not all-slack and not near optimal first computes exact DSE weights with one BTRAN per row, which is $\Theta(M^2)$ work on this LP. On our generated LPs, this initialization is a median 0.4\% of the solve.

\subsection{On realistic LPs}
\label{sec:empirical-gap}

To test whether the same effects occur on realistic LPs, we solved 381 LPs of \Dtw (25 classes), each from up to 53 starts: 42 from our decoder, ten high-overlap starts and the all-slack basis (Appendix~\ref{app:within}).

\textbf{Overlap and pivots rank starts poorly; instructions rank them almost perfectly.} Within an LP, the pivot count explains 71\% of the variance of log time, adding overlap explains nothing more, and the CPU instruction count explains 98.5\% (Figure~\ref{fig:gap}b; Table~\ref{tab:within}). Choosing a start by either proxy is costly. Among the decoder's 42 starts, the one with the highest overlap takes $1.28\times$ the time of the fastest, and the one with the fewest pivots $1.09\times$. Even when one start has both higher overlap and fewer pivots than another, it is still more than 5\% costlier in 12\% of such pairs, the empirical counterpart of Corollary~\ref{cor:proxy}.

\begin{figure}[t]
\centering
\includegraphics[width=\linewidth]{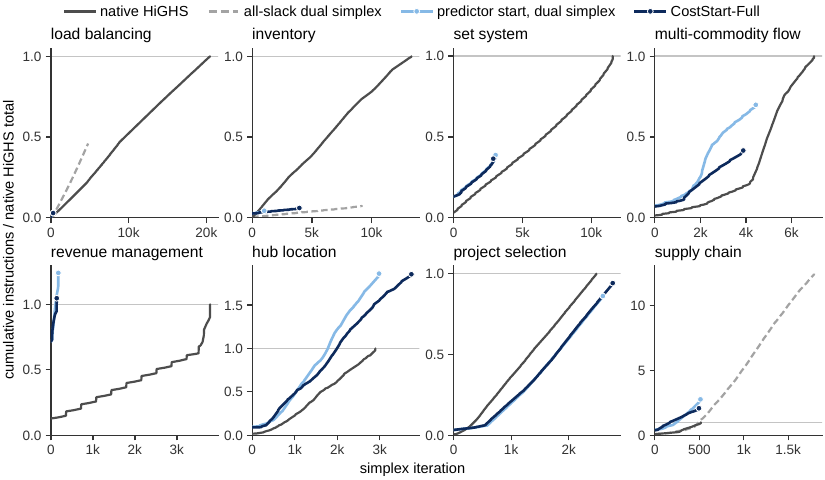}
\caption{Cumulative instructions against iteration, normalized by the native solver's total, on eight of the 20 traced LPs of \Dtw (the LPs at evenly spaced ranks of the native solver's iteration count). The value at iteration 0 is the setup; each curve ends at the solve's total, so a learned curve's end point is its total relative to the native solver. The all-slack dual simplex is drawn only where its path differs from the native solver's.}
\label{fig:itertrace}
\end{figure}

\textbf{The removed iterations are the cheap ones.} An iteration from a learned start costs about 1.3 to 1.5 times as much as one from the all-slack start. To see why, we recorded the instruction count of every iteration on 20 LPs of \Dtw (Figure~\ref{fig:itertrace}; Appendix~\ref{app:within}). Suppose the dual simplex from the all-slack basis takes $N_R$ iterations and the same algorithm from the predictor's start takes $N_L < N_R$. We compare the work per iteration of the first $N_R - N_L$ iterations of the all-slack path with that of its last $N_L$ iterations, excluding setup. On the 17 LPs to which this comparison applies, the early iterations cost $0.44\times$ [0.28, 0.66] as much as the late ones, and the learned path's own iterations cost about as much as the late ones ($0.89\times$ [0.71, 1.17]). A learned start thus skips the cheap early part of the path and pays roughly the late-path price for each of its own iterations. This is why its mean cost per iteration is higher although its total work is lower. Kernel statistics point, as in the theorem, to denser solve vectors, and no static property of the start predicts this cost.

\textbf{Implications for the method.} These findings shape the method in three ways. First, starts of similar closeness to $B^\star$ differ widely in cost, so the training label should be measured cost. Second, the cost landscape is rough: one random basis swap changes the cost by 4\% at the median and by more than 30\% one time in ten, so cost should enter through a choice among a few well-separated options rather than through gradients on every score. Third, a start near $B^\star$ trades fewer iterations for costlier ones, so the options should include ways to make each iteration cheaper.

\section{\ours: supervised prediction, cost-trained decisions}
\label{sec:method}

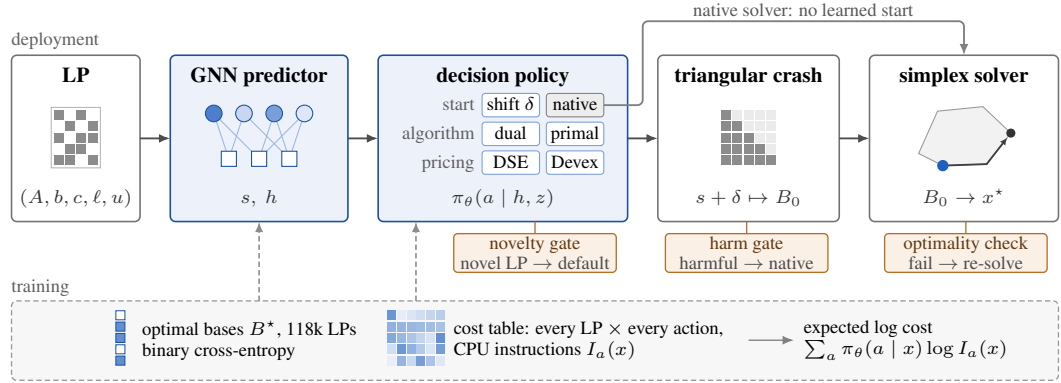
\begin{figure}[t]
\centering
\begin{tikzpicture}[x=1mm, y=1mm,
  font=\fontsize{8}{9.5}\selectfont,
  every node/.style={inner sep=0pt, outer sep=0pt},
  sz7/.style={font=\fontsize{7}{8}\selectfont},
  box/.style={draw=black!55, fill=white, rounded corners=2pt, line width=0.6pt},
  learn/.style={box, draw=csblue!85!black, fill=csblue!8},
  title/.style={anchor=north, font=\fontsize{8}{9.5}\selectfont\bfseries},
  mline/.style={anchor=base, sz7, text=black!80},
  chip/.style={draw=csblue!60, fill=white, rounded corners=1.2pt, line width=0.4pt,
    sz7, rowtext, inner xsep=1.3pt, minimum height=3.1mm, minimum width=7.6mm, anchor=west},
  chipn/.style={chip, draw=black!50, fill=black!7},
  rowlab/.style={sz7, rowtext, text=black!62, anchor=east},
  rowtext/.style={text height=1.75mm, text depth=0.45mm},
  tag/.style={draw=csorange!85!black, fill=csorange!12, rounded corners=1.5pt, line width=0.5pt,
    sz7, align=center, inner ysep=1.1pt, minimum width=21.5mm, anchor=north},
  stem/.style={draw=csorange!85!black, line width=0.5pt},
  seclab/.style={sz7, text=black!62, anchor=base west},
  arr/.style={-{Latex[length=1.5mm, width=1.2mm]}, line width=0.7pt, draw=black!70},
  tarr/.style={-{Latex[length=1.4mm, width=1.1mm]}, line width=0.6pt, draw=black!45,
    dash pattern=on 2pt off 1.2pt},
]
\definecolor{csblue}{HTML}{1F5FBF}
\definecolor{csorange}{HTML}{D9822B}
\draw[box]   (0,-11)     rectangle (17,11);
\draw[learn] (21,-11)    rectangle (44.5,11);
\draw[learn] (48.5,-11)  rectangle (81.5,11);
\draw[box]   (85.5,-11)  rectangle (109.5,11);
\draw[box]   (113.5,-11) rectangle (138.5,11);
\node[title] at (8.5,9.4)   {LP};
\node[title] at (32.75,9.4)  {GNN predictor};
\node[title] at (65,9.4)    {decision policy};
\node[title] at (97.5,9.4)  {triangular crash};
\node[title] at (126,9.4)   {simplex solver};
\node[mline] at (8.5,-8.6)   {$(A,b,c,\ell,u)$};
\node[mline] at (32.75,-8.6) {$s,\ h$};
\node[mline] at (65,-8.6)    {$\pi_\theta(a\mid h,z)$};
\node[mline] at (97.5,-8.6)  {$s+\delta \mapsto B_0$};
\node[mline] at (126,-8.6)   {$B_0 \to x^\star$};
\begin{scope}[shift={(5.65,-3.5)}]
  \foreach \i/\j in {0/0,1/0,3/0,1/1,2/1,0/2,2/2,3/2,1/3,3/3,0/4,2/4}
    \fill[black!45] (\i*1.5,\j*1.5) rectangle ++(1.2,1.2);
  \draw[black!30, line width=0.3pt] (-0.4,-0.4) rectangle (6.1,7.6);
\end{scope}
\begin{scope}[shift={(32.75,0.3)}]
  \foreach \k/\x in {1/-6,2/-2,3/2,4/6} \coordinate (v\k) at (\x,3);
  \foreach \k/\x in {1/-4,2/0,3/4} \coordinate (c\k) at (\x,-3);
  \foreach \a/\b in {1/1,1/2,2/1,2/3,3/2,4/2,4/3,3/3}
    \draw[csblue!45, line width=0.4pt] (v\a) -- (c\b);
  \foreach \k/\f in {1/80,2/25,3/60,4/15}
    \filldraw[draw=csblue!85!black, fill=csblue!\f, line width=0.4pt] (v\k) circle (1.1);
  \foreach \k in {1,2,3}
    \filldraw[draw=csblue!85!black, fill=white, line width=0.4pt] ($(c\k)+(-1,-1)$) rectangle ++(2,2);
\end{scope}
\begin{scope}[shift={(93.9,-3.5)}]
  \foreach \i in {0,...,4} \foreach \j in {0,...,4} {
    \pgfmathparse{int(\i+\j<=4 ? 1 : 0)}
    \ifnum\pgfmathresult=1 \fill[black!45] (\i*1.5,\j*1.5) rectangle ++(1.2,1.2);
    \else \fill[black!8] (\i*1.5,\j*1.5) rectangle ++(1.2,1.2); \fi }
\end{scope}
\begin{scope}[shift={(126,0.2)}]
  \coordinate (p1) at (-6,-2.2); \coordinate (p2) at (-3.5,3.4); \coordinate (p3) at (2.5,4);
  \coordinate (p4) at (6.2,0.5); \coordinate (p5) at (3,-3.6);  \coordinate (p6) at (-2.8,-3.8);
  \draw[black!40, fill=black!5, line width=0.4pt] (p1)--(p2)--(p3)--(p4)--(p5)--(p6)--cycle;
  \draw[black!75, line width=0.7pt, -{Latex[length=1mm,width=0.9mm]}, shorten >=0.75mm] (p6) -- (p5) -- (p4);
  \fill[csblue] (p6) circle (0.75);
  \fill[black!80] (p4) circle (0.6);
\end{scope}
\node[rowlab] at (61.25,4.4)  {start};
\node[rowlab] at (61.25,0.6)  {algorithm};
\node[rowlab] at (61.25,-3.2) {pricing};
\node[chip]  (ch1) at (62.25,4.4) {shift $\delta$};
\node[chipn, right=0.9mm of ch1] (chn) {native};
\node[chip]  (ch2) at (62.25,0.6) {dual};
\node[chip, right=0.9mm of ch2] {primal};
\node[chip]  (ch3) at (62.25,-3.2) {DSE};
\node[chip, right=0.9mm of ch3] {Devex};
\draw[arr] (17,0)    -- (21,0);
\draw[arr] (44.5,0)  -- (48.5,0);
\draw[arr] (81.5,0)  -- (85.5,0);
\draw[arr] (109.5,0) -- (113.5,0);
\draw[arr, draw=black!50, rounded corners=2pt]
  (chn.east) -- (83.5,4.4) -- (83.5,15.5) -- (126,15.5) -- (126,11);
\node[sz7, text=black!65, anchor=south] at (104.75,16.1) {native solver: no learned start};
\node[seclab] at (0,12.4) {deployment};
\foreach \x/\name/\eff in {69.25/novelty gate/novel LP $\to$ default,
                           97.25/harm gate/harmful $\to$ native,
                           126.25/optimality check/fail $\to$ re-solve} {
  \draw[stem] (\x,-11) -- (\x,-12.6);
  \node[tag] at (\x,-12.6) {{\color{csorange!45!black}\name}\\[-0.5pt]{\color{black!70}\eff}};
}
\node[seclab] at (0,-20.3) {training};
\draw[black!40, fill=black!3, rounded corners=2pt, line width=0.5pt, dash pattern=on 2pt off 1.2pt]
  (0,-32) rectangle (138.5,-21.5);
\begin{scope}[shift={(13.8,-30)}]
  \foreach \j/\f in {0/1,1/0,2/1,3/1,4/0} {
    \ifnum\f=1 \filldraw[draw=csblue!85!black, fill=csblue!70, line width=0.3pt] (0,\j*1.5) rectangle ++(1.2,1.2);
    \else \filldraw[draw=csblue!85!black, fill=white, line width=0.3pt] (0,\j*1.5) rectangle ++(1.2,1.2); \fi }
\end{scope}
\node[sz7, align=left, anchor=west] at (17.3,-26.7) {optimal bases $B^\star$, 118k LPs\\binary cross-entropy};
\begin{scope}[shift={(49.7,-30)}]
  \foreach \i in {0,...,5} \foreach \j in {0,...,4} {
    \pgfmathsetmacro{\v}{int(mod(37*\i+11*\j*\j+7*\i*\j,60)+10)}
    \fill[csblue!\v] (\i*1.3,\j*1.5) rectangle ++(1.1,1.2); }
\end{scope}
\node[sz7, align=left, anchor=west] at (58.5,-26.7) {cost table: every LP $\times$ every action,\\CPU instructions $I_a(x)$};
\draw[arr, draw=black!50, line width=0.5pt] (97.5,-26.7) -- (103.5,-26.7);
\node[sz7, align=left, anchor=west] at (104.8,-26.7) {expected log cost\\$\sum_a \pi_\theta(a\mid x)\log I_a(x)$};
\draw[tarr] (32.75,-21.5) -- (32.75,-11);
\draw[tarr] (53.5,-21.5)  -- (53.5,-11);
\end{tikzpicture}%
\caption{\ours at deployment (top) and in training (bottom). Blue boxes are learned and orange tags are guardrails, all charged in the timings; the solver is unmodified. From the predictor's embedding $h$ and score statistics $z$, the policy picks one of up to 20 actions: a shift $\delta$ of the scores $s$ with an algorithm and pricing rule, or the native solver. Training uses optimal bases for the predictor and instruction counts of every action for the policy.}
\label{fig:system}
\end{figure}

An overview of the \ours system is shown in Figure~\ref{fig:system}, which includes a supervised predictor (\S\ref{sec:generator}), a decision policy trained on measured cost (\S\ref{sec:head}), and three guardrails (\S\ref{sec:safety}). Algorithm~\ref{alg:deploy} gives the full per-LP procedure at deployment. The pipeline adds 30 to 50\,ms to a 20,000-variable LP, all of it charged.

\subsection{Supervised predictor and decoder}
\label{sec:generator}

\textbf{Predictor.} We encode the LP as a graph with a node per structural column, per slack column and per constraint, and an edge for every nonzero of $A$ \citep{gasse2019exact}. An 8-layer GATv2 network \citep{brody2022attentive} with 1.8M parameters scores every column of $[A\ I]$ for membership in the optimal basis that the solver returns, trained with binary cross-entropy first on \Trs and then on 3,457 LPs at the 20k size. It follows the supervised approach of \citet{fan2023smart} with a deeper, attention-based network in place of their graph convolution (Appendix~\ref{app:generator}); its forward pass also yields an embedding $h$ of the LP.

\textbf{Decoder.} A triangular crash builds a basis greedily from the scores (Algorithm~\ref{alg:crash}). It repeatedly makes the structural column $j$ basic in the row $i$ with the largest positive gain $(s_j - s_{n+i})/\deg(j)^{1/2}$, where $s_{n+i}$ is the score of row $i$'s slack and $\deg(j)$ the number of entries of column $j$. After each choice, it removes every candidate that would break triangularity. Rows that receive no structural column keep their slack. The selected columns therefore form a triangular submatrix with nonzero diagonal, so the basis is always nonsingular and needs no repair. Adding a constant $\delta$ to every structural score tends to admit more structural columns for $\delta > 0$ and fewer for $\delta < 0$, which gives the decision policy below one of its levers.

\begin{algorithm}[t]
\caption{Triangular crash with a structural shift}
\label{alg:crash}
\small
\begin{algorithmic}[1]
\Require scores $s \in \mathbb{R}^{n+m}$ (structural columns $1..n$, slacks $n+1..n+m$), shift $\delta$, exponent $p$
\State $E \gets \{(j,i) : |a_{ij}| \ge 10^{-8}\max_k |a_{kj}|\}$; \quad $\deg(j) \gets |\{i : (j,i) \in E\}|$
\State $g_{ji} \gets (s_j + \delta - s_{n+i})/\deg(j)^{p}$ for $(j,i) \in E$; \quad $\mathcal{B} \gets \emptyset$
\While{some $(j,i) \in E$ has $g_{ji} > 0$}
  \State take $(j,i) = \arg\max g_{ji}$ (ties: smaller slack score, then smaller index)
  \State $\mathcal{B} \gets \mathcal{B} \cup \{(j, i)\}$ \Comment{column $j$ becomes basic in row $i$}
  \State remove from $E$ every edge of column $j$ and every edge whose row lies in $\mathrm{supp}(a_j)$
\EndWhile
\State every row $i$ without a structural column keeps its slack $n+i$ basic
\State set nonbasic bound statuses from the bound types and the signs of the scores
\State \Return basis $\mathcal{B}$ with statuses
\end{algorithmic}
\end{algorithm}

\begin{algorithm}[t]
\caption{\ours at deployment, for one LP}
\label{alg:deploy}
\small
\begin{algorithmic}[1]
\Require LP $x = (A, b, c, \ell, u)$, stored after presolve; action set $\mathcal{A}$; predictor, policy ensemble, novelty and harm gates
\State $s, h \gets$ predictor$(x)$ \Comment{column scores and LP embedding, one forward pass}
\State $z \gets$ statistics of $s$ and of the size of $x$
\If{nearest-neighbor distance of $[h, z]$ to the training LPs $>$ threshold} \Comment{novelty gate}
  \State $a \gets$ (shift 0, dual simplex, the solver's default pricing)
\Else
  \State $a \gets \arg\max_{a' \in \mathcal{A}} \frac{1}{E}\sum_{e=1}^{E} \log \pi_{\theta_e}(a' \mid h, z)$ \Comment{ensemble of $E$ policies}
\EndIf
\If{$a$ is the native action}
  \State solve $x$ with the native solver; \Return
\EndIf
\State $B_0 \gets$ triangular crash on $s + \delta_a \mathbf{1}_{\text{struct}}$ (Algorithm~\ref{alg:crash})
\If{harm gate is used \textbf{and} predicted saving of $B_0$ over the native solver $< 12\%$}
  \State solve $x$ with the native solver; \Return
\EndIf
\State \texttt{setBasis}$(B_0)$; run the solver with the algorithm and pricing rule of $a$
\If{sum of dual infeasibilities $> 10^{-6}(1 + |f|)$} \Comment{optimality check}
  \State continue from the returned basis with primal and dual feasibility tolerances $10^{-9}$
  \State if not optimal, solve $x$ afresh with the native solver
\EndIf
\end{algorithmic}
\end{algorithm}

\subsection{Decision policy}
\label{sec:head}

A predictor of $B^\star$ is trained and judged only by closeness to $B^\star$; nothing in its objective sees what each iteration costs. The decision policy chooses, per LP, among a few well-separated actions that also change the setup and the work per iteration, and it is trained on the measured cost of every action. We study three nested action sets of 9, 11 and 20 actions.

\textbf{Shifts ($\mathcal{A}_1$).} Nine shifts $\delta \in \{0, \pm 0.5, \pm 1, \pm 2, \pm 3\}$ of the structural scores, each decoded by the crash and solved with the dual simplex. The policy over $\mathcal{A}_1$, with the harm gate below as a fallback to the native solver, is \oursS.

\textbf{Algorithm ($\mathcal{A}_2$).} The nine shifts plus two strategy actions: the native solver, for LPs on which no learned start pays off, and the unshifted start continued with primal simplex, for starts that are nearly primal feasible. The policy over $\mathcal{A}_2$ is \oursA.

\textbf{Pricing ($\mathcal{A}_3$).} $\mathcal{A}_2$ plus the nine shifts with Devex pricing \citep{harris1973pivot} forced from the first iteration. By default, HiGHS uses DSE edge weights and switches to Devex on its own mainly for a near-optimal supplied basis or after at least $0.1(m+n)$ iterations of costly DSE solves. Devex drops DSE's extra solve $\tau$ and its exact weights. It therefore trades more iterations for cheaper ones, and which effect dominates varies across classes. Of 13 alternative solver options we probed, Devex gave the largest gain (Appendix~\ref{app:periter}). The policy over $\mathcal{A}_3$ is \oursF.

\textbf{Cost signal.} The label of action $a$ on LP $x$ is $I_a(x)$, the CPU instructions the solver executes to install the start and solve, read from a hardware counter. Unlike wall time, this count is unaffected by machine load: repeated solves of the same LP differ by a median of only 0.005\%. Building the table of costs for $\mathcal{A}_3$ takes about 16 times the work of one cold solve per training LP.

\textbf{Policy and loss.} The policy $\pi_\theta(a \mid x) = \mathrm{softmax}(f_\theta(h(x), z(x)))$ is a two-layer MLP on $h$ and 12 statistics $z$ of the scores and the LP's size, trained on the exact expected log cost relative to the unshifted start,
\begin{equation}
\label{eq:loss}
\mathcal{L}(\theta) \;=\; \frac{1}{|\mathcal{X}|}\sum_{x \in \mathcal{X}} \sum_{a \in \mathcal{A}} \pi_\theta(a \mid x)\, \mathrm{clip}\big(\log I_a(x) - \log I_0(x),\, -3,\, 3\big),
\end{equation}
where $I_0$ is the cost of the unshifted start. Up to the clip, the exponential of this loss is the geometric-mean instruction ratio, the training counterpart of the ratio we report (Appendix~\ref{app:loss}). Because every action's cost is known for every training LP, the loss is computed exactly and has no sampling variance. Because the label is total cost, the policy prefers an action with more iterations whenever each iteration becomes cheap enough. At deployment we play the argmax of an ensemble of policies; the policy runs in under 2\,ms.

\textbf{Why not train the predictor on the cost?}
\label{sec:negatives}
Every attempt we made to train the predictor itself on the solver's cost, including policy gradients through a stochastic crash, perturbed optimizers and end-to-end training, failed to beat continued supervision at matched budget (Appendix~\ref{app:negatives}). The per-LP cost gradients point in unrelated directions: the cosine between the summed gradients of two disjoint halves of a batch is near 0, against 0.51 for the supervised loss. We therefore freeze the predictor and let cost enter only through the decision policy.

\subsection{Guardrails}
\label{sec:safety}

The policy is trained only on LPs of the training classes, and the decoded start can be costly even within them. Three guardrails limit the resulting losses; all of them are included in the charged time.

\textbf{Novelty gate.} Outside the training classes, the policy confidently chooses harmful actions. If the input $[h, z]$ is farther from every training LP than the 0.99 quantile of in-distribution nearest-neighbor distances, the policy plays the default action (unshifted start, dual simplex, the solver's default edge weights).

\textbf{Harm gate.} A gradient-boosted regressor predicts, from cheap features of the LP and the decoded start and without class labels, whether the learned start will save at least 12\% over the native solver; if not, the native solver runs. \oursA and \oursF can choose the native solver themselves, so for them the gate is needed mainly where the novelty gate takes the decision away from the policy, but it is still evaluated and charged on every LP.

\textbf{Optimality check.} One badly scaled class has optimal duals of order $10^{-7}$, the same size as the solver's dual feasibility tolerance. From a learned start, the solver can end at a basis with many wrong-signed reduced costs, each below that tolerance, and accept it as optimal although its objective is up to 20\% too high. The solver's sum of dual infeasibilities, unlike its count, includes these small values. If the sum exceeds $10^{-6}(1 + |f|)$ for objective $f$, we continue from the returned basis with feasibility tolerances of $10^{-9}$, and fall back to the native solver if that does not reach optimality. On development data, this rule flagged all 12 wrong solves and none of 2,358 correct ones (Appendix~\ref{app:check}).

\section{Experimental setup}
\label{sec:setup}

\textbf{Data.} All LPs come from SyntheticLPs, a library of parametric generators for 45 operations-research problem classes (Appendix~\ref{app:Syn_generator}). Each class has one or more formulations, and each formulation can be generated in a primal or a dual form. The dual form has a transposed constraint matrix and looks very different to the solver, so we count the two forms as separate variants. The 26 classes we train on have 134 variants: 67 formulations, each in a primal and a dual form. Formulations with integer or binary variables enter through their LP relaxations, so every instance is an LP. Every LP is presolved once and stored, and every method solves the stored LP. Set names give the typical number of variables; these are the sizes requested from the generator, and after presolve the median is 18,000 variables in \Trtw and 88,000 in \Tei. Table~\ref{tab:data} lists all sets.
\begin{itemize}
\item \textbf{Training.} The predictor is trained on \Trs (114,874 LPs, 50 to 2,000 variables) and then on 3,457 LPs at the 20k size. The policies and the harm gate are trained on \Trtw (9,511 LPs, 15,000 to 25,000 variables).
\item \textbf{Development.} \Dtw (15,000 to 25,000 variables) and \Dfif (40,000 to 60,000) are used for model selection, ablations and the design of the optimality check.
\end{itemize}
Test LPs are generated from fresh seeds, and an LP is admitted only if no training or development LP shares its content hash, generator call or seed (Appendix~\ref{app:data}).

\textbf{Test sets.}
\begin{itemize}
\item \Ttw: 1,356 LPs with 15,000 to 25,000 variables from 25 training classes. This is the main test set.
\item \Tfif (382 LPs, 40,000 to 60,000 variables) and \Tei (81 LPs, 60,000 to 100,000): larger instances of the same classes.
\item \Tva: the 296 LPs of \Ttw from ten formulations and six dual forms that a second, separately trained pipeline never sees.
\item \Tnc: 234 LPs from six classes absent from training.
\item \Tad: LPs from five of these unseen classes, with 125 for adaptation and 616 for evaluation.
\end{itemize}
An LP was admitted to a test set only if the all-slack dual simplex needed at least 2,000 iterations or 1\,s on it in a separate run at generation time. LPs are grouped by that run's time into windows W1 (up to 10\,s), W2 (10 to 120\,s) and W3 (over 120\,s), with time limits of 120, 600 and 3,600\,s.

\textbf{Protocol.} Every method solves each LP in a fresh solver instance on one pinned core; \Ttw was timed twice and averaged. Charged time includes features, inference, the policy, decoding, installing the basis, the solve, gates, checks and fallbacks. The predictor runs on one GPU and is charged at its wall time (on CPU cores alone it would make \oursS about 30\% slower on solves under 10\,s; Appendix~\ref{app:cost}). Failures (time limits, errors, objectives off by more than $10^{-5}$ relative to a reference) count at twice the time limit. We report iterations relative to the all-slack dual simplex, which runs the same algorithm as the warm starts, and time relative to the native solver. The main statistic is the class-balanced geometric mean of per-LP ratios (log ratios averaged within each class, then across classes), with 95\% bootstrap intervals over classes and then generator instances.

\textbf{Methods.} Baselines are the native solver, the all-slack dual simplex, the stock solver (all default options, presolve on, applied to the stored presolved LP) and a classical crash (inverse-degree scores through our decoder). We reimplemented \citet{fan2023smart} and \citet{lu2025learning} (0.2M and 1.2M parameters) and trained them on our data, both at their original training budget and at ours. Their checkpoints are chosen by their own loss (ours by solve time on development data), and their scores are decoded with our crash, which suits them better than their own decoders on our data (Appendix~\ref{app:literature}).

\section{Results}
\label{sec:results}

\begin{table}[t]
\centering\footnotesize
\caption{Main results on unseen LPs of the training classes: every learned start removes half or more of the iterations, but only our systems save substantial time. Class-balanced geometric means (below 1 is better) of iterations relative to the all-slack dual simplex and of charged time relative to the native solver, with 95\% intervals and classes won/lost. The class-blind gate is refit without the test LP's class; the per-class lookup plays each class's cheapest training action (Appendix~\ref{app:lookup}). On \Tfif and \Tei, methods run without the optimality check, and the larger policies are reported in Table~\ref{tab:general}a. More rows: Appendix~\ref{app:results}.}
\label{tab:main}
\setlength{\tabcolsep}{2.8pt}
\begin{tabular}{@{}l ccc ccc c@{}}
\toprule
& \multicolumn{3}{c}{\Ttw (1,356 LPs)} & \multicolumn{3}{c}{\Tfif (382 LPs)} & \Tei (81)\\
\cmidrule(lr){2-4}\cmidrule(lr){5-7}\cmidrule(l){8-8}
Method & iter. & time [95\% CI] & W/L & iter. & time [95\% CI] & W/L & time\\
\midrule
all-slack dual simplex & 1.00 & 1.07 \ci{0.87}{1.30} & 12/13 & 1.00 & 0.99 \ci{0.76}{1.30} & 12/13 & 0.96\\
stock solver, presolve on & 0.85 & 1.10 \ci{1.06}{1.16} & 3/22 & -- & -- & -- & --\\
classical crash & 1.22 & 1.76 \ci{1.19}{2.60} & 9/16 & -- & -- & -- & --\\
Fan et al.\ & 0.49 & 1.10 \ci{0.89}{1.35} & 15/10 & 0.50 & 1.02 \ci{0.78}{1.37} & 15/10 & 0.96\\
Lu \& Yan & 0.47 & 0.94 \ci{0.70}{1.22} & 14/11 & 0.50 & 0.85 \ci{0.56}{1.23} & 19/6 & 0.72\\
\midrule
predictor & 0.32 & 0.68 \ci{0.50}{0.91} & 20/5 & 0.33 & 0.61 \ci{0.39}{0.89} & 19/6 & 0.60\\
\oursS & 0.30 & 0.61 \ci{0.46}{0.76} & 23/2 & 0.32 & 0.53 \ci{0.34}{0.72} & 23/2 & 0.52\\
\quad with class-blind harm gate & 0.32 & 0.68 \ci{0.52}{0.85} & 18/7 & -- & 0.58 \ci{0.40}{0.79} & 18/7 & 0.58\\
\oursA & 0.25 & 0.51 \ci{0.38}{0.65} & 23/2 & \multicolumn{4}{c}{see Table~\ref{tab:general}}\\
\oursF & 0.28 & 0.49 \ci{0.37}{0.63} & 24/1 & \multicolumn{4}{c}{see Table~\ref{tab:general}}\\
\quad with harm gate (deployed) & 0.29 & 0.52 \ci{0.39}{0.65} & 23/2 & \multicolumn{4}{c}{see Table~\ref{tab:general}}\\
\midrule
per-class lookup (uses class labels) & -- & 0.54 \ci{0.41}{0.67} & -- & \multicolumn{4}{c}{--}\\
\bottomrule
\end{tabular}
\end{table}

\subsection{Unseen LPs of the training classes}
\label{sec:headline}

Every learned start removes half or more of the iterations (Table~\ref{tab:main}), but the published methods end up within 10\% of the native solver in time on \Ttw. Relative to the native solver, our predictor alone takes $0.68\times$ the time, \oursS $0.61\times$ and \oursF $0.49\times$; \oursF is faster than the native solver on 24 of the 25 classes. On the larger LPs of \Tfif and \Tei, \oursS still takes about half the time. The 95\% intervals of \oursS, and of the larger policies with the harm gate, exclude parity on all three test sets. We regard \oursF with both gates and the optimality check as the deployed configuration. At the training size the harm gate costs 6\%, but larger LPs need it (\S\ref{sec:generalization}).

\textbf{Where the gain comes from.} Most of our advantage over the published methods at the level of the start comes from the predictor; more training, another seed or the published methods' own decoders do not close the gap (Appendix~\ref{app:literature}). The decision policies then add to the predictor in different ways:
\begin{itemize}
\item \oursF lowers the predictor's time by a further 28\%, which makes it 1.9 to 2.2 times as fast as the published methods on \Ttw. On the larger test sets, \oursS is 1.4 to 1.9 times as fast as them.
\item \oursS lowers the predictor's time by 10\%, but that gain comes from its harm gate rather than its shift policy. The gate and the shift policy act on the same LPs, so with the gate present, the shift policy adds nothing on the test sets (Appendix~\ref{app:mechanism}). Part of the gate's value comes from recognizing the training classes: refit without the test LP's class, it loses most of its value (Table~\ref{tab:main}).
\end{itemize}

\textbf{A per-class lookup control.} A lookup that plays, for each class, the action with the lowest training cost is a strong control, although it needs class labels that our methods never see (Appendix~\ref{app:lookup}). It beats \oursS, which only chooses between the learned start and the native solver; a two-action lookup restricted to that choice ties \oursS. \oursF is 9\% faster than the lookup overall, which is not significant, and takes $0.76\times$ [0.64, 0.87] its time on solves of 10 to 120\,s.

\subsection{Fewer iterations and cheaper iterations}
\label{sec:periter}

\begin{figure}[t]
\centering
\includegraphics[width=.7\linewidth]{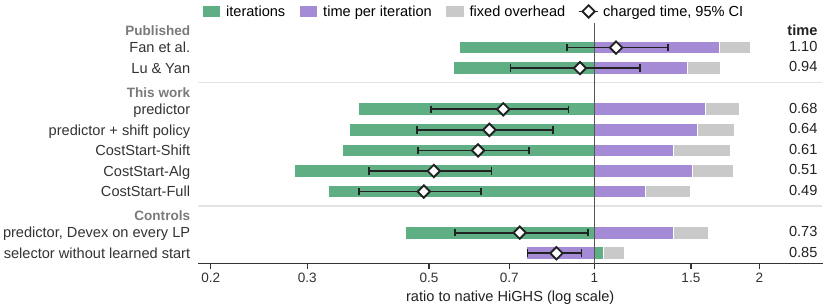}
\caption{Every learned start trades fewer iterations for costlier ones. Charged time on \Ttw relative to the native solver, split by Eq.~(\ref{eq:decomp}) into iterations, time per iteration and fixed overhead (class-balanced geometric means). On the log axis, segments left of 1 lower the time, segments right of 1 raise it, and the diamond (95\% interval) is their sum. ``Predictor + shift policy'' is \oursS without the harm gate.}
\label{fig:decomp}
\end{figure}

Figure~\ref{fig:decomp} applies the decomposition (\ref{eq:decomp}) to \Ttw. Every learned start cuts the iterations of the native solver by at least 40\%. For the published methods and our predictor, however, each iteration becomes 1.5 to 1.7 times as costly, and the published methods remove too few iterations to compensate. Our two larger policies reach similar times by different routes. \oursA gains mainly through fewer iterations, and \oursF mainly through cheaper ones: \oursF takes 15\% more iterations than \oursA, each 18\% cheaper, a net difference that is not significant (Appendix~\ref{app:periter}). The pricing choice has to be made per LP and together with the start. Forcing Devex on every LP, or training a selector with the same loss over the solver's own starts and pricing rules without the learned start, is slower than the predictor alone.

\textbf{Which cost to train on.} The training label matters most when the actions change the cost per iteration. We therefore trained the \oursF policy with the same recipe on log iteration counts instead of log instructions, and scored both policies offline on held-out training LPs (Table~\ref{tab:label-full}). On solves under 10\,s, the iteration-trained policy removes 15\% more pivots than the instruction-trained one but needs $1.40\times$ [1.18, 1.73] its instructions; on longer solves the two tie. The difference lies in the actions chosen: the iteration-trained policy sends far more LPs to primal simplex and far fewer to Devex, and as a result its iterations cost 1.65 times as much. Shift policies show the same pattern in charged time (Figure~\ref{fig:mech}a; Appendix~\ref{app:objective}), and on training LPs the shift closest to the known $B^\star$ gains 4\%, against 7\% for a cost-trained policy. Three retrained pipelines reproduce \oursS on \Ttw to within about 1\% (Appendix~\ref{app:seeds}).

\begin{table}[t]
\centering\footnotesize
\caption{Training label of the \oursF policy, offline on held-out training LPs: instructions of the chosen actions relative to the predictor's unshifted start, split into iterations and work per iteration, and the share of LPs routed to Devex and to primal simplex. The selection split was used to design the policy and is not held out.}
\label{tab:label-full}
\scriptsize
\setlength{\tabcolsep}{3.5pt}
\begin{tabular}{@{}llccccc@{}}
\toprule
LPs & Label & instructions [95\% CI] & iterations & per iteration & Devex & primal\\
\midrule
held-out, W1 (1,531) & instructions & 0.669 \ci{0.589}{0.751} & 0.892 & 0.749 & 51\% & 8\%\\
 & iterations & 0.937 \ci{0.817}{1.103} & 0.759 & 1.234 & 12\% & 36\%\\
held-out, W2 (316) & instructions & 0.668 \ci{0.526}{0.826} & 0.729 & 0.916 & 29\% & 20\%\\
 & iterations & 0.679 \ci{0.489}{0.998} & 0.601 & 1.129 & 5\% & 51\%\\
selection split of \Dtw (275 of 283) & instructions & 0.682 \ci{0.570}{0.797} & 0.868 & 0.786 & 60\% & 9\%\\
 & iterations & 0.960 \ci{0.798}{1.170} & 0.723 & 1.327 & 15\% & 46\%\\
\midrule
\multicolumn{7}{@{}p{0.95\linewidth}@{}}{Iteration label relative to instruction label, paired (instructions; iterations; per iteration): held-out W1 1.402 \ci{1.183}{1.727}; 0.851 \ci{0.769}{0.933}; 1.647 \ci{1.364}{2.040}. Held-out W2 1.016 \ci{0.767}{1.335}; 0.824 \ci{0.684}{0.950}; 1.233 \ci{1.002}{1.649}. Selection split 1.408 \ci{1.152}{1.771}; 0.834 \ci{0.742}{0.917}; 1.689 \ci{1.352}{2.152}.}\\
\bottomrule
\end{tabular}
\end{table}

\subsection{Generalization and deployment}
\label{sec:generalization}

\begin{table}[t]
\centering\scriptsize
\caption{\ours keeps most of its gain on larger LPs and held-out variants; new classes need adaptation. Charged time relative to the native solver (class-balanced). (a) Solves under 120\,s, with the optimality check inline (failures in parentheses); Table~\ref{tab:main} also includes longer solves. (b) With 95\% intervals; the last column is relative to the reference named in each block.}
\label{tab:general}
\setlength{\tabcolsep}{2pt}
\begin{minipage}[t]{0.40\linewidth}\vspace{0pt}
\centering
\begin{tabular}{@{}lcc@{}}
\toprule
(a) Larger LPs & \Tfif & \Tei\\
 & (309) & (63)\\
\midrule
all-slack dual simplex & 0.94 & 0.90\\
predictor & 0.65 (1) & 0.83 (1)\\
\quad Devex on every LP & 0.66 & 0.78 (1)\\
\oursS & 0.58 & 0.67\\
\oursA & 0.60 (2) & 0.83 (1)\\
\quad + harm gate & 0.55 & 0.66\\
\oursF & 0.56 (1) & 0.81 (1)\\
\quad + harm gate & 0.54 & 0.65\\
\bottomrule
\end{tabular}
\end{minipage}\hfill
\begin{minipage}[t]{0.58\linewidth}\vspace{0pt}
\centering
\begin{tabular}{@{}lll@{}}
\toprule
(b) Unseen types & vs.\ native & vs.\ reference\\
\midrule
\multicolumn{3}{@{}l}{\Tva, 10 formulations and 6 dual forms (ref.: deployed predictor)}\\
\quad predictor trained without these & 0.71 \ci{0.33}{1.41} & 1.33 \ci{1.13}{1.64}\\
\quad + own shift policy, deployed gate & 0.56 \ci{0.30}{0.86} & 1.05 \ci{0.67}{1.50}\\
\multicolumn{3}{@{}l}{\Tnc, 6 unseen classes (ref.: predictor)}\\
\quad predictor & 1.42 \ci{0.95}{2.30} & 1.00\\
\quad \oursS & 1.21 \ci{0.97}{1.62} & 0.85 \ci{0.62}{1.08}\\
\quad classical crash & 0.67 \ci{0.30}{1.18} & --\\
\multicolumn{3}{@{}l}{\Tad, 5 unseen classes (ref.: predictor, where the policy abstains)}\\
\quad adapted, 25 labeled LPs per class & 1.04 \ci{0.73}{1.48} & 0.91 \ci{0.84}{0.98}\\
\quad online bandit, no labels & 0.93 \ci{0.79}{1.09} & 0.82 \ci{0.64}{0.97}\\
\bottomrule
\end{tabular}
\end{minipage}
\end{table}

\textbf{Larger LPs.} Training LPs were generated at up to 25,000 variables, yet the predictor keeps most of its advantage at up to 100,000 (Table~\ref{tab:main}). The decision policy transfers only in part. The novelty gate sends most larger LPs to the default action, where the policy cannot choose the native solver, so the harm gate becomes necessary. Table~\ref{tab:general}a reports a separate run over the solves under 120\,s, where savings are smaller, with the optimality check inline; the check removes the wrong optima that \oursS returns on \Tfif without it. With both gates and the check, \oursF and \oursA are the fastest methods and have no failures.

\textbf{New variants.} On \Tva, a predictor retrained without these ten formulations and six dual forms is 33\% slower than the deployed one. Its own shift policy adds nothing there, and the deployed harm gate, which did see these formulations and dual forms in training, recovers most of the loss. On development data, held-out dual forms keep most of the predictor's gain and held-out formulations about a quarter (Appendix~\ref{app:hc}).

\textbf{New classes.} On six classes that no model saw, every learned start of ours is slower than the native solver, while the classical crash is faster (Table~\ref{tab:general}b). The novelty gate flags every such LP, but its fallback, the unshifted start, was chosen by leave-one-class-out validation on the training classes and proves wrong for new classes. The supported scope is therefore new instances, sizes and variants of the training classes.

For new classes, we tested two ways to adapt from the solver's own cost (Appendix~\ref{app:fewshot}). On \Tad, the unadapted policy abstains and falls back to the predictor's start.
\begin{itemize}
\item \textbf{With labels.} Retraining the policy with 25 labeled LPs per new class brings its time to $0.91\times$ [0.84, 0.98] that of the unadapted policy. This does not count the cost of labeling, about ten times the plain solves of the labeled LPs, so it pays off only for streams of thousands of LPs per class.
\item \textbf{Without labels.} A contextual bandit over shifts and solver strategies, learning from the cost of its own solves, ends at $0.82\times$ [0.64, 0.97] the predictor's time and $0.93\times$ that of the native solver over streams of 64 to 149 LPs per class. The same bandit does not pay off on the shorter streams of \Tnc (26 to 91 LPs per class).
\end{itemize}

\subsection{The public benchmark of Lu \& Yan}
\label{sec:luyan}

Finally, we evaluate on the four public datasets of \citet{lu2025learning} under their protocol: one model per dataset, a 7:3 split (ours, since theirs is unpublished) and 80 test LPs in total (Table~\ref{tab:luyanbench}; Appendix~\ref{app:luyanbench}). We fine-tune our predictor per dataset. For each dataset, cross-validation on the training LPs then chooses between a trained decision policy and the best single action, over an action set that also includes their top-$m$ decoder. The result takes $0.72\times$ [0.68, 0.75] (bootstrap over LPs) the time of the dual simplex, against $1.03\times$ for our reimplementation of Lu \& Yan and $0.92\times$ for Fan et al.; Lu \& Yan report $0.73\times$ on their own split. With only 19 to 70 training LPs per dataset, most of the gain comes from choosing the best single action per dataset ($0.73\times$) rather than from per-LP decisions.

\begin{table}[t]
\centering\scriptsize
\caption{Lu \& Yan's public benchmark under their protocol (one model per dataset, 7:3 split, presolved LPs): geometric-mean ratios over test LPs against the dual simplex with presolve off (iterations / fully charged time), with the geometric mean over the four datasets and its 95\% interval. Test LPs per dataset: Mirp1 9, Mirp2 18, Anonymous 23, Load balancing 30; native solve times 9.6, 30, 1.2 and 3.0\,s. All rows come from one paired panel except the two reported ones.}
\label{tab:luyanbench}
\setlength{\tabcolsep}{3.5pt}
\begin{tabular}{@{}lcccc c@{}}
\toprule
Method (decoder) & Mirp1 & Mirp2 & Anonymous & Load bal. & Geomean [95\%]\\
\midrule
Fan et al.\ reimplemented (top-$m$) & 0.58 / 0.72 & 1.01 / 1.14 & 0.92 / 1.62 & 0.61 / 0.69 & 0.76 / 0.98 \ci{0.93}{1.04}\\
Fan et al.\ reimplemented (top-$m$ repair) & 0.55 / 0.69 & 0.78 / 0.92 & 0.87 / 1.69 & 0.56 / 0.68 & 0.68 / 0.92 \ci{0.86}{0.99}\\
Lu \& Yan reimplemented (top-$m$) & 0.70 / 1.10 & 0.92 / 1.07 & 0.71 / 1.37 & 0.60 / 0.68 & 0.72 / 1.02 \ci{0.95}{1.10}\\
Lu \& Yan reimplemented (top-$m$ repair) & 0.71 / 1.12 & 0.80 / 0.99 & 0.70 / 1.48 & 0.58 / 0.70 & 0.69 / 1.03 \ci{0.96}{1.12}\\
our predictor, fine-tuned (top-$m$) & 0.43 / 0.64 & 1.00 / 1.07 & 0.55 / 1.05 & 0.59 / 0.67 & 0.61 / 0.83 \ci{0.78}{0.88}\\
our predictor, fine-tuned (top-$m$ repair) & 0.35 / 0.57 & 0.99 / 1.06 & 0.56 / 1.27 & 0.55 / 0.69 & 0.57 / 0.85 \ci{0.80}{0.91}\\
best single action per dataset & 0.39 / 0.54 & 0.96 / 1.00 & 0.55 / 1.04 & 0.51 / 0.51 & 0.57 / 0.73 \ci{0.67}{0.80}\\
\ours decision policy, per dataset & 0.35 / 0.53 & 0.98 / 1.02 & 0.65 / 0.94 & 0.51 / 0.51 & \textbf{0.58 / 0.72} \ci{0.68}{0.75}\\
\midrule
Lu \& Yan, as reported (their decoder) & \multicolumn{4}{c}{--} & 0.60 / 0.73\\
Fan et al., as reported (their decoder) & \multicolumn{4}{c}{--} & 0.74 / 0.89\\
\bottomrule
\end{tabular}
\end{table}

\section{Limitations and conclusion}
\label{sec:conclusion}

\textbf{Limitations.}
\begin{itemize}
\item \textbf{New classes.} The learned components do not transfer to classes absent from training. Adaptation recovers part of the loss but needs labels or a stream of about a hundred LPs per class.
\item \textbf{Data and solver.} All main results use one generator library and one solver, and test LPs were admitted by the dual simplex's solve time. The generated LPs have a narrower range of coefficients and less degeneracy than the Netlib collection (Appendix~\ref{app:Syn_generator}).
\item \textbf{Presolve.} Our LPs are presolved before every method runs. On the raw files of the public benchmark, the solver's own presolve is faster than every warm start, ours included, so learned starts help only where presolve has already been applied or is not used (Appendix~\ref{app:luyanbench}).
\item \textbf{Inference hardware.} Our timings charge a GPU forward pass. On CPU cores alone, the slower forward pass makes \oursS 29 to 30\% slower on solves under 10\,s (Appendix~\ref{app:cost}).
\item \textbf{Scope of the theorem.} The theorem concerns solve-vector work under one pivot rule, and a per-iteration cost linear in the number of rows bounds its ratio (Appendix~\ref{app:discussion}).
\end{itemize}

\textbf{Conclusion.} Learned bases save a smaller share of time than of iterations, because the iterations they remove are cheap and closeness to the optimum does not rank starts by cost. \ours keeps supervised prediction for the basis, trains the decisions about its use on measured cost, and roughly halves the time of the native solver on unseen LPs of the training classes, including larger ones.

\clearpage
\bibliography{cleaned_references}
\bibliographystyle{iclr2027_conference}

\clearpage
\appendix
\section{Proof of Theorem~\ref{thm:chain} and Corollary~\ref{cor:proxy}}
\label{app:proofs}

We prove Theorem~\ref{thm:chain} and Corollary~\ref{cor:proxy}, state the cost model under which solve-vector work bounds total work, and report how HiGHS behaves on the chain LP. Throughout, $M = t^4$ is the number of path vertices $1, \ldots, M$ as in the theorem, $R = M - 1$ is the number of equality rows (the main text's $m$), $q = t^2$ and $r = t^3$; all arithmetic is exact. The variables $v_j$ are explicit surplus columns $-e_j$ that play the role of slacks; nothing relies on the solver's built-in slack (logical) columns.

\textbf{Solve-vector work.} A dual simplex iteration with leaving row $p$ forms the BTRAN result $\rho_p = B^{-\top}e_p$, the FTRAN result $\alpha = B^{-1}a_{\mathrm{in}}$ of the entering column $a_{\mathrm{in}}$ and, under DSE pricing, the second FTRAN result $\tau = B^{-1}\rho_p$. The solve-vector work $K$ of a path is the total number of nonzeros of these three vectors over its pivots. The vectors are determined by the basis and the pivot, not by the LU factors, so $K$ is intrinsic to the sequence of bases, and it is a lower bound on the work of any implementation that forms them. It omits PRICE (the pivotal row $\rho_p^\top A_N$), the choice of the leaving row (CHUZR), the choice of the entering column with the bound-flipping ratio test (CHUZC and BFRT), factor updates and reinversion. On the chain LP every row has at most three nonzeros, and the arithmetic model below bounds these omitted steps by $O(K) + O(\log M)$ per pivot.

\subsection{The chain LP as a minimum-cost-flow problem}

Index nodes $0, \ldots, R$ and use the convention $+1$ at an arc's head and $-1$ at its tail, deleting the row of node 0. Then $u_1$ is the arc $0 \to 1$, $u_j$ ($2 \le j \le R$) the arc $j-1 \to j$, $u_M$ the arc $R \to 0$, and $v_j$ the arc $j \to 0$, with unit cost on every $u$ arc and zero cost on every $v$ arc. The LP is therefore an uncapacitated minimum-cost transshipment problem on a path with return arcs to a hub, with $R$ rows, $2M - 1$ columns and $3R$ nonzeros, and its constraint matrix is a reduced node-arc incidence matrix, hence totally unimodular. The separation of Theorem~\ref{thm:chain} is thus $\Theta(\sqrt{R})$ in the number of rows. Two cautions apply. Total unimodularity alone does not make every selection of columns a basis or free of fill; the selections used here are rooted spanning trees, and in the stated slot order their basis matrices are upper bidiagonal, which gives both properties. The right-hand side is rational, so basic solutions need not be integral; multiplying $b$ by $100M$ gives integer supplies with identical DSE choices and solve-vector supports. The construction is an existence result within this class and makes no claim about typical network-flow instances.

\subsection{Basis family, dual feasibility and the unique optimum}

Set $e_0 = e_M = 0$ in $\mathbb{R}^R$. The columns of $u_j$ and $v_j$ are $e_j - e_{j-1}$ and $-e_j$. For $S \subseteq \{1, \ldots, R\}$, put $u_j$ in slot $j$ if $j \in S$ and $v_j$ otherwise; $u_M$ stays nonbasic. The basis matrix is $B(S) = U(S)\Sigma(S)$ with $U(S) = I - \sum_{j \in S, j \ge 2} e_{j-1}e_j^\top$ and $\Sigma_{jj} = +1$ if $j \in S$, $-1$ otherwise. It is upper bidiagonal with unit-magnitude diagonal and at most $2R - 1$ nonzeros, so $L = I$, $U = B(S)$ is an exact no-fill factorization and $|\det B(S)| = 1$.

Selected edges join consecutive vertices into components. If $[h, v]$ is a component, back-substitution gives $u_j = \sum_{s=j}^{v-1} b_s > 0$ for $j \in S$ in $[h, v]$ and $v_j = -\sum_{s=j}^{v-1} b_s < 0$ for a missing boundary $j$ whose right component is $[j+1, v]$. Hence exactly the missing-boundary basics are primal infeasible.

Let $\lambda$ solve $B(S)^\top\lambda = c_B$ with $\lambda_0 = 0$. A basic surplus sets $\lambda_j = 0$; a basic structural gives $\lambda_j = 1 + \lambda_{j-1}$, so $\lambda_j$ is the length of the selected run ending at $j$. A missing nonbasic $u_j$ has reduced cost $1 + \lambda_{j-1} = a$, the length of its left component; a nonbasic $v_j$ at a selected edge has reduced cost $\lambda_j > 0$; $u_M$ has reduced cost $1 + \lambda_R > 0$. Every basis in the family is strictly dual feasible for every positive right-hand side. The constraints imply $u_j \ge u_{j+1} + b_j$ and $u_M \ge 0$; since every $u_i$ has positive cost the unique optimum is $u_j^\star = \sum_{s \ge j} b_s$, $u_M^\star = 0$, $v^\star = 0$, with $f^\star = \sum_j j b_j$ and the unique nondegenerate optimal basis $S^\star = \{1, \ldots, R\}$. Every missing edge requires a separate basis exchange.

\subsection{Legal pivots and exact solve-vector counts}

Let missing edge $j$ separate a left component of length $a$ and a right component $[j+1, v]$ of length $b$. Each inverse block of $U(S)$ is upper triangular with all upper entries equal to one, and $\Sigma$ only changes row signs; component $[h, v]$ corresponds to inverse indices $[\max(1, h-1), v-1]$. For the missing row $j$,
\[
B(S)^{-\top}e_j = -\sum_{s=j}^{v-1} e_s, \qquad \gamma_j = b, \qquad \delta_j = \sum_{s=j}^{v-1} b_s .
\]
Pricing this BTRAN result against the nonbasic columns gives coefficient $-1$ at $u_j$, $+1$ at the next missing $u_v$ (or $u_M$), nonnegative coefficients at surplus columns inside the block and zero elsewhere; the leaving $v_j < 0$ rises to its bound, so only negative coefficients are eligible and $u_j$ is the unique entrant; its reduced cost is $a > 0$, the dual step is $a$, and the new basis is exactly $B(S \cup \{j\})$. Every insertion order is a legal improving dual path with no ratio-test ties. Column $e_j - e_{j-1}$ has one source in each adjacent inverse block, reaching row $j$ on the right and $a - \chi$ rows on the left, where $\chi = \mathbf{1}\{h = 1\}$; the supports are disjoint. Two inverse rows have nonzero inner product exactly when they lie in the same block, since row $i$ of a block is $\sigma_i(0, \ldots, 0, 1, \ldots, 1)$ and the products cannot cancel. Therefore
\[
\nnz(B^{-1}a_{u_j}) = a + 1 - \chi, \qquad \nnz(B^{-\top}e_j) = b, \qquad \nnz(B^{-1}B^{-\top}e_j) = b,
\]
and a pivot contributes $a + 2b + 1 - \chi$ to $K$. Its DSE score is $\delta_j/\sqrt{b}$, which with a unit right-hand side equals $\sqrt{b}$; the merge order is therefore forced by the rule, not chosen.

\subsection{One right-hand side forces all paths strictly}

For even $q \ge 4$, $M = q^2$, $\varepsilon = 1/100$ and seeds $E = \{q-1, 2q-1, \ldots, M-1\}$ (all odd), define for odd $i$ and even $j$
\[
z_i = \begin{cases} 3 - i/M & i \in E \\ 0 & i \notin E \end{cases}, \qquad b_i = 1 + \varepsilon z_i, \qquad b_j = 1 - \varepsilon z_{j+1} + \varepsilon (M - j)/M .
\]
Every entry lies in $[1 - 3\varepsilon, 1 + 3\varepsilon]$ and $b_j + b_{j+1} = 2 + \varepsilon(M - j)/M$ for even $j$.

\textbf{Cold path.} From $S_0 = \emptyset$ (the all-slack basis), exact DSE forms blocks of $q$ consecutive vertices from left to right: within block $k$ it inserts $kq - 1, kq - 2, \ldots, (k-1)q + 1$ and then, for $k \ge 2$, inserts $(k-1)q$ to attach the block. Suppose the previous blocks are complete; every remaining right component is a singleton with score $b_j$; seed entries exceed $1 + 2\varepsilon$, nonseeds are at most $1 + \varepsilon$, and seed values decrease with index, so $kq - 1$ uniquely wins. It creates a right component of size two, and during the cascade exactly one missing boundary has a nonsingleton right component, with score at least $(1 - 3\varepsilon)\sqrt 2 > 1 + 3\varepsilon$, the maximum of any singleton; so it wins uniquely until the cascade reaches the prefix, including at the attachment step. Induction gives the whole path and its strictness.

\textbf{Warm B.} Omit only $F = \{M - 2r, M - 2r + 2, \ldots, M - 2\}$: components $(M - 2r, 2, \ldots, 2)$. Every missing boundary is even with score $(2 + \varepsilon(M-j)/M)/\sqrt 2$, strictly decreasing in $j$; inserting the leftmost joins the prefix to one pair and changes no other score, so all missing edges are inserted left to right with unique winners.

\textbf{Warm A.} Select $S_A = \{1, \ldots, M' - 1\}$ with $M' = M - r = (q - t)q$: components $(M', 1, \ldots, 1)$. This is the cold basis after its first $q - t$ blocks are complete and attached, so the same right-hand side and rule execute the last $t$ cold blocks with all their strict winners: each block grows from its seed by $q - 1$ internal pivots and then attaches, $tq = r$ pivots in all. (This compares exact algebraic DSE states; it does not assume that a solver restarted from a saved basis preserves any numerical history.)

\subsection{Counts and the corollary}

Every cold block contributes internal symmetric work $2 + \cdots + q = q(q+1)/2 - 1$, and attachments contribute the same total, so $W_0 = q^3 + q^2 - 2q$ and, with $\chi = 1$ exactly $q$ times, $K_0 = \tfrac32 W_0 + (M - 1) - q = (3q^3 + 5q^2 - 8q - 2)/2$. Warm B's merged lengths are $M - 2r + 2, \ldots, M$, so $N_B = r$ and $K_B = r(M - r + 3)$. For warm A each internal block pivot has $a = 1$, $\chi = 0$ and $b = 1, \ldots, q - 1$, and the attachment of block $k$ has left length $M' + (k-1)q$, right length $q$ and $\chi = 1$, giving $N_A = r$ and $K_A = 2t^5 - t^4/2 + 5t^3/2 - 2t$. Both overlaps equal $1 - r/(M-1)$, and substituting $M = t^4$, $q = t^2$, $r = t^3$ gives the limits in the theorem. Both warm starts select $u_1$ and every later inserted column has index at least two, so after $k$ pivots each basis has $2R - r - 1 + k$ nonzeros: the sparsity profiles agree. The valid dual bound $f(S) = \sum_j b_j\lambda_j$ gives both warm starts a relative gap of at most $4(1 + 3\varepsilon)r/[(1 - 3\varepsilon)(M-1)] = O(1/t)$ against a cold gap of one.

\textbf{Where warm B's work goes.} Each pivot of warm B attaches a pair to the correct prefix, whose length grows from $M - 2r$ to $M - 2$. Its BTRAN result and DSE FTRAN result have $b = 2$ nonzeros each, while the FTRAN result of the entering column spans the whole prefix ($a + 1 - \chi$ nonzeros; in network terms, the tree path closed by the entering arc). Almost all of $K_B$ is therefore the FTRAN of the entering column, which every dual simplex iteration forms whatever its pricing rule: the separation is a property of the path, and DSE only fixes the order of the merges. The bases along the two warm paths have equal numbers of nonzeros throughout, so what differs is the sparsity of $B^{-1}$ along the path. Warm B also costs $K_B/K_0 \sim 2t/3 = \Theta(M^{1/4})$ times the cold start, although it needs only $r/(M-1) \approx 1/t$ of the cold pivots.

\textbf{Robustness.} For a fixed basis row $R_i$ the score $s_i(b) = (-R_ib)_+/\|R_i\|_2$ is 1-Lipschitz in $b$; cold winner gaps are at least $\varepsilon/q$, cascade gaps at least $(1 - 3\varepsilon)\sqrt 2 - (1 + 3\varepsilon)$, warm-B gaps at least $\sqrt 2\varepsilon/M$, and all basic sign margins at least $1 - 3\varepsilon$. The eligible entrant and the dual ratio test depend only on $A$ and $c$, not on $b$. Hence every perturbation with $\|h\|_2 < \varepsilon/(\sqrt 2 M)$ preserves all three paths and all counts. The radius shrinks with $M$ and no uniform radius is claimed; at $t = 16$ it is about $10^{-7}$, the order of HiGHS' default feasibility tolerances. This completes Theorem~\ref{thm:chain}.

Corollary~\ref{cor:proxy}. The warm-B proof holds for any $1 \le r' < M/2$. With $r' = r - 1$ the start has overlap $1 - (r-1)/(M-1)$ and $r - 1$ pivots, so it is the unique overlap maximizer and pivot minimizer of the pair, and its work is $K'_B = (r-1)(M - r + 4) \ge t^7/2$ against $K_A \le 3t^5$, a ratio of at least $t^2/6 > 1$ for all $t \ge 4$. \qed

\textbf{Arithmetic model.} The chain matrix has at most three entries per row. With support-sensitive triangular solves, row-wise pricing over the rows in the support of the BTRAN result, updates proportional to the charged supports, and factor maintenance that changes one diagonal sign and at most one superdiagonal entry per pivot, forming the three vectors and updating costs $\Theta(K)$. Initial values, reduced costs and exact row norms follow from component recurrences in $O(M)$; a max-heap of infeasible boundaries gives exact full row selection in $O(\log M)$ per pivot. The total is $\Theta(M + K) + O(N \log M)$, and $K_A = \Theta(M^{5/4})$, $K_0 = \Theta(M^{3/2})$, $K_B = \Theta(M^{7/4})$ dominate. The unbounded ratio is a statement about this model. An implementation with a per-iteration term linear in the number of rows (a full CHUZR scan, column-wise PRICE or dense solves), or with refactorization after a fixed number of updates, bounds the ratio of total work by a constant, although the constant can be large. Computing every exact initial DSE weight by a separate BTRAN instead of by the recurrences materializes $\Theta(M^2)$ entries on both warm starts, which masks the separation (below). An independent exact checker reproduces 20 complete strict DSE trajectories across the four starts (46,934 pivots through 20,735 rows) and four dense exact paths with all inverse matrices, reduced costs, solve-vector supports and factor profiles.

\subsection{Initialization and the HiGHS realization}

The theorem counts continuation work after exact initial weights are in place, and it assumes the linear-time initialization that these bidiagonal bases allow. If every exact initial DSE weight is instead computed by a separate BTRAN, the initial output counts are $J_A = \binom{M-r}{2} + r$ and $J_B = \binom{M-2r}{2} + 3r$ against $J_0 = M - 1$ for the cold start. Both warm counts grow like $M^2/2$ and dominate the continuation, so $(J_B + K_B)/(J_A + K_A) \to 1$: in this cost model, overlap and pivot count miss the setup as well as the remaining path.

HiGHS computes exact initial DSE weights, one BTRAN per row, for any nonslack start that is not already dual feasible and near optimal; for such a near-optimal start it uses Devex instead. This holds in its default and DSE modes. In its default mode it can also switch from DSE to Devex during a solve in which the extra DSE solve stays much denser than the other solves.

\begin{table}[t]
\centering\footnotesize
\caption{The chain LP of Theorem~\ref{thm:chain} at $t = 16$ (65,535 rows) in HiGHS 1.13.1: wall time in seconds of installing the start and solving (one thread, presolve off, dual simplex; one repetition at this size). Warm A and warm B share overlap 0.94 and 4,096 pivots; the cold start takes 65,535. The first two settings also turn off cost perturbation and HiGHS' less-infeasible DSE row choice; the defaults keep automatic edge weights, scaling and that row choice.}
\label{tab:chain-highs}
\setlength{\tabcolsep}{6pt}
\begin{tabular}{@{}lcccc@{}}
\toprule
Setting & cold & warm A & warm B & warm B / warm A\\
\midrule
Devex forced, scaling off & 1.21 & 0.156 & 4.21 & 27.0\\
DSE requested, scaling off & 1.48 & 34.0 & 34.6 & 1.02\\
HiGHS defaults & 26.4 & 36.1 & 34.9 & 0.97\\
\bottomrule
\end{tabular}
\end{table}

Table~\ref{tab:chain-highs} gives the three starts at $t = 16$ under three settings. With presolve off, Devex forced from the first iteration and scaling off, the solve-vector supports that HiGHS reports match the exact DSE paths, and warm A and warm B, which share overlap $0.94$ and 4,096 pivots, take 0.16 and 4.2\,s, against 1.2\,s from the cold start: warm B needs 1/16 of the cold pivots and is still 3.5 times slower than the cold start. In the runs with scaling off, HiGHS' continuation instructions per pivot on the cold and warm-A paths fit $\approx 14\text{k} + 730\,t^2$, a fixed cost of about 14k instructions per pivot plus a term proportional to $K/N$, and on warm B $\approx 140\,M$, so the solver follows the model's scaling. The measured ratio of 27 is below the model's $K_B/K_A = 121$ at $t = 16$ because the solver's cost is not proportional to $K$ alone: the fixed per-pivot cost weighs most on warm A's short pivots, and the times include the setup. With DSE requested, the exact initialization makes both warm solves about 23 times slower than cold. Under the full defaults the paths change as well as the initialization: both warm continuations cost $\Theta(M)$ per pivot and their ratio stays near $1.8$, and the warm totals are within 4\% of each other at 1.3 to 1.4 times cold, whose own solve takes 26\,s under these settings. On our real LPs, whose crash bases are triangular and sparse, the exact initial weights cost a median 0.4\% of the solve (Appendix~\ref{app:periter}).

\subsection{Metric alignment of the policy's loss}
\label{app:loss}

Let $Q$ sample classes uniformly and then LPs from each class, and let $L_\delta = \log(C_\delta/C_{\text{ref}})$ be the log cost ratio of action $\delta$ against a common positive reference. For a policy whose draw is conditionally independent of solver randomness, $\mathbb{E}_Q \sum_\delta \pi_\theta(\delta \mid x) L_\delta$ equals the expected log cost ratio of the sampled deployment, and its exponential is the class-balanced geometric-mean ratio; at fixed $x$ the convex combination is minimized by putting all mass on a minimum-cost action, and a common policy-independent baseline (here $\log C_0$) preserves risk differences. Eq.~(\ref{eq:loss}) is this identity with instructions as $C$, the unshifted start as reference, a clip on the per-action loss, and uniform rather than class-balanced weights over the training LPs (the training set is roughly class-balanced by construction). The argmax deployment differs from the sampled policy; a one-line bound gives $\mu_{h(x)}(x) - \mu_\star(x) \le |\mathcal{A}| \sum_\delta \pi(\delta \mid x)(\mu_\delta(x) - \mu_\star(x))$ for the argmax $h(x)$, a population statement only.

\section{Within-LP study and landscape measurements}
\label{app:within}
\label{app:landscape}

This appendix gives the design and the full statistics of the within-LP study of \S\ref{sec:empirical-gap} (Table~\ref{tab:within}), and the measurements behind the cost per iteration of learned starts and the roughness of the cost landscape.

\textbf{Within-LP panel.} The panel has 381 LPs from the W1 part of \Dtw (25 classes; 179 have a primal or dual twin in the predictor's training data, and results without them are released with the code). Each LP is solved from up to 53 starts: 42 from the deployed decoder (33 shifts of the structural scores in $[-4, 4]$ at degree exponent $p = 0.5$, and the nine-shift grid at $p = 1.0$), ten high-overlap starts (sequential-softmax crash samples on $B^\star$ membership at temperatures 0.02 to 0.3), and the all-slack start; the median LP has 52 distinct starts, 17,329 in all. Every start of an LP was solved on one pinned core of an exclusive node with the deployed dual-simplex settings and a 120\,s limit, with two timing repetitions in seeded random order; iteration counts were identical across repetitions and the median repeat difference in wall time was 1.1\%. Overlap is measured against the nearest of the known optimal bases of the LP (the reference $B^\star$ and the distinct optimal final bases reached by the 53 solves, a median of 45 per LP; 78\% of LPs have more than one). Five of the 17,329 solves failed and are kept at twice the time limit (PAR2). Rank statistics are class-balanced with 2,000-replicate bootstrap intervals over classes and then LPs. The perfect-proxy regrets in \S\ref{sec:empirical-gap} take the best-time member among proxy ties. Table~\ref{tab:within} lists all statistics.

\begin{table}[h]
\centering\footnotesize
\caption{Within-LP statistics on the 381-LP panel of \Dtw (25 classes, 17,329 distinct starts). Rank and variance statistics are computed within each LP and pooled; intervals are class-balanced bootstraps over classes and then LPs.}
\label{tab:within}
\setlength{\tabcolsep}{4pt}
\begin{tabular}{@{}ll@{}}
\toprule
Statistic & Value [95\% CI]\\
\midrule
\multicolumn{2}{@{}l}{\emph{Kendall $\tau$ with time}}\\
\quad overlap & 0.31 \ci{0.22}{0.41}\\
\quad pivot count & 0.65 \ci{0.59}{0.71}\\
\quad instruction count & 0.91 \ci{0.89}{0.93}\\
\multicolumn{2}{@{}l}{\emph{Kendall $\tau$, other pairs}}\\
\quad overlap and pivot count & 0.45 \ci{0.36}{0.55}\\
\quad overlap and instructions & 0.30 \ci{0.20}{0.40}\\
\quad pivot count and instructions & 0.63 \ci{0.56}{0.71}\\
\midrule
\multicolumn{2}{@{}l}{\emph{Variance of log time explained ($R^2$)}}\\
\quad overlap alone & 0.145\\
\quad $\log(1 + \text{pivots})$ & 0.706\\
\quad $\log(1 + \text{pivots})$ and overlap & 0.707\\
\quad log instructions & 0.985\\
\midrule
\multicolumn{2}{@{}l}{\emph{Matched pairs: overlap within 0.005, pivots within 2\% (26,908 pairs)}}\\
\quad share with time ratio $> 1.25$ & 0.092 \ci{0.043}{0.152}\\
\quad share with time ratio $> 2$ & 0.019 \ci{0.000}{0.055}\\
\quad quantiles 50/75/90/99\% & 1.03 / 1.08 / 1.17 / 1.72\\
\quad largest ratio & 11.6\\
\midrule
\multicolumn{2}{@{}l}{\emph{Proxy-dominance reversals: higher overlap and fewer pivots, yet costlier}}\\
\quad dominant start $> 5\%$ costlier & 0.122 \ci{0.082}{0.170}\\
\quad dominant start $> 25\%$ costlier & 0.054 \ci{0.029}{0.084}\\
\midrule
\multicolumn{2}{@{}l}{\emph{Cost of the proxy-optimal start, relative to the fastest of 42 decoder starts}}\\
\quad highest overlap & 1.275 \ci{1.18}{1.39}\\
\quad fewest pivots & 1.089 \ci{1.054}{1.130}\\
\quad best under both (exists on 31\% of LPs) & 1.087 \ci{1.043}{1.149}\\
\bottomrule
\end{tabular}
\end{table}

\begin{figure}[h]
\centering
\includegraphics[width=0.52\linewidth]{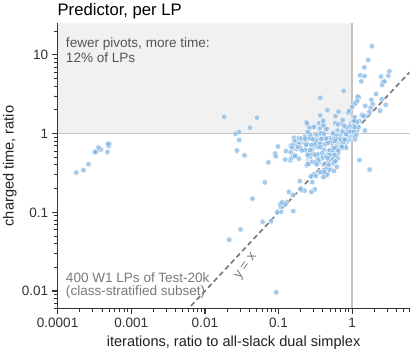}
\caption{Per-LP iteration and charged-time ratios of the predictor against the all-slack dual simplex on a class-stratified subset of 400 W1 LPs of \Ttw (two timing repetitions averaged). Points in the upper left take fewer pivots and more time.}
\label{fig:perlp}
\end{figure}

\textbf{Per-pivot cost.} Figure~\ref{fig:perlp} shows the gap between iterations and time per LP on this 400-LP subset of \Ttw: on most LPs the predictor saves a far larger share of pivots than of time. On two further development panels, 874 LPs at the 20k size and 477 LPs at 40k to 60k variables (the latter including 23 duplicates of training LPs), the variance of $\log(T_{\text{predictor}}/T_{\text{all-slack dual}})$ decomposes by covariance share into iterations (70\% and 78\%), cost per iteration (29\% and 9\%) and start overhead (1\% and 13\%). Pivots from the learned start cost 1.50 times as much as pivots from the all-slack start (geometric mean) on both panels, and 1.34 times on the 283-LP selection split of \Dtw. The learned basis has $0.80$ to $0.84$ times the nonzeros of $B^\star$, and basis density does not track the ratio of the cost per iteration (Spearman $+0.07$). In instrumented solves, the cumulative solve-vector work (Appendix~\ref{app:proofs}), estimated from HiGHS' kernel statistics as calls times mean result density and dominated by the DSE FTRAN and the BTRAN, explains paired instruction differences beyond iterations, overlap and initial fill (within-LP leave-one-class-out $R^2$ from 0.51 to 0.88), while no static property of the start predicts it (leave-one-class-out $R^2 \le 0$).

\textbf{Per-iteration traces.} To see directly which iterations a learned start removes, we recorded the instructions of every iteration through HiGHS' simplex callback (the traced solves reproduce the plain solves' iterations and objectives exactly, and the callback's overhead, typically 0.3\%, is subtracted) on 20 LPs of \Dtw, one or two per class, for the native solver, the all-slack dual simplex, the predictor's start with the dual simplex and \oursF. For a reference path with $N_R$ iterations and a learned path with $N_L < N_R$, we compare the reference's work per iteration over its first $N_R - N_L$ iterations with its work per iteration over its last $N_L$, and the learned path's work per iteration with the latter; setup is excluded. Against the all-slack dual simplex, the first iterations cost $0.44\times$ \ci{0.28}{0.66} as much as the last ones for the predictor's start (below 1 on 15 of 17 LPs) and $0.45\times$ \ci{0.26}{0.72} for \oursF (14 of 16), and the learned paths' iterations cost $0.89\times$ \ci{0.71}{1.17} and $0.81\times$ \ci{0.62}{1.06} as much as the reference's last ones; against the native solver the figures are $0.42$ and $0.38$, and $0.92$ and $0.72$. The learned starts skip the cheap early part of the path and pay about the price of its late part for each of their own iterations, which is why their mean cost per iteration is higher although their total work is lower. Installing a learned basis (\texttt{setBasis}, factorization and initial weights) costs about five times the native solver's setup, 11\% of the learned solves' instructions against 1.2\% for the native solver. Figure~\ref{fig:itertrace} in \S\ref{sec:empirical-gap} shows eight of the traces, including three classes where the learned starts lose. The test does not apply where a learned start takes more iterations, and a few LPs go the other way.

\textbf{Roughness and basin.} This panel has 312 LPs at the 20k size that the predictor never saw (24 classes), with 33 candidate starts each. On it, repeat solves are bitwise identical in iterations and objective and agree to $5\times10^{-6}$ in log instructions at the median, so all variance among candidates is real.
Swap statistics are measured on 60 of the 312 LPs and direct installations on 200.
\begin{itemize}
\item One random valid swap from the predictor's start changes log instructions by a median of 0.042 (90th percentile 0.28); 24\% of swaps are better and 24\% worse by more than 5\%. $|\Delta|$ grows only slowly with the number of swaps, while the mean drifts upward.
\item The best of eight stochastic-crash samples takes $0.86\times$ the instructions of the deterministic start, but the mean sample is no better than it. A cross-fitted selector on cheap features of the start captures 34\% of the gap to the best of the 33 candidates within the training classes and about 10\% outside them.
\item An exact optimal basis installed directly (HiGHS repairs singular sets on install; charged) solves in $0.051\times$ \ci{0.036}{0.075} the predictor's instructions, with 86\% of LPs solving in zero iterations; such a basis is primal and dual feasible, so HiGHS also skips the initial DSE weights and the cost perturbation. The same basis through the crash costs $0.595\times$, because the crash drops a median of 221 of its basic columns.
\item Random swaps of 0.1\%, 1\%, 5\% and 20\% of the optimal basis cost $0.30$, $0.67$, $1.47$ and $2.87\times$.
\item Model-generated complete bases sit about 6\% (about 900 swaps) from an optimum and cost $1.03$ to $1.09\times$ installed directly.
\end{itemize}

\section{Predictor, decoder and shift policy}
\label{app:generator}
\label{app:head}

This appendix details the predictor and its training (\S\ref{sec:generator}), the decoder (Algorithm~\ref{alg:crash} in \S\ref{sec:generator}), and the shift policy with its training, ablations (Table~\ref{tab:shift-ablations}) and controls (Table~\ref{tab:controls}). Appendix~\ref{app:mechanism} explains where the shift policy's gain arises. Table~\ref{tab:hyper} collects all settings.

\textbf{Graph and features.} The LP is a graph with one node per structural column, per slack column and per constraint, and an edge for every nonzero of $A$, in the style of \citet{gasse2019exact}. Node and edge features are normalized objective coefficients, bounds, right-hand sides, degrees and matrix entries. The target is membership in the optimal basis that the solver returns from the all-slack start, which is one optimal basis among possibly several.

\textbf{Predictor training.} Stage 1 trains on the 114,874 LPs of \Trs (37 classes, median 830 columns and 612 rows after presolve) in three chained runs of 2, 4 and 6 epochs (43k steps, 8.8 H100-hours) with TF32 matrix multiplication. The loss per LP is the basis-membership cross-entropy over all $n + m$ columns of $[A\ I]$, divided by $m$, plus $0.2$ times the cross-entropy of an auxiliary dual-tightness output and $0.2$ times that of the nonbasic-bound outputs. Stage 2 trains on 3,457 LPs of \Trtw (26 classes, median 18.6k columns) for 6 passes at a constant learning rate, in fp32 with deterministic kernels. A checkpoint is saved every 37 steps and selected on the 283-LP selection split of \Dtw by the paired wall time of \texttt{setBasis} and the solve against the model after stage 1 (best at step 333, ratio 0.847). Ablations:
\begin{itemize}
\item Overlap of the crash with $B^\star$ rises from 0.645 (a 4-layer model) to 0.787, and precision at $|B^\star|$ to 0.76; the crash on $B^\star$'s own membership reaches 0.93.
\item Depth: 4 vs.\ 8 layers at matched data give overlap 0.730 vs.\ 0.742. Twelve layers add 0.004 at 1.2 to 1.3 times the forward cost and take $1.043\times$ \ci{0.95}{1.14} the charged time of 8 layers on \Dtw (one timing repetition). GATv2, TransformerConv and GINE backbones agree within 0.005 at depth 4.
\item Length of stage 1: 2, 6, 12, 18 and 24 epochs give overlap 0.761, 0.780, 0.787, 0.791 and 0.793; charged time saturates at 12 epochs (18 vs.\ 12: $1.018\times$ \ci{0.97}{1.08}).
\item Retraining stage 2 on the deduplicated 9,511 LPs of \Trtw at the same number of LP visits gives overlap 0.787 and precision 0.763, identical to the deployed checkpoint.
\end{itemize}

\textbf{Deployment.} The forward pass runs one LP at a time (microbatching perturbs scores by $10^{-6}$ and flips crash decisions) on a fused GATv2 kernel with TF32 matrix multiplication: 8-layer inference falls from 19.5 to 6.4\,ms median. TF32 changes 87\% of starts, so every policy, gate and panel uses the deployed path.

\textbf{Decoder.} Algorithm~\ref{alg:crash} gives the triangular crash. Consider the bipartite graph between structural columns and rows in which an edge $(j, i)$ exists when column $j$ has an entry in row $i$ of at least $10^{-8}$ times the column's largest entry; its gain is $g_{ji} = (s_j + \delta - s_{n+i})/\deg(j)^{p}$ with $p = 0.5$, where $s_{n+i}$ is the score of row $i$'s slack and $\deg(j)$ the number of eligible entries of column $j$. Selecting column $j$ in row $i$ removes every remaining edge of column $j$ and every edge whose row lies in the support of column $j$, so the selected columns form a triangular submatrix with nonzero diagonal, and the basis is nonsingular without repair. A shift $\delta$ raises column $j$'s gains by $\delta/\deg(j)^p$, which tends to admit more structural columns for $\delta > 0$ and fewer for $\delta < 0$, and tilts the choice toward low-degree columns. A finer grid of shifts, and a second per-LP decoder dimension (the exponent $p$ alongside $\delta$), raise the oracle but not the learned policy (Table~\ref{tab:shift-ablations}).

\textbf{Decoder tie rules and statuses.} Within a column, candidate rows are ordered by the score of their slack, ascending, with ties broken by row index; heap ties go to the smaller column index. Nonbasic statuses: lower-only bound $\to$ lower, upper-only $\to$ upper, free $\to$ zero, fixed $\to$ lower, boxed $\to$ upper if and only if the logit of the bound-status output is positive. The reference Python and the native implementation produce identical selections. Statuses are not a confound for any result: forcing a start's statuses to those of another method leaves instruction counts bitwise unchanged on a development panel.

\textbf{Shift policy: training.} The cost of an action is the number of user-space instructions, read from \texttt{perf\_event} around \texttt{setBasis} and the solve and summed; on the cost table, repeated solves agree to a median of $5\times10^{-5}$ and at most $6\times10^{-4}$ in log instructions. A solve is a failure if it is not optimal, if its basis was changed on install, or if its objective misses the reference by more than $10^{-6}$ relative; a failure costs $\log(2 \cdot 120\,\text{s} \cdot 5\times10^9)$. The policy is an MLP $84 \to 128 \to 128 \to 9$ (29k parameters) with SiLU and dropout 0.1, whose last layer is initialized so that the untrained policy plays $\delta = 0$; the two larger policies have 11 and 20 outputs and ten-seed ensembles (Appendix~\ref{app:periter}). The policy stage costs 1.7\,ms per LP (3.9\,ms with the novelty gate). The training data are the 7,699 W1 LPs of \Trtw times 9 actions, drawn from a 131,706-solve table over 14,634 LPs that took 56 minutes on 64 cores (about 60 core-hours; 71 solves were failures); \oursA adds the costs of the native solver and primal simplex and the W2 LPs, 9,349 LPs in all. Training uses AdamW (learning rate $10^{-3}$, weight decay $10^{-4}$), batches of 256 LPs and 100 epochs; the epoch is chosen on the 283-LP selection split of \Dtw (a held-out fifth of the training LPs gives the same numbers), with three seeds and the argmax of the mean log-probability. The loss clip at 3 keeps $20\times$ blow-ups visible; clipping at 1 hid them, and that policy played catastrophic shifts on six bin-packing LPs of a development panel.

\textbf{Shift policy: ablations.} Table~\ref{tab:shift-ablations} lists the ablations. No change of the recipe moves the result outside seed noise, a finer or two-dimensional action grid raises the oracle but not the policy, and the policy does not transfer to held-out classes without the novelty gate (on load balancing it plays $\pm 3$ on every LP). Adding proxies of the formulation and of recall to the inputs helps within the training classes by two points and does not transfer, and a margin rule on the ensemble is flat from 0 to 4.

\begin{table}[h]
\centering\footnotesize
\caption{Ablations of the shift policy: instructions relative to the unshifted predictor (class-balanced geometric means), on the 283-LP selection split of \Dtw for the first two blocks. Single seeds unless noted. The action-set rows compare with the nine shifts on the same LPs. Transfer rows hold out each class, formulation or dual form of the training LPs in turn (three seeds, epoch chosen on a held-out fifth of the training LPs); the worst held-out classes are knapsack (2.10), telecom (2.17) and load balancing (2.01).}
\label{tab:shift-ablations}
\setlength{\tabcolsep}{5pt}
\begin{tabular}{@{}ll@{}}
\toprule
Change & Ratio\\
\midrule
\multicolumn{2}{@{}l}{Recipe}\\
\quad main recipe, 5 seeds (mean $\pm$ sd) & $0.896 \pm 0.013$\\
\quad\quad ensemble of the 5 seeds & 0.900 \ci{0.839}{0.961}; pivots 0.920; summed instructions 0.905\\
\quad hidden width 64 / 256 / 512 / 1024 & 0.900 / 0.900 / 0.904 / 0.925\\
\quad learning rate $3\times10^{-4}$ / $3\times10^{-3}$ & 0.896 / 0.898\\
\quad loss clip 1 / 6 / none & 0.903 / 0.891 / 0.888\\
\quad dropout 0 / 0.3 & 0.893 / 0.902\\
\quad training LPs 1k / 3k / 7k (per seed) & 0.919 / 0.913 / 0.898\\
\quad ensemble of ten seeds & 0.879\\
\multicolumn{2}{@{}l}{Action set}\\
\quad shifts in $[-1, 1]$ (5 actions) & 0.954 (oracle 0.833)\\
\quad shifts in $[-4, 4]$ in steps of 0.25 (33 actions) & 0.917 vs.\ 0.912 for the nine shifts\\
\quad nine shifts $\times$ degree exponent $p \in \{0.5, 1\}$ & 0.883 vs.\ 0.889 (oracles 0.718 vs.\ 0.767)\\
\multicolumn{2}{@{}l}{Transfer}\\
\quad leave one class out & 1.144 to 1.171\\
\quad leave one formulation out (39 formulations) & 1.027 \ci{0.91}{1.19}; with novelty gate 0.996 (99\% abstention)\\
\quad leave one form (primal or dual) out & 1.154; with novelty gate 1.000\\
\bottomrule
\end{tabular}
\end{table}

\begin{table}[h]
\centering\footnotesize
\caption{Controls for the shift policy on the 283-LP selection split of \Dtw (22 classes): instruction ratios to the unshifted predictor. The oracle picks the best of the nine shifts per LP.}
\label{tab:controls}
\setlength{\tabcolsep}{5pt}
\begin{tabular}{@{}lc@{}}
\toprule
Rule & Ratio\\
\midrule
best fixed shift chosen on training data & 1.000 ($\delta = 0$)\\
every other fixed shift & 1.005 to 1.090\\
best fixed retune of the decoder (shift and degree exponent) & 0.967\\
best constant shift per class, chosen in hindsight & 0.926\\
learned policy, score statistics only (no embedding) & 0.975\\
learned policy, embedding only & 0.897\\
learned policy as deployed (3-seed ensemble; single seeds $0.896 \pm 0.013$) & 0.898\\
learned policy, 10-seed ensemble & 0.879\\
oracle over the nine shifts & 0.767\\
learned policy, leave-one-class-out & 1.144 to 1.171\\
\quad with the novelty gate & 1.000\\
\bottomrule
\end{tabular}
\end{table}

\textbf{What the learning contributes.} Table~\ref{tab:controls} shows that the shift policy's gain is a learned per-LP decision and not a retuning of the decoder. Every fixed shift is worse than none, the best decoder retune gains 3\%, the best constant per class chosen in hindsight gains 7\%, and the learned per-LP policy gains 10\% (12\% with ten seeds), 41\% of the oracle's log gain. The predictor's embedding carries the signal. On the predictor side, training length and the stage-2 fine-tune reduce solve time and greater depth does not: on \Dtw the model trained only on \Trs is at parity with the native solver, the fine-tune takes it to $0.82\times$, and a 12-layer network at matched budget is 4\% slower than the 8-layer one.

\subsection{Mechanism: a learned control of structural recall}
\label{app:mechanism}

\begin{figure}[h]
\centering
\includegraphics[width=0.85\linewidth]{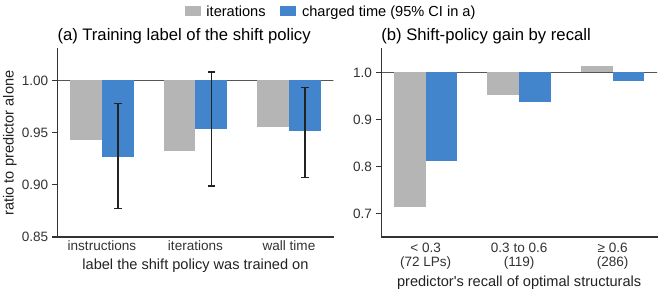}
\caption{(a) Shift policies that differ only in their training label, on \Dtw (charged time, three repetitions; Appendix~\ref{app:objective}): the policy trained on iterations removes more pivots and saves less time. (b) The shift policy's gain over the predictor as a function of the predictor's recall of the optimal structural columns, on a development panel of 477 LPs at 40k to 60k variables that includes 23 duplicates of training LPs, with a shift policy trained on training data that include duplicates. Panel (b) only locates where the gain arises; no reported comparison uses it.}
\label{fig:mech}
\end{figure}

A uniform shift helps although \S\ref{sec:gap} found the cost landscape rough, because the landscape is rough at the scale of single swaps but smooth along one coarse direction, the number of structural columns the crash admits. The predictor's losses concentrate on LPs where it admits too few of the optimal structural columns. Recall, the fraction of $B^\star$'s structural columns that a start contains, is the start property whose effect transfers best across classes (leave-one-class-out $R^2$ of 0.3), and it interacts with degeneracy: degenerate LPs with recall below 0.3 run at 1.3 to 2.3 times the all-slack dual simplex, while degenerate LPs with adequate recall behave like nondegenerate ones ($0.7\times$). Alternative optima are not the explanation (recall against the nearest of six optima found rises by only 0.02), and neither is the crash's stopping rule: the predictor scores the right kind of columns highly, but too few of them. The shift policy's gain is concentrated on these LPs (Figure~\ref{fig:mech}b). Where the predictor's recall is below 0.3, the policy runs in $0.81\times$ its time with $0.71\times$ its iterations, against $0.98\times$ where recall exceeds 0.6. The shift policy is thus a learned per-LP control of how many structural columns to commit to, chosen against what that commitment costs the solver.

The harm gate targets the same LPs, which is why it subsumes most of the shift policy's gain in \oursS: fit to the cost of the unshifted start, it hands many of them to the native solver before the shift is used. On the test sets, with the gate present, adding the shift policy changes the result by less than 4\% (\Ttw: no change; \Tfif and \Tei: 2\% slower; an 80-LP W2 subset of \Ttw: 3.5\% faster). Without the gate, the shift policy takes $0.94\times$ the predictor's time on \Ttw but moves it by at most one point on \Tfif and \Tei, because the novelty gate sends most of those LPs to the unshifted start.

\section{Decision policies with algorithm and pricing actions}
\label{app:periter}

This appendix covers the design of the action sets $\mathcal{A}_2$ and $\mathcal{A}_3$ (Table~\ref{tab:actionsets}), the training of \oursA and \oursF and their development readouts (Table~\ref{tab:devreadouts}), the probe of solver options behind the pricing action (Table~\ref{tab:probe}), and the decomposition of Figure~\ref{fig:decomp} (Table~\ref{tab:decomp}).

\textbf{Design of the action sets.} We include the native solver and primal simplex because, on a held-out-class panel, no learned start beats the solver's own choice of algorithm on some classes. Per LP, the oracle over $\mathcal{A}_2$ reaches $0.67\times$ the unshifted start's instructions on the 275 LPs of the probe below, against $0.76\times$ over the shifts alone; adding dual simplex from the all-slack basis, which nearly duplicates the native solver, or the classical crash made the learned policy worse on held-out training LPs. The oracle over $\mathcal{A}_3$ reaches $0.54\times$, and its gain over $\mathcal{A}_2$ lies almost entirely in the work per iteration (iterations stay at $0.72\times$ while the work per iteration falls from $0.92\times$ to $0.74\times$). The decoder offers no such lever: every crash start is exactly triangular, and static properties of the start explain about 2\% of the within-LP variance of work per iteration beyond the iteration count, so the policy must reach the cost per iteration through a solver option. Every design choice was made on a fixed held-out fifth of the training LPs (Table~\ref{tab:actionsets}).

\begin{table}[h]
\centering\footnotesize
\caption{Design of the action sets on a held-out fifth of the training LPs: instructions relative to the predictor (class-balanced geometric means). For \oursA, W1 is reported on two sub-splits. Configuration only: a selector among the solver's own starts under both pricing rules, without the learned start.}
\label{tab:actionsets}
\setlength{\tabcolsep}{6pt}
\begin{tabular}{@{}lcc@{}}
\toprule
Action set & W1 & W2\\
\midrule
\multicolumn{3}{@{}l}{Design of \oursA}\\
\quad shifts only & 0.906 / 0.935 & 0.851\\
\quad shifts and native solver & 0.810 / 0.781 & 0.792\\
\quad shifts, native solver and primal simplex ($\mathcal{A}_2$, chosen) & 0.783 / 0.754 & 0.605\\
\quad $\mathcal{A}_2$ plus all-slack dual simplex and classical crash & 0.796 / 0.773 & 0.648\\
\multicolumn{3}{@{}l}{Design of \oursF}\\
\quad start and algorithm ($\mathcal{A}_2$) & 0.785 & 0.645\\
\quad plus the nine shifts under Devex ($\mathcal{A}_3$, chosen) & 0.679 & 0.666\\
\quad plus native solver under Devex & 0.690 & 0.850\\
\quad shifts under both rules, without the strategy actions & 0.732 & 1.017\\
\quad configuration only & 1.166 & 1.895\\
\bottomrule
\end{tabular}
\end{table}

\textbf{Training of \oursA and \oursF.} The training rows of \oursA are the nine shifts of the shift policy's table plus the two strategy actions (the native solver, and the unshifted start with primal simplex), each relative to the same collection's unshifted dual-simplex solve, clipped at $\pm 3$ with failures at $+3$, on 9,349 LPs of \Trtw from W1 and W2. The table of \oursF adds, per training LP, the nine shifts under Devex and the native solver under both pricing rules (the latter used only by the configuration-only control), bridged to the existing rows through each collection's own unshifted solve under default pricing; the table covers 7,699 W1 and 1,501 W2 LPs, 9,197 of them with every action. Both policies use the shift policy's MLP, inputs and loss, with ten-seed ensembles, 150 epochs, learning rate $3\times10^{-3}$, dropout 0.1, batches of 256 LPs and window weight 1. Instruction counts agree across CPU models. The 20-action table took 184k solves and 373 core-hours of loaded single-thread time, $15.7$ times the instructions of one all-slack dual-simplex solve per LP ($20.1\times$ in W1, $14.5\times$ in W2); the eleven actions of $\mathcal{A}_2$ account for $8.5\times$ and the Devex shifts for $7.2\times$.

\textbf{Design readouts.} For \oursA, sixty training epochs undertrained the policy (0.870 / 0.883 / 0.719 in the layout of Table~\ref{tab:actionsets}); hidden width 256 and 300 epochs were within seed spread (three disjoint ten-seed ensembles differ by about 0.02). Leave-one-class-out without the novelty gate gives 1.196; with the gate it gives 1.000 when the gate falls back to the unshifted start, 1.452 to the native solver and 2.266 to primal simplex, so the fallback stays the unshifted start. The held-out fifth used for these readouts holds 1,531 W1 LPs of 25 classes and 316 W2 LPs of 11: \oursF took $0.669\times$ \ci{0.589}{0.751} and $0.668\times$ \ci{0.526}{0.826} the predictor's instructions and \oursA $0.822\times$ and $0.615\times$; paired, \oursF took $0.813\times$ \ci{0.733}{0.897} the instructions of \oursA in W1, with iterations $1.00\times$ and work per iteration $0.81\times$, and $1.086\times$ \ci{0.946}{1.359} in W2.

\begin{table}[h]
\centering\footnotesize
\caption{Development readouts of the two larger policies: charged time relative to the predictor (class-balanced, 95\% intervals where available). The \Dtw panel has 207 W1 LPs of 23 classes and 45 W2 LPs of 8; \Dfif has 290 W1 LPs of 23 classes and 164 W2 LPs of 16. No training LP is larger than 25,000 variables.}
\label{tab:devreadouts}
\setlength{\tabcolsep}{2.5pt}
\begin{tabular}{@{}lcccc@{}}
\toprule
& \multicolumn{2}{c}{\Dtw panel} & \multicolumn{2}{c}{\Dfif}\\
\cmidrule(lr){2-3}\cmidrule(l){4-5}
Method & W1 & W2 & W1 & W2\\
\midrule
\oursF & 0.761 \ci{0.618}{0.912} & 0.524 \ci{0.409}{0.666} & 0.891 \ci{0.813}{0.971} & 0.857 \ci{0.719}{0.993}\\
\quad + check and harm gate & -- & -- & 0.751 & 0.682\\
\oursA & 0.801 & 0.462 & 0.941 \ci{0.88}{1.00} & 0.856 \ci{0.73}{0.97}\\
\quad + check and harm gate & -- & -- & 0.796 & 0.679\\
predictor, Devex on every LP & 0.980 & 1.189 & 0.865 \ci{0.730}{1.019} & 1.372 \ci{0.943}{2.326}\\
configuration only & 1.164 & 0.904 & -- & --\\
\bottomrule
\end{tabular}
\end{table}

Table~\ref{tab:devreadouts} gives the charged development readouts. On the \Dtw panel, \oursF's gain in W1 splits into iterations $0.94\times$, time per iteration $0.77\times$ \ci{0.67}{0.88} and fixed overhead $1.06\times$. \oursA was designed on the selection split, and the W1 part of this panel shares 98 LPs with it, so the selection split and that W1 part are development readouts only; the clean evaluations are W2 of the panel, \Dfif and the test sets. On \Dfif, with the harm gate applied to both, \oursA takes $0.98\times$ and $0.97\times$ the gated predictor's time on W1 and W2. On W2 of \Dfif, \oursF takes $1.000\times$ \ci{0.955}{1.051} the time of \oursA and $0.604\times$ \ci{0.360}{0.872} that of the predictor with Devex on every LP. On W1 it takes $0.923\times$ \ci{0.853}{0.996} the instructions of \oursA and ties Devex on every LP ($1.031\times$ \ci{0.852}{1.244} in time), because the novelty gate sends 152 of the 290 LPs to its fallback under default pricing. A variant whose fallback uses Devex helped on W1, where it was designed ($0.889\times$ the instructions of \oursF), but not on W2 ($1.076\times$ \ci{0.954}{1.224} in time), so we keep the default fallback.

\textbf{Harm gate with the policies.} When combined with \oursA or \oursF, the harm gate is evaluated on every LP, on the start the policy deployed (for the primal and native actions, the unshifted start); if it predicts too small a saving, the native solver runs instead, and the policy's start, the gate features and the gate evaluation are charged. The gate was trained on the unshifted start with dual simplex and default pricing, so it judges primal and Devex choices with a model fit to another action; on \Ttw many of its overrides reverse primal-simplex choices (137 of 286 for \oursA, 117 of 156 for \oursF). \oursA and \oursF can choose the native solver themselves. Adding the harm gate slows them at the training size (\oursA $0.54\times$ instead of $0.51\times$ native, \oursF $0.52\times$ instead of $0.49\times$ on \Ttw), because their own native action covers most LPs the gate would hand over while the gate's features and evaluation are charged. At larger sizes the gate is needed, because the novelty gate sends most LPs to the unshifted start, where the policy never gets to choose the native solver.

\textbf{Pricing on the test and development panels.} On \Ttw, \oursF takes $0.72\times$ \ci{0.63}{0.81} the predictor's time and $0.80\times$ \ci{0.70}{0.90} that of \oursS. Against \oursA it is at $0.96\times$ \ci{0.88}{1.05}, with 15\% more iterations that each cost 18\% less. By window (Table~\ref{tab:windows-test20k}), \oursF takes $0.55$, $0.36$ and $0.33\times$ the native solver on W1, W2 and W3, against $0.60$, $0.33$ and $0.22\times$ for \oursA, so its gain is concentrated on short solves. Development data showed the same pattern, on held-out training LPs ($0.81\times$ the instructions of \oursA, iterations unchanged) and on the \Dtw panel ($0.95\times$ \ci{0.83}{1.08} the time of \oursA, $0.81\times$ its time per iteration). The predictor's start with Devex on every LP takes $1.07\times$ the predictor's time on \Ttw, and a selector trained with the same loss to choose among the solver's own starts and the two pricing rules takes $1.26\times$; \oursF takes $0.67\times$ and $0.57\times$ their times. Table~\ref{tab:actions-test20k} lists the actions the three systems take on \Ttw.

\textbf{Probe of solver options.} On 275 LPs of the selection split, the predictor's start was solved under 14 option settings, the default included (Table~\ref{tab:probe}). Devex was the only one that clearly helped. Devex relative to default pricing ranges by class from 0.51 (knapsack, a pure per-iteration effect) to 2.39 (project selection). Wall-time ratios measured under 39 concurrent workers track the instruction ratios with a per-LP correlation of 0.987. Setup under default pricing costs $1.12\times$ the Devex setup; the difference, the cost of the exact initial DSE weights, is a median 0.4\% of the solve (90th percentile 3.6\%). At this size Devex's gain is therefore cheaper iterations rather than avoided weight initialization, unlike on the chain LP of Appendix~\ref{app:proofs}. Per-LP oracles relative to the predictor: shifts 0.764; shifts under both rules 0.589; $\mathcal{A}_2$ 0.667; $\mathcal{A}_3$ 0.538, a paired $0.807\times$ \ci{0.736}{0.879} beyond $\mathcal{A}_2$ with iterations unchanged (0.726 to 0.724) and work per iteration from 0.919 to 0.743.

\begin{table}[h]
\centering\footnotesize
\caption{Probe of HiGHS options on the predictor's start (275 LPs of the selection split of \Dtw): instructions, iterations and work per iteration relative to HiGHS' default options (class-balanced geometric means), and failures. HiGHS' default dual edge weights use DSE with automatic switches to Devex.}
\label{tab:probe}
\setlength{\tabcolsep}{5pt}
\begin{tabular}{@{}lcccc@{}}
\toprule
Option & instructions & iterations & work per iteration & failures\\
\midrule
Devex pricing & 0.910 & 1.067 & 0.853 & 0\\
Dantzig pricing & 1.846 & 2.073 & 0.891 & 7\\
DSE without the switch to Devex & 1.057 & 0.982 & 1.076 & 0\\
price strategy 1 / 2 & 1.010 / 0.994 & 1.000 & 1.010 / 0.994 & 0\\
cost perturbation off / $\times 5$ & 1.160 / 1.046 & 1.121 / 1.035 & 1.035 / 1.011 & 1 / 0\\
factor pivot threshold 0.01 / 0.5 & 1.007 / 1.014 & 1.005 / 1.011 & 1.003 / 1.003 & 0\\
less-infeasible DSE check off & 0.992 & 0.981 & 1.011 & 0\\
scaling off & 1.084 & 1.098 & 0.988 & 1\\
update limit 100 & 1.149 & 1.011 & 1.136 & 0\\
Devex and cost perturbation off & 1.115 & 1.263 & 0.883 & 1\\
\bottomrule
\end{tabular}
\end{table}

\textbf{Headroom inside the decoder and what determines work per iteration.} On the within-LP panel of Appendix~\ref{app:within}:
\begin{itemize}
\item Oracle headroom. The instruction oracle over the nine shifts is $0.947\times$ \ci{0.919}{0.970} the iteration oracle, with $1.04\times$ the iterations and $0.91\times$ the work per iteration, and over the 42 decoder starts it is $0.909\times$ (1.10 and 0.83). Within an LP, log iterations and log work per iteration correlate with a median of $-0.54$, so starts with fewer iterations tend to have costlier ones.
\item Triangularity and setup. Every one of the 17,329 starts is exactly triangular (no nucleus after singleton peeling, so no fill beyond the basis), and setup (installing and factorizing the start and initializing the pricing weights) is a median 3.1\% of the solve's instructions. Within an LP, log work per iteration has $R^2 = 0.433$ against log iterations alone and $0.457$ with the structural count, basis nonzeros, the nonzeros of a COLAMD LU factorization and the measured factorization cost added.
\item Noise floor. Calling \texttt{setBasis} with the same basic set in three selection orders gives bytewise identical statuses. Setup instructions change by a median 0.3\% under random row and column relabelings of the LP, and whole-solve instructions by a median 2.1\% through tie-breaking; the better of two relabelings of one start is $0.971\times$ the other, a noise floor for any oracle gain.
\end{itemize}
Offline, on the 283-LP selection split, the shift policy's gain over the predictor splits into iterations $0.907\times$ and work per iteration $0.974\times$ \ci{0.916}{1.036}, that of \oursA into $0.808\times$ and $0.999\times$, and the gain of the instruction-trained shift policy over the iteration-trained one into iterations $1.010\times$ and work per iteration $0.932\times$ \ci{0.872}{0.991}; in the charged label ablation (Appendix~\ref{app:objective}) the last split is $1.012\times$ and $0.957\times$ \ci{0.925}{0.987}.

\textbf{Decomposition on \Ttw.} Figure~\ref{fig:decomp} uses the per-LP rows of the \Ttw panel with the check inline, restricted to LPs where both the method and the native solver are valid and take at least one iteration (all 1,356 for every method shown except the configuration-only control, 1,353); intervals are two-level bootstraps over classes and LPs. Table~\ref{tab:decomp} gives four paired contrasts split the same way by Eq.~(\ref{eq:decomp}). The analysis script is released with the code.

\begin{table}[h]
\centering\footnotesize
\caption{Paired contrasts on \Ttw, split by Eq.~(\ref{eq:decomp}) into iterations, time per iteration and fixed overhead (class-balanced geometric means of the first method relative to the second, 95\% intervals where computed).}
\label{tab:decomp}
\setlength{\tabcolsep}{4pt}
\begin{tabular}{@{}lcccc@{}}
\toprule
Contrast & time & iterations & time per iteration & fixed\\
\midrule
\oursF vs.\ \oursA & 0.959 \ci{0.884}{1.046} & 1.149 \ci{1.051}{1.282} & 0.821 \ci{0.730}{0.909} & 1.017\\
\oursA vs.\ \oursS & 0.831 \ci{0.744}{0.915} & 0.820 & 1.084 & 0.934\\
\oursF vs.\ \oursS & 0.796 \ci{0.696}{0.903} & 0.942 & 0.890 \ci{0.815}{0.961} & 0.950\\
predictor vs.\ Lu \& Yan & 0.725 \ci{0.637}{0.823} & 0.671 & 1.078 \ci{0.946}{1.212} & 1.002\\
\bottomrule
\end{tabular}
\end{table}

\section{Solver-feedback training of the predictor: negative results}
\label{app:negatives}

We tried to train the predictor itself on the solver's cost. All routes below use the same predictor, cost signal and development data, and none beats continued supervision at matched data and compute. Ratios are instructions or charged time against the stated reference; for routes (ii), (iii), (viii) and (ix) the reference is the predictor in instructions, against which the shift policy reaches 0.927 on the same data. We group the routes into four families.

\textbf{Gradient estimates of the solver cost through the decoder.}
\begin{enumerate}[label=(\roman*),leftmargin=2em]
\item Exact policy gradient through the stochastic crash (a Gumbel-max sequential softmax over edges whose path log-probability is computed in one \texttt{logcumsumexp} over removal steps; unbiasedness verified against the enumerated law). On 14,351 LPs at the 20k size, with 256 LPs $\times$ 8 samples per step and 139k solves, the deterministic start reaches 1.004 to 1.036 and the samples 0.996 to 1.016 relative to the predictor's own start. On 123k small LPs with 2,048 LPs per step, it gains 2\% on held-out small LPs, which does not transfer to the 20k size. The cosine between the summed gradients of two disjoint halves of 256 LPs is 0.00 to 0.12 for the policy gradient at temperatures 0.02 to 0.2, with and without rank shaping, against 0.51 for the cross-entropy; the batch signal-to-noise ratio is 0.97 to 1.25 against 2.9.
\item Per-edge policy gradient with 4,096 LPs per step: 1.012, with a batch signal-to-noise ratio of at most 1 at every batch size and no direction better than a random one.
\item Evolution strategies on a ten-dimensional continuous decoder action: 0.886 at best ($0.955$ to $0.980\times$ the shift policy, every interval containing 1, one training run), with summed instructions at the shift policy's level.
\item Sampled deployment of the stochastic crash (one seeded sample at $\tau = 0.02$ to $0.05$): no gain and added tail risk.
\end{enumerate}

\needspace{5\baselineskip}\textbf{Targets from measured best starts.}
\begin{enumerate}[label=(\roman*),leftmargin=2em,resume]
\item Perturbed Fenchel--Young losses toward measured best starts (exact maximum-weight basis oracle or the crash): best 0.921 against 0.798 for matched supervision on a weaker predictor; harmful as a second stage on the strong predictor (1.11 to 1.17).
\item Teacher-path imitation 1.14 to 1.18; set-valued DAgger 0.95, then degrading; two-round commitment-conditioned decoding $\pm 2\%$.
\end{enumerate}

\textbf{Learned cost models.}
\begin{enumerate}[label=(\roman*),leftmargin=2em,resume]
\item End-to-end fine-tuning of the encoder through the shift policy's loss, starting from a trained policy (equal data and solves): 0.947 to 0.957 against 0.902 for the policy-only control at 50 steps under a cross-entropy anchor, 0.933 at learning rate $5\times10^{-6}$ without one, 0.918 with a trust region on the scores; 0.971 when the policy is trained from scratch.
\item A start-conditioned critic that re-ranks a pool of shifted starts: 0.913 ($0.985\times$ \ci{0.963}{1.007} the shift policy).
\item A critic-gradient adapter on the scores: 0.99 to 1.01, because the adapter exploits the critic's errors.
\item A basis-conditioned GNN verifier over 9 to 33 candidate starts captures 65\% of the oracle gap in instructions within the training classes and $-4\%$ outside them. Solver-in-the-loop racing over candidates is worse than no racing, because signals from a prefix of the solve are uninformative until 20 to 40\% of the path and entering a candidate costs 2\% of a solve.
\end{enumerate}

\textbf{Richer basis models and other variants.}
\begin{enumerate}[label=(\roman*),leftmargin=2em,resume]
\item A determinantal basis model ($P(B) \propto \det([A\ I]_B)^2 \prod_{j \in B} e^{s_j}$, exact Cauchy--Binet likelihood, convex in the scores) and a masked discrete diffusion model over complete bases both learn their objectives (the diffusion model's conditional AUC rises from 0.69 to 0.80 with the revealed fraction) but yield no better starts, because the crash discards most of the value of exactness (an exact optimum through the crash: $0.60\times$; installed directly: $0.05\times$) and model samples land outside the 2 to 3\% basin.
\item A crash that preserves primal feasibility removes primal phase 1 (iterations $0.55$ to $0.69\times$ the native solver) and takes $1.01$ to $1.03\times$ the native solver in charged time, because the remaining pivots run on a dense basis (FTRAN 45 to 329\,$\mu$s per call).
\item Random node features for symmetry breaking: $-0.002$ overlap.
\item Algorithm-aware continuation (primal, dual or the solver's choice from the learned start): the oracle over all four methods lifts the ceiling of the two-way harm gate by 3.5\% at the 20k size.
\end{enumerate}

\section{Guardrails}
\label{app:gates}

This appendix details the three guardrails of \S\ref{sec:safety}: the novelty gate and the harm gate, which decide before the solve, and the optimality check (Appendix~\ref{app:check}), which acts after it. Table~\ref{tab:hyper} lists their settings.

\textbf{Novelty gate.} The policy is learned within the training classes; outside them it chooses harmful actions (leave-one-class-out, 1.14 to 1.17 times the predictor's time for the shift policy and 1.20 for \oursA), and thresholding its confidence does not help, because the ensemble assigns high probability to harmful actions. The nearest-neighbor test on the standardized $[h, z]$ vector abstains on 98 to 100\% of held-out-class LPs and on 2\% of in-class ones, brings the held-out-class ratio to $1.000$ of the predictor, and costs 2\,ms. Abstaining to the native solver instead of the unshifted start would be worse on the training classes held out one at a time ($1.45\times$); \S\ref{sec:generalization} shows that this choice does not carry over to new classes.

\textbf{Harm gate.} A gradient-boosted regressor predicts $\log(I_{\text{learned}}/I_{\text{native}})$ from 73 deployment-time features: the LP's size and sparsity statistics, statistics of the predicted scores, and the composition of the decoded start. It uses no class labels. When the predicted saving is below 12\%, the native solver runs instead. The model and threshold were chosen by leave-one-class-out cross-validation on the training classes. To measure how much of the gate's value comes from recognizing the training classes, we also evaluate a class-blind variant: for each test LP, a gate refit with that LP's class removed from its training data.

\needspace{5\baselineskip}\textbf{Gate details.}
\begin{itemize}
\item Novelty gate: standardized 84-dimensional $[h, z]$ against the 7,699 training rows, threshold $2.466$, the 0.99 quantile of nearest-neighbor distances from a random fifth of the training rows to the rest. Every quantile in $[0.5, 0.99]$ gives leave-one-class-out $1.000$, and 0.99 was chosen for the smallest in-class abstention (2\%).
\item Harm gate: models compared by leave-one-class-out cross-validation on 3,466 training rows (ridge, cost-weighted logistic regression, and boosted trees with squared or Huber loss); boosted trees with Huber loss were chosen (150 trees of depth 3, learning rate 0.05, leaves of at least 25 rows). The threshold is chosen from 201 quantiles of the out-of-fold predictions to minimize the class-balanced time-weighted loss: keep the learned start if and only if $\hat y \le -0.128$. Its tree evaluation costs 3.3\,ms per LP and is charged.
\item On \Dtw the harm gate turns bin packing from $1.75\times$ the native solver to $1.05\times$ and hub location from $1.98\times$ to $0.99\times$. Fit without the LP's class, it loses most of its gain (Table~\ref{tab:main}), so its skill within the training classes is partly recognition of those classes, which is why the novelty gate is reported beside it.
\item A budgeted restart (abort the learned solve after $4m$ iterations and run the native solver) removes one stall on the development panel and is available as a tail guard.
\end{itemize}

\subsection{Optimality check}
\label{app:check}

\textbf{The failure.} On one badly scaled class, load balancing in its primal form, a learned start that is not dual feasible can end at a basis that HiGHS accepts as optimal but whose objective is up to 20\% too high. These LPs minimize one makespan variable subject to 47k inequality rows with coefficients up to $10^4$; the optimum is about $0.003$ and every optimal dual is of order $10^{-7}$. HiGHS returns a basis that is optimal within its absolute dual feasibility tolerance of $10^{-7}$: about 8,800 reduced costs have the wrong sign, each below the tolerance, and their sum is about $4\times10^{-4}$. On this class the tolerance is thus of the size of the optimal duals themselves. The likely mechanism is that HiGHS starts the dual simplex in phase 2 when the dual infeasibilities of the start are small, absorbs them by shifting costs, and after removing the shifts finds no wrong-signed reduced cost above tolerance, so it triggers no cleanup. Native starts, which are dual feasible on this class, never failed this way, so the check is charged to the learned start.

\textbf{The check.} HiGHS' sum of dual infeasibilities, unlike its count, includes values below tolerance, which is why it detects these failures. If the sum exceeds $10^{-6}(1 + |f|)$ for objective value $f$, the solve continues from the returned basis with the primal and dual feasibility tolerances tightened to $10^{-9}$; if that does not reach optimality, the native solver solves the LP afresh. The threshold compares a sum of reduced costs with an objective value, so it is an empirical rule and not a certificate; it rests on a separation of more than three orders of magnitude on development data. Scaling the objective of this class would be a solver-side remedy; the check needs no class-specific setting.

\textbf{Development evaluation.} On a development panel of 477 LPs at 40k to 60k variables (including 23 duplicates of training LPs), the solver reported 2,370 solves as optimal. Twelve of them were wrong, every one a learned start on a primal-form load-balancing LP, as were the wrong-objective terminations on \Tfif. Candidate checks on these 2,370 solves:
\begin{itemize}
\item HiGHS' sum of dual infeasibilities with threshold $10^{-6}(1 + |f|)$ flags 12 of the 12 wrong solves and 0 of the 2,358 correct ones. The statistic is at least $1.9\times10^{-4}$ on every wrong solve and at most $5.8\times10^{-8}$ on every correct one, so any threshold from $10^{-7}$ to $10^{-5}$ works. It reads one solver statistic (8\,$\mu$s).
\item A zero-iteration re-check at $10^{-9}$ tolerances flags the 12 wrong solves and 97 correct ones, and continuing from those false positives breaks four correct solves.
\item An explicit KKT residual takes 73\,ms and flags 137 correct solves.
\end{itemize}
Continuation at $10^{-9}$ from the returned basis recovered the reference optimum on all 12 wrong solves in 1.7 to 67\,s, in every case less than a native solve (77 to 144\,s); the fallback to a native solve never triggered. On W2 of \Dfif the check reduces the predictor's failures from 6 to 1 (a time limit) and its summed time ratio from 1.385 to 0.863, while the class-balanced geometric mean moves from 0.619 to 0.582. The rule was chosen on development data. The \Tfif and \Tei columns of Table~\ref{tab:main} are reported without the check (the failures of \oursS there are the load-balancing wrong optima plus one solver error on \Tei); the \Ttw panel and Table~\ref{tab:checked} apply it inline.

\section{Deployment procedure and method details}
\label{app:deploy}

Algorithm~\ref{alg:deploy} (\S\ref{sec:safety}) gives the per-LP procedure of \ours at deployment. Every step, including the ones whose result is discarded, is included in the charged time. \oursS uses the action set $\mathcal{A}_1$ and always evaluates the harm gate; \oursA and \oursF use $\mathcal{A}_2$ and $\mathcal{A}_3$ and are reported with and without the harm gate. Table~\ref{tab:hyper} collects the settings of every learned component and of the solver.

\begin{table}[h]
\centering\footnotesize
\caption{Settings of the learned components and of the solver. Selection was always on development data or a held-out fifth of the training LPs, never on test data.}
\label{tab:hyper}
\setlength{\tabcolsep}{4pt}
\begin{tabular}{@{}>{\raggedright\arraybackslash}p{3.6cm}>{\raggedright\arraybackslash}p{10.0cm}@{}}
\toprule
Component & Setting\\
\midrule
Predictor & 8-layer GATv2 on the variable-constraint graph, 1.8M parameters, 72-dimensional pooled embedding; fused kernel with TF32, one LP per forward pass (6.4\,ms median at 20k variables on an H100)\\
\quad stage 1 & 114,874 LPs of \Trs, 12 epochs (43k steps, 32 LPs per step, edge budget 400k), AdamW at $10^{-3}$ without weight decay, warm-up and cosine decay, gradient clip 5; loss: basis-membership cross-entropy divided by $m$ plus $0.2\times$ auxiliary dual-tightness and bound-status terms\\
\quad stage 2 & 3,457 LPs of \Trtw, 6 passes, 48 LPs per step, AdamW at $3\times10^{-4}$, clip 1; checkpoint every 37 steps, chosen by paired solve time on the 283-LP selection split of \Dtw\\
Decoder & entry threshold $10^{-8}$ of the column maximum, degree exponent $p = 0.5$, shifts $\{0, \pm0.5, \pm1, \pm2, \pm3\}$\\
Policies & MLP $84 \to 128 \to 128 \to |\mathcal{A}|$, SiLU, dropout 0.1, AdamW (weight decay $10^{-4}$), batch 256 LPs, loss clip $\pm 3$, failures at $+3$, last layer initialized to the unshifted start\\
\quad \oursS & $|\mathcal{A}_1| = 9$; 7,699 training LPs; learning rate $10^{-3}$, 100 epochs, 3 seeds\\
\quad \oursA, \oursF & $|\mathcal{A}_2| = 11$, $|\mathcal{A}_3| = 20$; 9,349 and 9,197 training LPs; learning rate $3\times10^{-3}$, 150 epochs, 10 seeds\\
Novelty gate & standardized 84-dimensional $[h, z]$; threshold 2.466, the 0.99 quantile of nearest-neighbor distances from a random fifth of training LPs to the rest\\
Harm gate & gradient-boosted regression trees on 73 features, Huber loss, 150 trees of depth 3, learning rate 0.05, leaves of at least 25 rows; keep the learned start if the predicted log ratio to the native solver is at most $-0.128$\\
Optimality check & sum of dual infeasibilities $> 10^{-6}(1 + |f|)$; continuation from the returned basis at primal and dual feasibility tolerances $10^{-9}$, then the native solver if needed\\
Solver & HiGHS 1.13.1, one thread, presolve off, all other options at their defaults (primal and dual feasibility tolerances $10^{-7}$) unless listed here. Learned starts: dual simplex (\texttt{simplex\_strategy} 1) with HiGHS' default dual edge weights, or Devex forced (\texttt{simplex\_dual\_edge\_weight\_strategy} 1) for the Devex actions; primal simplex for the primal action. Native solver: \texttt{simplex\_strategy} 0 (choose)\\
\bottomrule
\end{tabular}
\end{table}

\begingroup

\sisetup{detect-all, table-number-alignment = center}

\section{The SyntheticLPs generator library}
\label{app:Syn_generator}

Empirical work on linear programming is limited by the supply of instances. The public collections, such as Netlib and MIPLIB, are small and fixed. Netlib contains roughly one hundred usable LPs: enough for a benchmark table, but not enough to train a model, to measure how a method's behavior changes with problem size, or to hold out a test set that is disjoint from every instance used during development. Random-matrix generators, which sample $A$, $b$ and $c$ from some distribution, give an unlimited supply of instances, but a dense Gaussian $A$ has none of the sparsity, block structure or network incidence patterns of practical LPs, and conclusions drawn from such instances transfer poorly.

All LPs in this paper come from \textsc{SyntheticLPs.jl}\footnote{\url{https://github.com/flixpar/SyntheticLPs.jl}}, which samples problems rather than matrices: a transportation plan, a refinery blend, a radiotherapy fluence map, a crew roster. Each generator implements a standard formulation with the constraint structure of its domain, parameterized so that its data, and with them its dimensions, are redrawn for every instance. The matrix structure comes from the model, and only the numbers are sampled. This appendix describes the design of the library (\S\ref{sec:gen-design}), the instance space and generation cost (\S\ref{sec:gen-space}), reproducibility and feasibility control (\S\ref{sec:gen-repro}), a structural comparison with Netlib (\S\ref{sec:gen-struct}) and the limitations of the generated instances (\S\ref{sec:gen-limits}), and lists every generator (\S\ref{app:gen-catalogue}). All measurements in this appendix use HiGHS 1.13.1 through JuMP, with presolve off for every solve-based statistic, over \num{2364} generated instances.

\subsection{Design}
\label{sec:gen-design}

The library has two levels. A \emph{category} is a problem domain, such as transportation, radiotherapy or hub location; we call categories classes in the main text. A \emph{formulation} (called a variant in the library's code) is a concrete model within a category. Formulations are distinct models rather than parameter settings of a common model: within the transportation category, for example, the capacitated and fixed-charge formulations have different rows, not only different coefficients. The version used here provides 45 categories with 129 formulations, of which 71 contain integer or binary variables and 58 are pure LPs. Each formulation can be emitted in a primal or a dual form. Because the dual form has a transposed constraint matrix and looks very different to the solver, the main text counts the two forms of a formulation as separate variants. In this paper, formulations with integer or binary variables are used through their LP relaxations.

Every formulation follows a two-stage contract. A constructor \texttt{Variant(}\textit{target\_variables}, \textit{status}, \textit{seed}\texttt{)} performs all random sampling, drawing from a constructor-local random number generator seeded by \textit{seed}, and stores the sampled data in a plain struct. A second function, \texttt{build\_model}, maps that struct to a JuMP model and makes no random draws. This has three consequences:
\begin{itemize}
\item an instance is determined entirely by the tuple $(\textit{variant}, \textit{target\_variables}, \textit{status}, \textit{seed})$;
\item generation neither reads nor advances the caller's global random state, so instances can be generated in parallel without coordination and without perturbing other random components of a pipeline;
\item the sampled data can be inspected as a typed object before any model is built; a transportation instance, for example, exposes its supplies, demands and costs, which an MPS file does not retain.
\end{itemize}

\subsection{Instance space and generation cost}
\label{sec:gen-space}

Each of the 129 formulations can be combined with several independent options (Table~\ref{tab:gen-axes}), so the number of distinct models is effectively unbounded. Changing the seed changes the data. Changing \texttt{target\_variables} re-derives the underlying dimensions, such as the numbers of sources and destinations, periods and machines, or beams and voxels, rather than padding a fixed skeleton. Dualization transposes the whole structure. In this paper we use the primal and dual forms of each formulation and the LP relaxations of the formulations with integer variables.

\begin{table}[t]
  \centering
  \caption{Orthogonal axes of the generated instance space. The final three are
  post-\texttt{build\_model} transforms applied centrally, so they act uniformly
  across all formulations rather than being reimplemented per generator.}
  \label{tab:gen-axes}
  \small
  \begin{tabular}{@{}llr@{}}
    \toprule
    Axis & Effect on the emitted model & Cardinality \\
    \midrule
    \texttt{variant}               & the formulation itself                        & \num{129} \\
    \texttt{target\_variables}     & dimensions re-derived, not padded             & wide range \\
    \texttt{seed}                  & fresh data at fixed size and structure        & unbounded \\
    \texttt{feasibility\_status} & feasible / infeasible / unknown               & 3 \\
    \texttt{relax\_integer}        & a MIP, or its LP relaxation                   & 2 \\
    \texttt{bounds\_\allowbreak to\_\allowbreak constraints} & bounds kept, or emitted as explicit rows    & 2 \\
    \texttt{dualize}               & the primal, or its LP dual                    & 2 \\
    \bottomrule
  \end{tabular}
\end{table}

The \texttt{target\_variables} argument is a request, and the generator chooses internal dimensions to meet it. Table~\ref{tab:gen-scaling} shows how closely it does so across 2.6 orders of magnitude, with one instance per formulation at each target size. For most formulations the requested count is met exactly, because the generator solves for dimensions rather than approximating them. Build time is essentially independent of size, at roughly 19 instances per second on a single thread (about 1.6 million per core-day), because the work is proportional to the number of nonzeros and sparse models stay sparse as they grow. The exceptions are a few formulations that are dense by necessity, namely the regression family and both portfolio formulations, whose number of nonzeros grows as $O(mn)$ rather than $O(n)$. These take 10 to 13\,s per build at a few thousand variables and become impractical beyond that; they are among the formulations not measured at the larger targets in Table~\ref{tab:gen-scaling}.

\begin{table}[t]
  \centering
  \caption{Size fidelity and generation cost (\num{558} builds). The median
  ratio of achieved to requested variable count is exactly $1.000$ at every
  scale. Median build time is flat at ${\approx}\SI{52}{\milli\second}$ from
  \num{200} variables upward; the build time at \num{50} variables is omitted because
  it absorbs per-formulation JIT compilation on first call. Fewer formulations were
  measured at the largest targets (see text).}
  \label{tab:gen-scaling}
  \small
  \begin{tabular}{@{}r r S[table-format=1.3] r r S[table-format=1.3]@{}}
    \toprule
    Requested & Formulations & {Median ratio} & Within $\pm10\%$ & Within $\pm25\%$ & {Median build (s)} \\
    \midrule
    \num{50}     & 129 & 1.000 & 91\,\% & 99\,\%  & {--} \\
    \num{200}    & 129 & 1.000 & 96\,\% & 100\,\% & 0.052 \\
    \num{1000}   & 129 & 1.000 & 97\,\% & 98\,\%  & 0.052 \\
    \num{5000}   &  93 & 1.000 & 96\,\% & 97\,\%  & 0.052 \\
    \num{20000}  &  78 & 1.000 & 90\,\% & 91\,\%  & 0.053 \\
    \bottomrule
  \end{tabular}
\end{table}

\subsection{Reproducibility and feasibility control}
\label{sec:gen-repro}

\textbf{Determinism.} We serialized generated models to MPS and compared their SHA-256 digests. For all 129 formulations, identical inputs give byte-identical models, a different seed gives a different model, and the output neither depends on nor advances the caller's global random state. A dataset is therefore fully described by its seeds and configuration.

\textbf{Feasibility.} Each generator accepts the status \texttt{feasible}, \texttt{infeasible} or \texttt{unknown} and meets the request by construction, adjusting capacities, demands or budgets without calling a solver. We tested every formulation with three seeds and both definite statuses, solving each instance with presolve off: all 387 \texttt{feasible} requests ended optimal and all 387 \texttt{infeasible} requests ended infeasible. Many generators also attach a complete primal witness for feasible requests or a structural certificate of infeasibility that can be checked without a solver. Where a generator's feasibility logic is heuristic, an optional solver check rebuilds the instance with the next seed on a mismatch, walking through seeds deterministically so that a given seed and solver always give the same instance. The \texttt{unknown} status gives no guarantee; roughly a quarter of the instances generated with it are infeasible.

\subsection{Structural comparison with Netlib}
\label{sec:gen-struct}

To check that the generated instances resemble real LPs, we compared 283 synthetic instances at 400 variables with 109 Netlib instances, using the same extraction of $A$, $b$ and $c$, the same solver and presolve off (Table~\ref{tab:gen-netlib}). Sparsity, which naive generators reproduce worst, falls in the right range: density, nonzeros per column and per row, and the variability of column degrees are all within a factor of 1.7 of Netlib's medians. The incidence structure is inherited from the formulations and so needs no fitting.

\begin{table}[t]
  \centering
  \caption{Structural comparison against Netlib (medians). Sparsity statistics
  agree closely; numerical range and degeneracy remain lower than in real
  instances.}
  \label{tab:gen-netlib}
  \small
  \begin{tabular}{@{}l S[table-format=1.4] S[table-format=2.4] l@{}}
    \toprule
    Property (median) & {Netlib} & {Synthetic} & Assessment \\
    \midrule
    Matrix density                          & 0.0079 & 0.0132 & close \\
    Nonzeros per column                     & 4.23   & 3.71   & close \\
    Nonzeros per row                        & 8.45   & 5.14   & close \\
    Column-degree variability (CV)          & 0.55   & 0.43   & close \\
    \addlinespace
    Primal degeneracy                       & 0.196  & 0.094  & about half \\
    Simplex iterations per row              & 1.24   & 0.55   & less work \\
    Coefficient spread, $\log_{10}|A_{ij}|$ & 4.12   & 1.85   & narrower \\
    Variables absent from the objective     & 0.382  & 0.000  & gap \\
    \bottomrule
  \end{tabular}
\end{table}

Two gaps remain. First, synthetic coefficients span 1.85 orders of magnitude against Netlib's 4.12, so the generated instances exercise numerical robustness less. Second, in most formulations every variable has a nonzero objective coefficient, whereas real models contain auxiliary and balance columns that are absent from the objective (38\% of columns at the Netlib median). Both gaps come from drawing parameters independently from smooth continuous distributions, and both could be addressed in the sampling layer without changing any formulation. On the other hand, the generated instances span a wider range than Netlib on several features (Table~\ref{tab:gen-breadth}). They range from pure network matrices, whose entries are all $\pm 1$ and whose constraint matrix is totally unimodular, to fully dense statistical LPs, and from completely nondegenerate instances to ones that are 80\% degenerate.

\begin{table}[t]
  \centering
  \caption{Feature-space breadth, reported as the 5th--95th percentile interval
  over each corpus.}
  \label{tab:gen-breadth}
  \small
  \begin{tabular}{@{}lll@{}}
    \toprule
    Feature & Synthetic & Netlib \\
    \midrule
    Matrix density                          & 0.006--1.000 & 0.000--0.092 \\
    Nonzeros per column                     & 1.96--55.98  & 2.18--14.24 \\
    Column-degree variability (CV)          & 0.00--2.40   & 0.11--1.23 \\
    Coefficient spread, $\log_{10}|A_{ij}|$ & 0.00--5.75   & 0.00--8.65 \\
    Fraction of entries equal to $\pm 1$    & 0.000--1.000 & 0.094--1.000 \\
    Primal degeneracy                       & 0.000--0.800 & 0.000--0.841 \\
    \bottomrule
  \end{tabular}
\end{table}

\subsection{Limitations}
\label{sec:gen-limits}

The following limitations bear on how results on these instances should be read.
\begin{itemize}
\item The numerical range is narrower than in real LPs (1.85 against 4.12 decades), so results that depend on numerical conditioning should not treat these instances as representative.
\item Primal degeneracy is about half of Netlib's. It is present where the combinatorics force it, as in assignment, covering and routing polytopes, and largely absent in continuous allocation models, where real data would introduce it through repeated and rounded values.
\item Dense formulations do not scale past a few thousand variables (\S\ref{sec:gen-space}).
\item Multi-period data are drawn with noise that is independent over time, so constraints that couple periods, such as ramping, storage and inventory balance, bind less tightly than under the persistent deviations of real demand.
\item Netlib is an imperfect reference: its instances are larger than those measured here (median 1,376 against 400 columns), so part of the difference in iterations per row reflects size rather than realism.
\end{itemize}

\subsection{Full generator catalog}
\label{app:gen-catalogue}

Table~\ref{tab:gen-catalogue} lists all 45 categories and 129 formulations of the version used here, and how each category is used in this paper. ``Train, test'': one of the 25 classes of \Trtw with test instances in \Ttw. ``Train'': in \Trtw, but without test instances. ``Unseen'': one of the six classes of \Tnc, absent from all training data (\Tad uses five of them). At the 20k size, the 26 training categories contribute 67 of their formulations; the others do not produce LPs of that size. The predictor's first training stage on \Trs uses 37 categories.

\begingroup
\small
\setlength{\tabcolsep}{4pt}
\begin{longtable}{@{}l r l >{\raggedright\arraybackslash}p{0.48\linewidth}@{}}
  \caption{The categories and formulations of SyntheticLPs.jl and their use in this paper.}
  \label{tab:gen-catalogue}\\
  \toprule
  Category & \# & Use & Formulations \\
  \midrule
  \endfirsthead
  \toprule
  Category & \# & Use & Formulations \\
  \midrule
  \endhead
  \bottomrule
  \endfoot
  \texttt{airline\_crew} & 1 &  & standard \\
  \texttt{assignment} & 2 &  & standard, workload\_balance \\
  \texttt{bin\_packing} & 2 & train, test & heterogeneous, standard \\
  \texttt{blending} & 3 & unseen & equipment\_batches, multi\_product, standard \\
  \texttt{container\_loading} & 2 & unseen & standard, two\_dimensional\_bin\_packing \\
  \texttt{crop\_planning} & 1 & train, test & standard \\
  \texttt{cutting\_stock} & 4 &  & due\_dates, integer\_patterns, setup\_cost, standard \\
  \texttt{diet\_problem} & 3 & unseen & food\_groups, nutrient\_bounds, standard \\
  \texttt{energy} & 7 & train, test & dc\_opf, optimal\_transmission\_switching, ramping, reserves, standard, storage, transmission \\
  \texttt{facility\_location} & 3 & train, test & p\_median, standard, two\_echelon \\
  \texttt{feed\_blending} & 1 &  & standard \\
  \texttt{graph\_optimization} & 6 & train, test & generalized\_independent\_set, independent\_set, map\_labeling, quasi\_clique, vertex\_coloring, vertex\_cover \\
  \texttt{hub\_location} & 8 & train, test & budgeted\_backbone, capacitated, compact\_single\_allocation, hub\_covering, hub\_network, multiple\_allocation, p\_hub\_median, r\_allocation \\
  \texttt{inventory} & 4 & train, test & lot\_sizing, multi\_echelon, multi\_item, standard \\
  \texttt{inverse\_optimization} & 8 & train, test & classical\_normalized, linf, market\_clearing, noisy\_observations, restricted\_optimal\_value, shortest\_path, shortest\_path\_layered, standard \\
  \texttt{job\_shop\_scheduling} & 1 & unseen & standard \\
  \texttt{knapsack} & 4 & train, test & bounded, mixed\_integer\_set, multidimensional, standard \\
  \texttt{land\_use} & 1 & train, test & standard \\
  \texttt{load\_balancing} & 2 & train, test & discrete\_placement, standard \\
  \texttt{maritime\_inventory\_routing} & 1 & train, test & standard \\
  \texttt{multi\_commodity\_flow} & 3 & train, test & binary\_capacity, integer\_flow, standard \\
  \texttt{network\_flow} & 2 &  & generalized\_flow, standard \\
  \texttt{neural\_network\_verification} & 1 &  & relu\_big\_m \\
  \texttt{nurse\_scheduling} & 1 & train, test & standard \\
  \texttt{operating\_room\_scheduling} & 6 & train, test & benchmark\_loading, case\_sequencing, elective\_assignment, master\_surgical\_schedule, robust\_elective, weekly\_planning \\
  \texttt{portfolio} & 2 &  & cvar, tracking\_error \\
  \texttt{process\_planning} & 5 & unseen & campaign, capacity\_expansion, hydrogen\_network, mode\_switching, refinery \\
  \texttt{product\_mix} & 1 &  & standard \\
  \texttt{production\_planning} & 1 &  & standard \\
  \texttt{project\_selection} & 1 & train & standard \\
  \texttt{radiotherapy} & 5 &  & beam\_angle\_selection, mean\_tail\_dose, minmax\_deviation, robust\_fluence, weighted\_deviation \\
  \texttt{regression} & 4 &  & basis\_pursuit, chebyshev, lad, quantile \\
  \texttt{resilient\_network\_design} & 1 & unseen & standard \\
  \texttt{resource\_allocation} & 1 &  & standard \\
  \texttt{revenue\_management} & 2 & train, test & standard, stochastic\_overbooking \\
  \texttt{scheduling} & 1 & train, test & standard \\
  \texttt{set\_system} & 4 & train, test & combinatorial\_auction, set\_cover, set\_packing, set\_partitioning \\
  \texttt{stochastic\_program} & 1 & train, test & standard \\
  \texttt{supply\_chain} & 5 & train, test & carbon, multi\_product, network\_planning, single\_source, standard \\
  \texttt{telecom\_network\_design} & 1 & train, test & standard \\
  \texttt{transportation} & 6 & train, test & balanced, capacitated, emission\_constrained, fixed\_charge, standard, transshipment \\
  \texttt{tsp} & 8 & train, test & assignment\_relaxation, asymmetric, flow, multiple\_salespersons, precedence, prize\_collecting, standard, time\_windows \\
  \texttt{unit\_commitment} & 1 & train, test & standard \\
  \texttt{vehicle\_routing} & 1 & train, test & cvrp \\
  \texttt{workforce\_shift\_scheduling} & 1 &  & covering \\
\end{longtable}
\endgroup
\endgroup

\section{Data, test sets and protocol}
\label{app:data}

This appendix lists the datasets (Table~\ref{tab:data}) and describes how the splits are kept disjoint, the timing protocol and the machines.

\begin{table}[h]
\centering\scriptsize
\caption{Datasets. Sizes are requested numbers of variables. Every test set is disjoint from every training and development set by content hash, generator call and seed; the classes of \Tnc and \Tad were never used in training.}
\label{tab:data}
\setlength{\tabcolsep}{3.5pt}
\begin{tabular}{@{}llrl>{\raggedright\arraybackslash}p{4.3cm}l}
\toprule
Set & Variables & LPs & Classes & Use\\
\midrule
\Trs & 50--2,000 & 114,874 & 37 & predictor, stage 1 \\
\Trtw & 15k--25k & 9,511 & 26 (25 with test LPs) & predictor stage 2 (3,457 LPs); policies and gates \\
\Dtw & 15k--25k & 450 & 25 & selection and ablations; half twinned with training \\
\quad selection split & 15k--25k & 283 of \Dtw & 22 & checkpoints, epochs and designs \\
\Dfif & 40k--60k & 454 & 23 & optimality check; policies at larger sizes \\
\Ttw & 15k--25k & 1,356 & 25 & main result \\
\Tfif & 40k--60k & 382 & 25 & larger instances \\
\Tei & 60k--100k & 81 & 24 & larger still \\
\Tva & 15k--25k & 296 of \Ttw & 10 formulations, 6 dual forms & held-out formulations and dual forms \\
\Tnc & 15k--25k & 234 & 6 unseen & transfer to new classes \\
\Tad & 15k--25k & 616 + 125 & 5 unseen & few-shot adaptation, online bandit \\
\bottomrule
\end{tabular}
\end{table}

\textbf{Tiers and test sets.} LPs are generated in tiers named by the requested number of variables: small (50 to 2,000), 20k (15,000 to 25,000) and 50k/80k (40,000 to 100,000). One fresh set at the 20k size is split by rule into \Ttw (the 25 training classes with test instances) and \Tnc (six classes no model saw). Two fresh sets at larger sizes give \Tfif and \Tei; a fixed rule excludes three formulations whose instances do not grow with the requested size. \Tad is drawn from an 800-LP panel of five held-out classes, split by a fixed rule into 125 adaptation and 616 evaluation LPs; evaluation LPs that share a generator instance with an adaptation LP are excluded. \Dtw is used for selection and ablations; about half of it has a primal or dual twin in the predictor's training data, so its levels against the native solver are optimistic, and we use it for paired comparisons among learned methods. Its 283-LP selection split (the W1 LPs of a fixed validation split, 22 classes) selected checkpoints, epochs and designs.

\textbf{Timing protocol.} Every method solves the same stored LP in a fresh solver instance on one pinned core after a discarded warm-up solve, in a seeded random order per LP. Panels at the 20k size were run twice (\Tnc three times) and averaged within LP; panels at larger sizes ran once on exclusive nodes. Instruction counts are recorded beside every time. A solve fails if it hits the time limit, returns an error, or ends with an objective more than $10^{-5}$ relative from a reference solve; failures stay in the data at twice the time limit, as in PAR2. The class-balanced geometric mean averages the per-LP log ratios within each class and then across classes; intervals come from 2,000 two-level bootstrap draws (classes, then generator instances within class, so that the primal and dual forms of one instance are resampled together). Because a geometric mean can hide a few very slow LPs, we also report summed time, classes won and lost, and failures.

\textbf{Generators and admission.} SyntheticLPs.jl provides 45 problem classes with 129 formulations in all (Appendix~\ref{app:gen-catalogue}), each in a primal and a dual form; we count the two forms of a formulation as separate variants. The 20k tier covers 26 classes and 67 formulations, or 134 variants, that produce LPs of that size (project selection has no test instances). Formulations with integer or binary variables enter through their LP relaxations, so every instance is an LP. Presolve runs once offline with HiGHS; the stored LP is what every method solves. Timed sets admit an LP only if the all-slack dual simplex takes at least 2,000 iterations or 1\,s on it, and the same solve's time assigns the windows W1 to W3.

\textbf{Disjoint splits.} Every generated LP carries a name-free content hash. Test LPs are generated from fresh seeds and admitted only if no training or development LP shares their content hash, generator call or seed, so no test LP has an exact copy or a primal-dual twin in the training or development data. The training data of the policies and gates (\Trtw, 9,511 LPs) likewise contain no copy or twin of a development LP. The predictor's 20k training data, however, contain transposed twins of 179 of the 381 W1 LPs of \Dtw. Levels of \Dtw against the native solver are therefore optimistic, and there we rely on relative comparisons among learned methods trained on the same data. Nothing was selected on any test set.

\textbf{Machines.} Panels ran on three kinds of machines: a 112-core Xeon 8480+ server for some development panels, a 128-core Xeon 8592+ node with four H100s, and exclusive Xeon 8168, 8268 and 6248R nodes of a CPU partition for the solves at 40k to 100k variables, with starts and gate features prepared on the GPU node and charged from there. Ratios are compared only within a panel, where every method of an LP ran on one machine; the CPU model is recorded per row.

\section{Full results, label ablation, lookup control, seeds and deployment cost}
\label{app:results}

\begin{figure}[h]
\centering
\includegraphics[width=\linewidth]{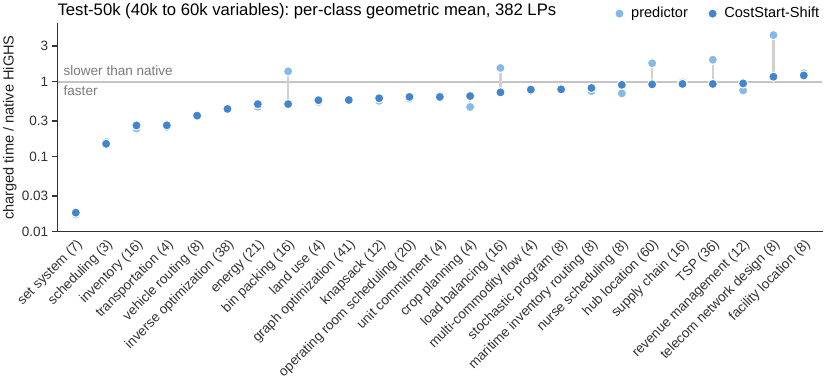}
\caption{Per-class charged time on \Tfif (ratio to the native solver, geometric mean over the class's LPs, all three windows; LP counts in parentheses). The predictor loses six classes; the shift policy and harm gate of \oursS turn four of them into wins, and the two remaining losses are classes where the primal simplex from a primal-feasible all-slack basis is hard to beat. Table~\ref{tab:percat-test50k} gives the same panel with more methods.}
\label{fig:percat}
\end{figure}

This appendix extends the results of \S\ref{sec:headline} and \S\ref{sec:periter}. Table~\ref{tab:mainfull} is the full version of Table~\ref{tab:main}. On \Ttw, Table~\ref{tab:windows-test20k} breaks the results down by window, Table~\ref{tab:percat-test20k} by class, Table~\ref{tab:actions-test20k} lists the actions taken and Table~\ref{tab:overhead} the charged overhead; Table~\ref{tab:percat-test50k} and Figure~\ref{fig:percat} give \Tfif by class. Appendices~\ref{app:objective}, \ref{app:lookup} and~\ref{app:seeds} give the training-label ablation, the per-class lookup control and the seed study, and Appendix~\ref{app:cost} the cost of training and of CPU-only inference.

\begin{table}[h]
\centering\scriptsize
\caption{Full version of Table~\ref{tab:main}, with the ablation rows of \oursS and the supervised control (the predictor trained for one more pass over the 20k training data). Iterations relative to the all-slack dual simplex, charged time relative to the native solver; same conventions as Table~\ref{tab:main}.}
\label{tab:mainfull}
\setlength{\tabcolsep}{3pt}
\begin{tabular}{@{}l ccc ccc c@{}}
\toprule
& \multicolumn{3}{c}{\Ttw (1,356 LPs)} & \multicolumn{3}{c}{\Tfif (382 LPs)} & \Tei\\
\cmidrule(lr){2-4}\cmidrule(lr){5-7}\cmidrule(l){8-8}
Method & iter. & time [95\% CI] & W/L & iter. & time [95\% CI] & W/L & time (81)\\
\midrule
all-slack dual simplex & 1.00 & 1.07 \ci{0.87}{1.30} & 12/13 & 1.00 & 0.99 \ci{0.76}{1.30} & 12/13 & 0.96\\
stock solver, presolve on & 0.85 & 1.10 \ci{1.06}{1.16} & 3/22 & -- & -- & -- & --\\
classical crash & 1.22 & 1.76 \ci{1.19}{2.60} & 9/16 & -- & -- & -- & --\\
Fan et al.\ (original budget) & 0.50 & 1.10 \ci{0.89}{1.38} & 15/10 & 0.52 & 1.00 \ci{0.78}{1.28} & 15/10 & 0.96\\
Fan et al.\ (full budget) & 0.49 & 1.10 \ci{0.89}{1.35} & 15/10 & 0.51 & 1.02 \ci{0.78}{1.37} & 15/10 & 0.96\\
Lu \& Yan (original budget) & 0.47 & 0.97 \ci{0.72}{1.29} & 14/11 & 0.46 & 0.83 \ci{0.52}{1.18} & 18/7 & 0.72\\
Lu \& Yan (full budget) & 0.47 & 0.94 \ci{0.70}{1.22} & 14/11 & 0.50 & 0.85 \ci{0.56}{1.23} & 19/6 & 0.72\\
\midrule
predictor & 0.32 & 0.68 \ci{0.50}{0.91} & 20/5 & 0.33 & 0.61 \ci{0.39}{0.89} & 19/6 & 0.60\\
\quad supervised control & 0.32 & 0.69 \ci{0.49}{0.93} & 20/5 & -- & -- & -- & --\\
predictor + shift policy & 0.31 & 0.64 \ci{0.48}{0.86} & 21/4 & 0.34 & 0.60 \ci{0.37}{0.89} & 19/6 & 0.60\\
predictor + harm gate & 0.30 & 0.61 \ci{0.46}{0.76} & 23/2 & 0.31 & 0.52 \ci{0.35}{0.70} & 23/2 & 0.51\\
\oursS (shift policy + harm gate) & 0.30 & 0.61 \ci{0.46}{0.76} & 23/2 & 0.32 & 0.53 \ci{0.34}{0.72} & 23/2 & 0.52\\
\quad with class-blind harm gate & 0.32 & 0.68 \ci{0.52}{0.85} & 18/7 & -- & 0.58 \ci{0.40}{0.79} & 18/7 & 0.58\\
\oursA & 0.25 & 0.51 \ci{0.38}{0.65} & 23/2 & \multicolumn{4}{c}{Table~\ref{tab:checked}}\\
\quad with harm gate & 0.26 & 0.54 \ci{0.41}{0.69} & 23/2 & \multicolumn{4}{c}{Table~\ref{tab:checked}}\\
\oursF & 0.28 & 0.49 \ci{0.37}{0.63} & 24/1 & \multicolumn{4}{c}{Table~\ref{tab:checked}}\\
\quad with harm gate & 0.29 & 0.52 \ci{0.39}{0.65} & 23/2 & \multicolumn{4}{c}{Table~\ref{tab:checked}}\\
\midrule
\multicolumn{8}{@{}l}{Controls for the pricing action}\\
predictor, Devex on every LP & 0.39 & 0.73 \ci{0.56}{0.97} & 18/7 & \multicolumn{4}{c}{--}\\
selector without learned start & 0.88 & 0.86 \ci{0.75}{0.96} & 16/9 & \multicolumn{4}{c}{--}\\
\bottomrule
\end{tabular}
\end{table}

\textbf{Robustness of the main result.} On \Ttw, leaving any one class out moves \oursF within $0.46$ to $0.54\times$ the native solver (median over left-out classes $0.54$), \oursA within $0.49$ to $0.56$, \oursS within $0.59$ to $0.68$ and the predictor within $0.64$ to $0.76$, against $0.90$ to $1.04$ for Lu \& Yan and $1.04$ to $1.14$ for Fan et al. The most influential class for every method is set systems, on which every learned start is far faster than the native solver; the upper ends of these ranges are the values without it. Failures are few: 3 of 382 LPs for \oursS on \Tfif, all wrong-objective terminations on badly scaled load-balancing LPs, against 4 time limits for the native solver. These panels ran without the optimality check, which removes all three (Table~\ref{tab:checked}); in the \Ttw panel, where it runs inline, none of our learned methods fails. The two classes of \Tfif that \oursS loses (Figure~\ref{fig:percat}), facility location and telecom network design ($1.2\times$), are ones where the primal simplex from a primal-feasible all-slack basis is cheap per pivot and the harm gate does not reliably predict this.

\textbf{Windows and summed time on \Tfif.}
\begin{itemize}
\item \oursS against the native solver by window: W1 (164 LPs) 0.621 \ci{0.47}{0.80}, W2 (145) 0.635 \ci{0.50}{0.80}, W3 (73) 0.407. Summed time: 0.762 pooled, and 0.37, 1.03 and 0.74 by window.
\item Leaving out any one class moves the class-balanced figure of \oursS within 0.513 to 0.612 (predictor with harm gate: 0.503 to 0.600); the median over left-out classes is 0.629. \oursS is faster than the native solver on 215 of 382 LPs, the predictor alone on 248.
\item Failures: native solver 4 (time limits), all-slack dual simplex 0, \oursS 3 (wrong objectives), predictor 5, Fan et al.\ 6, Lu \& Yan 4 to 5 (one unbounded status). On \Tei \oursS fails twice (one wrong objective, one solver error).
\end{itemize}

\IfFileExists{sections/generated/tab_windows_test20k.tex}{\begin{table}[htbp]
\centering\footnotesize
\caption{Charged time on \textsc{Test-20k} by window (the all-slack dual simplex's time at generation). Time: class-balanced geometric-mean ratio to the native solver with 95\% bootstrap interval (2,000 draws over classes, then generator instances) and classes won/lost; iterations: class-balanced geometric-mean ratio to the all-slack dual simplex of one plus the iteration count; summed: total charged time over that of the native solver; failures count at twice the window's time limit (120, 600 and 3,600\,s). Numbers are averaged within LP over two timing repetitions.}
\label{tab:windows-test20k}
\setlength{\tabcolsep}{5pt}
\begin{tabular}{@{}l c c c c r@{}}
\toprule
Method & time [95\% CI] & W/L & iterations & summed & failures\\
\midrule
\multicolumn{6}{@{}l}{W1 ($\le 10$\,s): 1,052 LPs, 24 classes}\\
\quad native solver (reference) & 1.00 & -- & 0.86 & 1.00 & 0\\
\quad all-slack dual & 1.04 \ci{0.85}{1.29} & 10/14 & 1.00 & 0.91 & 0\\
\quad stock solver, presolve on & 1.13 \ci{1.09}{1.18} & 1/23 & 0.86 & 1.03 & 0\\
\quad classical crash & 1.84 \ci{1.28}{2.72} & 9/15 & 1.27 & 2.65 & 1\\
\quad Fan et al. & 1.09 \ci{0.88}{1.37} & 14/10 & 0.46 & 1.20 & 0\\
\quad Lu \& Yan & 1.06 \ci{0.88}{1.29} & 13/11 & 0.48 & 1.19 & 0\\
\quad predictor & 0.74 \ci{0.59}{0.93} & 19/5 & 0.31 & 0.76 & 0\\
\quad \oursS & 0.67 \ci{0.56}{0.80} & 22/2 & 0.29 & 0.57 & 0\\
\quad \oursA & 0.60 \ci{0.50}{0.72} & 22/2 & 0.25 & 0.52 & 0\\
\quad \oursF & 0.55 \ci{0.46}{0.66} & 23/1 & 0.29 & 0.58 & 0\\
\midrule
\multicolumn{6}{@{}l}{W2 (10--120\,s): 267 LPs, 12 classes}\\
\quad native solver (reference) & 1.00 & -- & 0.77 & 1.00 & 0\\
\quad all-slack dual & 1.38 \ci{1.02}{1.94} & 5/7 & 1.00 & 1.08 & 0\\
\quad stock solver, presolve on & 1.00 \ci{0.93}{1.05} & 3/9 & 0.76 & 0.98 & 0\\
\quad classical crash & 1.62 \ci{0.91}{3.13} & 4/8 & 1.18 & 1.55 & 0\\
\quad Fan et al. & 1.21 \ci{0.94}{1.61} & 6/6 & 0.80 & 1.04 & 0\\
\quad Lu \& Yan & 0.70 \ci{0.28}{1.35} & 7/5 & 0.63 & 0.82 & 0\\
\quad predictor & 0.53 \ci{0.20}{1.01} & 10/2 & 0.48 & 0.72 & 0\\
\quad \oursS & 0.46 \ci{0.18}{0.80} & 8/4 & 0.44 & 0.61 & 0\\
\quad \oursA & 0.33 \ci{0.15}{0.61} & 11/1 & 0.32 & 0.45 & 0\\
\quad \oursF & 0.36 \ci{0.15}{0.64} & 11/1 & 0.35 & 0.51 & 0\\
\midrule
\multicolumn{6}{@{}l}{W3 ($> 120$\,s): 37 LPs, 5 classes}\\
\quad native solver (reference) & 1.00 & -- & 0.71 & 1.00 & 0\\
\quad all-slack dual & 1.85 \ci{0.94}{6.23} & 3/2 & 1.00 & 1.58 & 0\\
\quad stock solver, presolve on & 0.80 \ci{0.49}{1.14} & 3/2 & 0.54 & 0.95 & 0\\
\quad classical crash & 0.78 \ci{0.30}{4.16} & 4/1 & 0.49 & 0.71 & 0\\
\quad Fan et al. & 1.09 \ci{0.38}{4.48} & 4/1 & 0.64 & 0.99 & 0\\
\quad Lu \& Yan & 0.79 \ci{0.26}{3.96} & 4/1 & 0.53 & 0.93 & 0\\
\quad predictor & 0.80 \ci{0.25}{4.31} & 4/1 & 0.47 & 0.98 & 0\\
\quad \oursS & 0.82 \ci{0.33}{2.70} & 3/2 & 0.53 & 0.81 & 0\\
\quad \oursA & 0.22 \ci{0.12}{0.42} & 5/0 & 0.18 & 0.23 & 0\\
\quad \oursF & 0.33 \ci{0.12}{1.29} & 4/1 & 0.21 & 0.30 & 0\\
\midrule
\multicolumn{6}{@{}l}{all windows: 1,356 LPs, 25 classes}\\
\quad native solver (reference) & 1.00 & -- & 0.85 & 1.00 & 0\\
\quad all-slack dual & 1.07 \ci{0.87}{1.30} & 12/13 & 1.00 & 1.26 & 0\\
\quad stock solver, presolve on & 1.10 \ci{1.06}{1.16} & 3/22 & 0.85 & 0.98 & 0\\
\quad classical crash & 1.76 \ci{1.19}{2.60} & 9/16 & 1.22 & 1.34 & 1\\
\quad Fan et al. & 1.10 \ci{0.89}{1.35} & 15/10 & 0.49 & 1.04 & 0\\
\quad Lu \& Yan & 0.94 \ci{0.70}{1.22} & 14/11 & 0.47 & 0.90 & 0\\
\quad predictor & 0.68 \ci{0.50}{0.91} & 20/5 & 0.32 & 0.83 & 0\\
\quad \oursS & 0.61 \ci{0.46}{0.76} & 23/2 & 0.30 & 0.68 & 0\\
\quad \oursA & 0.51 \ci{0.38}{0.65} & 23/2 & 0.25 & 0.37 & 0\\
\quad \oursF & 0.49 \ci{0.37}{0.63} & 24/1 & 0.28 & 0.44 & 0\\
\bottomrule
\end{tabular}
\end{table}
}{}
\IfFileExists{sections/generated/tab_percat_test20k.tex}{\begin{table}[htbp]
\centering\footnotesize
\caption{Per-class charged time on \textsc{Test-20k} relative to the native solver: geometric mean over the class's LPs of the per-LP ratio, all three windows, failures at twice the time limit, lowest entry per row in bold. The last rows give the class-balanced geometric mean over the 25 classes, the classes won/lost against the native solver, and the failures. Every learned arm runs the optimality check inline; numbers are averaged within LP over two timing repetitions.}
\label{tab:percat-test20k}
\setlength{\tabcolsep}{4pt}
\begin{tabular}{@{}lr cccc ccc@{}}
\toprule
& & \multicolumn{3}{c}{baselines} & & \multicolumn{3}{c}{\ours}\\
\cmidrule(lr){3-5}\cmidrule(l){7-9}
Class & LPs & all-slack & Fan et al. & Lu \& Yan & predictor & Shift & Alg & Full\\
\midrule
bin packing & 50 & \textbf{0.46} & 0.96 & 1.00 & 0.83 & 0.54 & 0.50 & 0.74\\
crop planning & 11 & 1.05 & 0.74 & 0.80 & \textbf{0.68} & 0.92 & 0.70 & 0.69\\
energy & 104 & 0.74 & 0.79 & 0.81 & 0.58 & 0.64 & 0.53 & \textbf{0.50}\\
facility location & 37 & 2.92 & 2.07 & 1.63 & 1.39 & 1.19 & 1.02 & \textbf{0.91}\\
graph optimization & 156 & 1.26 & 1.09 & 1.12 & 0.76 & 0.69 & \textbf{0.36} & 0.63\\
hub location & 195 & 2.54 & 2.85 & 2.60 & 2.07 & 0.96 & 0.91 & \textbf{0.77}\\
inventory & 77 & 0.48 & 0.39 & 1.07 & 0.38 & 0.42 & 0.38 & \textbf{0.32}\\
inverse optimization & 143 & 1.00 & 0.73 & 0.56 & 0.38 & 0.36 & 0.28 & \textbf{0.27}\\
knapsack & 39 & 1.00 & 1.51 & 0.99 & 0.69 & 0.69 & 0.67 & \textbf{0.50}\\
land use & 13 & 1.02 & 1.66 & 1.56 & 0.80 & 0.91 & 0.80 & \textbf{0.59}\\
load balancing & 52 & 0.98 & 1.42 & 1.85 & 0.60 & \textbf{0.37} & 0.38 & 0.42\\
maritime inventory routing & 26 & 1.00 & 0.97 & 0.87 & 0.78 & \textbf{0.75} & \textbf{0.75} & 0.86\\
multi-commodity flow & 26 & 1.00 & 0.93 & 1.00 & 0.81 & 0.99 & 0.58 & \textbf{0.39}\\
nurse scheduling & 26 & 0.99 & 1.30 & 1.06 & 0.81 & \textbf{0.71} & \textbf{0.71} & 0.73\\
operating room scheduling & 65 & 1.00 & 0.84 & 0.80 & 0.70 & 0.76 & \textbf{0.61} & \textbf{0.61}\\
revenue management & 39 & 1.01 & 0.83 & 0.89 & 0.76 & 0.77 & 0.64 & \textbf{0.54}\\
scheduling & 6 & 0.82 & 0.90 & 0.70 & 0.36 & 0.33 & \textbf{0.30} & 0.37\\
set systems & 26 & 1.00 & 0.80 & 0.08 & \textbf{0.05} & 0.06 & \textbf{0.05} & \textbf{0.05}\\
stochastic programming & 26 & 1.01 & 0.81 & 0.71 & 0.52 & 0.69 & 0.53 & \textbf{0.37}\\
supply chain & 46 & 2.37 & 1.72 & 1.79 & 1.05 & 0.92 & 0.86 & \textbf{0.85}\\
telecom network design & 25 & 1.85 & 3.79 & 2.58 & 3.76 & \textbf{1.42} & 1.68 & 1.69\\
transportation & 13 & 0.57 & 0.53 & 0.42 & 0.52 & 0.52 & 0.31 & \textbf{0.29}\\
traveling salesman & 116 & 2.59 & 2.77 & 2.24 & 1.68 & 0.95 & 0.92 & \textbf{0.86}\\
unit commitment & 13 & 0.99 & 0.83 & 0.56 & 0.54 & 0.55 & \textbf{0.20} & \textbf{0.20}\\
vehicle routing & 26 & 0.57 & 0.72 & 0.67 & 0.43 & 0.38 & 0.44 & \textbf{0.34}\\
\midrule
class-balanced mean & 1,356 & 1.07 & 1.10 & 0.94 & 0.68 & 0.61 & 0.51 & \textbf{0.49}\\
classes won/lost & & 12/13 & 15/10 & 14/11 & 20/5 & 23/2 & 23/2 & 24/1\\
failures & & 0 & 0 & 0 & 0 & 0 & 0 & 0\\
\bottomrule
\end{tabular}
\end{table}
}{}
\IfFileExists{sections/generated/tab_actions_test20k.tex}{\begin{table}[htbp]
\centering\scriptsize
\caption{Actions taken on \textsc{Test-20k}, as shares of LPs (\%). Top: action groups by pricing rule and shift $\delta$; for \oursS the native column is the harm gate's fallback, and for all three systems $\delta = 0$ under DSE includes the LPs that the novelty gate returns to the unshifted start. Bottom: the shift chosen by the shift policy of \oursS before the harm gate, and the gate's fallback rate. Numbers are averaged within LP over two timing repetitions; actions and gate decisions are the same in every repetition.}
\label{tab:actions-test20k}
\setlength{\tabcolsep}{4pt}
\begin{tabular}{@{}ll r cc ccc ccc@{}}
\toprule
& & & & & \multicolumn{3}{c}{DSE pricing} & \multicolumn{3}{c}{Devex pricing}\\
\cmidrule(lr){6-8}\cmidrule(l){9-11}
System & window & LPs & native & primal & $\delta < 0$ & $\delta = 0$ & $\delta > 0$ & $\delta < 0$ & $\delta = 0$ & $\delta > 0$\\
\midrule
\oursF & W1 & 1,052 & 8.1 & 7.2 & 0.0 & 12.9 & 16.7 & 41.8 & 0.0 & 13.2\\
 & W2 & 267 & 16.1 & 25.8 & 0.0 & 7.9 & 27.3 & 20.2 & 0.0 & 2.6\\
 & W3 & 37 & 0.0 & 29.7 & 0.0 & 0.0 & 24.3 & 29.7 & 0.0 & 16.2\\
 & all & 1,356 & 9.4 & 11.5 & 0.0 & 11.6 & 19.0 & 37.2 & 0.0 & 11.2\\
\midrule
\oursA & W1 & 1,052 & 13.9 & 12.6 & 8.6 & 32.1 & 32.8 & -- & -- & --\\
 & W2 & 267 & 2.6 & 49.1 & 12.7 & 6.7 & 28.8 & -- & -- & --\\
 & W3 & 37 & 0.0 & 59.5 & 0.0 & 0.0 & 40.5 & -- & -- & --\\
 & all & 1,356 & 11.3 & 21.1 & 9.1 & 26.3 & 32.2 & -- & -- & --\\
\midrule
\oursS & W1 & 1,052 & 27.8 & -- & 18.1 & 16.3 & 37.8 & -- & -- & --\\
 & W2 & 267 & 44.2 & -- & 11.6 & 27.7 & 16.5 & -- & -- & --\\
 & W3 & 37 & 13.5 & -- & 0.0 & 83.8 & 2.7 & -- & -- & --\\
 & all & 1,356 & 30.6 & -- & 16.3 & 20.4 & 32.7 & -- & -- & --\\
\bottomrule
\end{tabular}

\vspace{6pt}
\begin{tabular}{@{}l r ccccccccc c@{}}
\toprule
& & \multicolumn{9}{c}{shift $\delta$ chosen by the shift policy} & harm gate\\
\cmidrule(lr){3-11}
window & LPs & $-3$ & $-2$ & $-1$ & $-0.5$ & $0$ & $+0.5$ & $+1$ & $+2$ & $+3$ & to native\\
\midrule
W1 & 1,052 & 20.2 & 0.0 & 0.6 & 3.7 & 27.9 & 1.5 & 0.0 & 0.0 & 46.1 & 27.8\\
W2 & 267 & 12.0 & 0.0 & 0.0 & 0.0 & 49.8 & 0.0 & 0.0 & 0.0 & 38.2 & 44.2\\
W3 & 37 & 0.0 & 0.0 & 0.0 & 0.0 & 97.3 & 0.0 & 0.0 & 0.0 & 2.7 & 13.5\\
all & 1,356 & 18.0 & 0.0 & 0.4 & 2.9 & 34.1 & 1.2 & 0.0 & 0.0 & 43.4 & 30.6\\
\bottomrule
\end{tabular}
\end{table}
}{}
\IfFileExists{sections/generated/tab_overhead.tex}{\begin{table}[htbp]
\centering\scriptsize
\caption{Charged overhead per LP on \textsc{Test-20k}: median wall time in milliseconds and median share of the LP's total charged time (\%), by window. Decode includes the policy (with its novelty gate) and the crash; the harm gate (features and model) is charged on every LP, including those it hands to the native solver, which then need no setBasis. Total is the sum of the five components. Numbers are averaged within LP over two timing repetitions.}
\label{tab:overhead}
\setlength{\tabcolsep}{4pt}
\begin{tabular}{@{}l rr rr rr rr@{}}
\toprule
 & \multicolumn{2}{c}{W1} & \multicolumn{2}{c}{W2} & \multicolumn{2}{c}{W3} & \multicolumn{2}{c}{all}\\
\cmidrule(lr){2-3}\cmidrule(lr){4-5}\cmidrule(lr){6-7}\cmidrule(l){8-9}
Component & ms & \% & ms & \% & ms & \% & ms & \%\\
\midrule
\multicolumn{9}{@{}l}{predictor}\\
\quad features & 9.1 & 1.28 & 9.0 & 0.05 & 10.4 & 0.01 & 9.1 & 0.89\\
\quad inference & 6.7 & 0.99 & 8.9 & 0.05 & 11.1 & 0.01 & 7.0 & 0.69\\
\quad decode & 8.6 & 1.31 & 10.6 & 0.06 & 19.0 & 0.01 & 9.0 & 0.90\\
\quad setBasis & 5.5 & 0.70 & 9.6 & 0.05 & 8.8 & 0.01 & 6.3 & 0.49\\
\quad harm gate & -- & -- & -- & -- & -- & -- & -- & --\\
\quad total overhead & 30.1 & 4.38 & 39.4 & 0.21 & 53.6 & 0.04 & 32.5 & 3.19\\
\midrule
\multicolumn{9}{@{}l}{\oursS}\\
\quad features & 8.8 & 1.50 & 9.0 & 0.07 & 10.4 & 0.01 & 8.9 & 1.07\\
\quad inference & 6.7 & 1.30 & 8.9 & 0.07 & 11.1 & 0.01 & 7.1 & 0.88\\
\quad decode & 19.7 & 3.84 & 23.1 & 0.19 & 39.9 & 0.02 & 20.5 & 2.61\\
\quad setBasis & 2.3 & 0.46 & 4.0 & 0.01 & 8.7 & 0.01 & 2.8 & 0.27\\
\quad harm gate & 10.0 & 1.91 & 11.0 & 0.09 & 16.5 & 0.01 & 10.2 & 1.23\\
\quad total overhead & 49.8 & 9.38 & 56.9 & 0.44 & 94.2 & 0.07 & 52.0 & 6.34\\
\midrule
\multicolumn{9}{@{}l}{\oursA}\\
\quad features & 8.7 & 1.67 & 9.0 & 0.11 & 10.3 & 0.05 & 8.9 & 1.16\\
\quad inference & 6.7 & 1.45 & 9.0 & 0.11 & 11.0 & 0.03 & 7.1 & 0.94\\
\quad decode & 11.8 & 2.43 & 13.1 & 0.16 & 21.3 & 0.05 & 12.1 & 1.63\\
\quad setBasis & 3.9 & 0.66 & 8.7 & 0.09 & 8.9 & 0.06 & 5.2 & 0.50\\
\quad harm gate & -- & -- & -- & -- & -- & -- & -- & --\\
\quad total overhead & 31.9 & 6.44 & 41.8 & 0.49 & 59.2 & 0.19 & 34.0 & 4.60\\
\midrule
\multicolumn{9}{@{}l}{\oursF}\\
\quad features & 8.7 & 1.98 & 9.0 & 0.07 & 10.2 & 0.02 & 8.9 & 1.31\\
\quad inference & 6.7 & 1.55 & 9.0 & 0.07 & 11.1 & 0.02 & 7.1 & 1.03\\
\quad decode & 11.6 & 2.83 & 13.0 & 0.11 & 21.3 & 0.04 & 11.9 & 1.72\\
\quad setBasis & 4.6 & 0.87 & 7.7 & 0.04 & 9.5 & 0.02 & 5.5 & 0.66\\
\quad harm gate & -- & -- & -- & -- & -- & -- & -- & --\\
\quad total overhead & 32.4 & 7.59 & 38.7 & 0.30 & 57.8 & 0.10 & 34.2 & 4.82\\
\bottomrule
\end{tabular}
\end{table}
}{}
\IfFileExists{sections/generated/tab_percat_test50k.tex}{\begin{table}[htbp]
\centering\footnotesize
\caption{Per-class charged time on \textsc{Test-50k} relative to the native solver: geometric mean over the class's LPs of the per-LP ratio, all three windows, failures at twice the time limit, lowest entry per row in bold. This panel runs without the optimality check; Fan et al.\ comes from a separate paired panel on the same LPs with its own native-solver solves. The last rows give the class-balanced geometric mean, the classes won/lost against the native solver, and the failures.}
\label{tab:percat-test50k}
\setlength{\tabcolsep}{4pt}
\begin{tabular}{@{}lr cccc c@{}}
\toprule
& & \multicolumn{3}{c}{baselines} & & \\
\cmidrule(lr){3-5}
Class & LPs & all-slack & Fan et al. & Lu \& Yan & predictor & \oursS\\
\midrule
bin packing & 16 & \textbf{0.34} & 0.76 & 0.77 & 1.38 & 0.50\\
crop planning & 4 & 1.04 & 0.58 & 0.63 & \textbf{0.46} & 0.65\\
energy & 21 & 0.56 & 0.54 & 0.64 & \textbf{0.47} & 0.50\\
facility location & 8 & 2.21 & 2.27 & 1.69 & 1.30 & \textbf{1.22}\\
graph optimization & 41 & 1.15 & 0.95 & 0.92 & 0.58 & \textbf{0.57}\\
hub location & 60 & 2.52 & 2.65 & 2.41 & 1.77 & \textbf{0.92}\\
inventory & 16 & 0.38 & 0.27 & 0.71 & \textbf{0.24} & 0.26\\
inverse optimization & 38 & 1.00 & 0.88 & 0.71 & \textbf{0.44} & \textbf{0.44}\\
knapsack & 12 & 0.99 & 0.95 & 0.99 & \textbf{0.56} & 0.60\\
land use & 4 & 1.01 & 1.14 & 0.99 & \textbf{0.54} & 0.57\\
load balancing & 16 & 1.17 & 3.24 & 2.80 & 1.53 & \textbf{0.72}\\
maritime inventory routing & 8 & 1.00 & 0.91 & 0.93 & \textbf{0.75} & 0.83\\
multi-commodity flow & 4 & 1.04 & 1.07 & 0.91 & \textbf{0.76} & 0.79\\
nurse scheduling & 8 & 0.78 & 1.13 & 0.82 & \textbf{0.70} & 0.91\\
operating room scheduling & 20 & 1.01 & 0.73 & 0.69 & \textbf{0.59} & 0.63\\
revenue management & 12 & 0.98 & 0.79 & 0.85 & \textbf{0.77} & 0.95\\
scheduling & 3 & 0.45 & 0.49 & 0.41 & 0.16 & \textbf{0.15}\\
set systems & 7 & 1.00 & 0.75 & 0.03 & \textbf{0.02} & \textbf{0.02}\\
stochastic programming & 8 & 1.03 & 1.09 & 0.96 & \textbf{0.78} & 0.79\\
supply chain & 16 & 2.85 & 2.06 & 1.80 & 0.98 & \textbf{0.94}\\
telecom network design & 8 & 2.42 & 4.03 & 4.94 & 4.20 & \textbf{1.17}\\
transportation & 4 & 0.31 & 0.39 & \textbf{0.25} & \textbf{0.25} & 0.26\\
traveling salesman & 36 & 3.36 & 3.22 & 2.82 & 1.97 & \textbf{0.94}\\
unit commitment & 4 & 0.99 & 0.89 & 0.69 & 0.66 & \textbf{0.63}\\
vehicle routing & 8 & 0.47 & 0.69 & 0.56 & 0.36 & \textbf{0.35}\\
\midrule
class-balanced mean & 382 & 0.99 & 1.02 & 0.85 & 0.61 & \textbf{0.53}\\
classes won/lost & & 12/13 & 15/10 & 19/6 & 19/6 & 23/2\\
failures & & 0 & 7 & 5 & 5 & 3\\
\bottomrule
\end{tabular}
\end{table}
}{}

\subsection{Training label of the decision policy}
\label{app:objective}

The label ablation of \S\ref{sec:periter} (Figure~\ref{fig:mech}a, Table~\ref{tab:labels}) trains three shift policies that differ only in the label, each relative to shift 0: log instructions (I; a failure counts as $\log(2 \cdot 120\,\text{s} \cdot 5\times 10^9)$), $\log(1 + \text{iterations})$ (N; failure $\log(1 + 2\times10^6)$), and log measured seconds under 64 concurrent workers (T; failure $\log 240$). Everything else follows the shift policy's recipe (nine actions, hidden width 128, dropout 0.1, AdamW at $10^{-3}$, 100 epochs, batch 256, clip 3, seeds 0 to 2, epoch chosen on a random fifth of the training LPs). The labels have standard deviations 0.350, 0.329 and 0.351 and within-LP Spearman correlations of 0.98 (I and T), 0.85 (I and N) and 0.83 (N and T). The charged panel has the 381 W1 LPs of \Dtw; starts are prepared once and are identical for the three policies, and the solves run on one exclusive node with three repetitions in seeded order, failures at 240\,s (PAR2), and no failures.

\begin{table}[h]
\centering\footnotesize
\caption{Training-label ablation on the 381 W1 LPs of \Dtw: three shift policies that differ only in the label (I: instructions, N: iterations, T: wall time under load). Charged time, iterations and summed time relative to the predictor (class-balanced, three timing repetitions), then paired contrasts between the policies.}
\label{tab:labels}
\setlength{\tabcolsep}{5pt}
\begin{tabular}{@{}lccc@{}}
\toprule
Label & time [95\% CI] & iterations & summed time\\
\midrule
instructions (I) & 0.927 \ci{0.877}{0.978} & 0.943 & 0.918\\
iterations (N) & 0.953 \ci{0.899}{1.008} & 0.932 & 0.947\\
wall time (T) & 0.951 \ci{0.907}{0.993} & 0.955 & 0.963\\
\midrule
I relative to N & 0.973 \ci{0.934}{1.010} & \multicolumn{2}{l}{by repetition 0.9735, 0.9735, 0.9719}\\
\quad on the 202 LPs without a twin in training & 0.939 \ci{0.886}{0.984} & &\\
I relative to T & 0.974 \ci{0.951}{0.996} & &\\
\bottomrule
\end{tabular}
\end{table}

Split by Eq.~(\ref{eq:decomp}), the gain of the instruction policy over the iteration policy lies entirely in the time per iteration (the split is in Appendix~\ref{app:periter}): on the offline label tables, the iteration policy's choices sit 4\% above their LP's mean work per iteration and the instruction policy's 3\% below it. The offline tables show the same without timing noise:
\begin{itemize}
\item On the selection split of \Dtw, the shift with the fewest iterations costs $1.06\times$ \ci{1.03}{1.09} the instruction-optimal one, and on 38\% of LPs choosing by iterations costs at least 1\% more than the instruction-optimal shift.
\item On the same split, the cheapest action in instructions among the iteration minimizers costs $1.058\times$ \ci{1.029}{1.094} the instruction oracle over the nine shifts and $1.096\times$ \ci{1.051}{1.159} over the 42 decoder starts; on the 7,699 training LPs, $1.048\times$ \ci{1.028}{1.072}. At least 1\% instruction regret is unavoidable for any iteration-optimal choice on 33 to 55\% of LPs.
\end{itemize}

\textbf{Label ablation for the full policy.} The shift actions leave little room for the label to matter (the offline gap above is about 6\%). The 20 actions of \oursF include actions that change the time per iteration, so we also trained the \oursF policy on log iteration counts relative to the unshifted start (clipped at $\pm3$, failures at $+3$), with its exact recipe (ten seeds, 150 epochs, the same held-out fifth of the training LPs). An overlap label was not possible, since the cost tables record no per-action overlap. Table~\ref{tab:label-full} in \S\ref{sec:periter} scores each policy's chosen actions in instructions relative to the predictor's unshifted start. Trained on iterations, the policy removes 15 to 18\% more iterations but routes only 5 to 15\% of LPs to Devex instead of 29 to 60\%, sending more to primal simplex; its iterations then cost 1.2 to 1.7 times as much, for 40\% more instructions on short solves and a tie on longer ones. This is the 20-action counterpart of the shift-policy result, and the per-iteration term is where the instruction label earns its gain.

\subsection{A per-class lookup control}
\label{app:lookup}

A natural question about a per-LP policy is whether a per-class rule would do as well. The lookup plays, for every LP of a class, the action with the lowest mean relative log instructions for that class in the training cost tables (the 20 actions of \oursF plus the classical crash). It uses the LP's class label at deployment, which none of our methods does, and it has no answer for a class absent from training. To evaluate it on \Ttw without new solves, we restrict it to actions that were run as full methods on that panel: the unshifted predictor start with DSE or with Devex, the native solver with or without Devex, the all-slack dual simplex and the classical crash. It then chooses the predictor's start with Devex for 12 training classes, with default pricing for 9, the native solver with Devex for hub location, telecom and TSP, the all-slack dual for bin packing and the classical crash for unit commitment. Table~\ref{tab:lookup} gives the results. The lookup is significantly faster than \oursS, and its two-action version (learned start or the native solver) ties it, which is further evidence that \oursS's gain over the predictor comes from deciding when to use the start at all. \oursF is 9\% faster than the lookup overall (not significant), $0.76\times$ its time on W2 and $0.50\times$ on W3, and $0.88\times$ \ci{0.79}{0.98} its instructions; offline on held-out training LPs, the \oursF policy takes $0.92\times$ \ci{0.81}{1.03} (W1) and $0.93\times$ \ci{0.83}{1.02} (W2) the instructions of the unrestricted 20-action lookup. The per-LP policy therefore pays off mainly on longer solves, and the lookup's strength itself supports the paper's premise that the algorithm and pricing decisions, trained on measured work, are where the time is.

\begin{table}[h]
\centering\footnotesize
\caption{Per-class lookup control on \Ttw (two timing repetitions averaged): charged time relative to the native solver by window, class-balanced, and paired contrasts. The lookup uses class labels; our methods do not.}
\label{tab:lookup}
\setlength{\tabcolsep}{4pt}
\begin{tabular}{@{}lcccc@{}}
\toprule
Method & W1 (1,052) & W2 (267) & W3 (37) & all [95\% CI]\\
\midrule
per-class lookup, restricted actions & 0.585 & 0.479 & 0.663 & 0.537 \ci{0.406}{0.672}\\
per-class lookup, two actions (start or native) & 0.658 & 0.480 & 0.802 & 0.608 \ci{0.460}{0.744}\\
best single action on training (predictor start) & 0.738 & 0.530 & 0.802 & 0.683 \ci{0.497}{0.899}\\
\oursS & 0.671 & 0.457 & 0.822 & 0.614 \ci{0.466}{0.763}\\
\oursA & 0.600 & 0.331 & 0.217 & 0.510\\
\oursF & 0.551 & 0.362 & 0.334 & 0.489 \ci{0.369}{0.619}\\
\midrule
\multicolumn{5}{@{}p{0.95\linewidth}@{}}{Paired, \oursF relative to the restricted lookup: W1 $0.942$ \ci{0.846}{1.049}, W2 $0.756$ \ci{0.640}{0.874}, all $0.911$ \ci{0.819}{1.016}. \oursS relative to the restricted lookup: all $1.144$ \ci{1.033}{1.265}; relative to the two-action lookup: $1.010$ \ci{0.929}{1.091}.}\\
\bottomrule
\end{tabular}
\end{table}

\subsection{Seeds}
\label{app:seeds}

\textbf{Predictor seeds.} A second seed of stage 2 matches the first within 0.01 in overlap, precision and paired time at matched steps (0.873 vs.\ 0.881 against the model after stage 1, at step 296). A second seed of stage 1 at 2 epochs reaches overlap 0.7575 and precision 0.7248, against 0.7608 and 0.7194. A replicate of stage 2 (from the same stage-1 initialization, on the deduplicated training set with a new seed), with its own action costs and shift policy, on \Dtw: its predictor takes $1.065\times$ \ci{1.01}{1.13} the deployed predictor's instructions, its shift policy $0.890\times$ \ci{0.83}{0.95} its own predictor's, and the replicate system $1.014\times$ the deployed system in instructions and $1.007\times$ in time.

\textbf{Full-pipeline seeds.} Stage 1, stage 2 (on the deduplicated training set), the nine-shift cost table and the gated shift policy were retrained with seeds 101, 202 and 303. On the 283-LP selection split of \Dtw, which also selected their checkpoints and epochs, the three predictors reach overlap 0.785 to 0.793 and $0.822$ to $0.840\times$ the native solver in charged time (geometric mean 0.830, standard deviation of the log ratio 0.011). With their own gated shift policies they take $0.733$, $0.770$ and $0.761\times$ (mean 0.755, standard deviation 0.026), with iterations $0.513\times$ and instructions $0.702\times$, against $0.728\times$ for the deployed predictor with its policy. Offline, the three shift policies take $0.849$, $0.903$ and $0.904\times$ their own predictor's instructions, and a policy fitted to another seed's predictor gains little ($0.977\times$), so policies are specific to their predictor. Single policy seeds give $0.896 \pm 0.013$ (5 seeds), and they are portable across hardware: six single policies re-costed on a different CPU model agree to 0.001.

On \Ttw (Table~\ref{tab:seeds}; all windows, both repetitions averaged, check inline), the pipelines of seeds 202 and 303 ran as paired panels anchored on the native solver, the deployed predictor and \oursS; seed 101 is part of the main panel. The seed spread of the full pipeline on the test set ($\pm 0.003$) is two orders of magnitude below the width of the class and LP intervals.

\begin{table}[h]
\centering\footnotesize
\caption{Full-pipeline seeds on \Ttw: charged time relative to the native solver (class-balanced, two repetitions averaged, check inline), and the shift policy's time relative to its own predictor. The deployed predictor, trained on 20k training data that include duplicates, takes 0.683, and \oursS 0.614.}
\label{tab:seeds}
\setlength{\tabcolsep}{6pt}
\begin{tabular}{@{}lccc@{}}
\toprule
& seed 101 & seed 202 & seed 303\\
\midrule
predictor & 0.714 & 0.711 & 0.708\\
\quad + own shift policy & 0.642 & 0.647 & 0.649\\
\quad + own shift policy, check and harm gate & 0.621 & 0.616 & 0.619\\
shift policy relative to its predictor & 0.90 & 0.91 \ci{0.85}{0.96} & 0.92 \ci{0.87}{0.96}\\
\bottomrule
\end{tabular}
\end{table}

\subsection{Deployment cost}
\label{app:cost}

\textbf{Inference on the CPU.} All panels run the predictor's forward pass on a GPU and charge its time. On CPU cores alone, an optimized forward pass (fused message passing, bf16 matrix units, compiled kernels) takes a median 93\,ms per LP of the 20k tier on 16 cores, against 8\,ms on the GPU, with starts whose solve cost is within 1 to 2\% of the GPU starts. Charged, this makes \oursS 29 to 30\% slower on solves under 10\,s (on a 400-LP W1 subset of \Ttw, from $0.70\times$ to $0.90\times$ the native solver), 3 to 5\% slower on longer solves, and 14\% slower on \Tfif pooled; a CPU deployment would need a forward pass of about 30\,ms to keep the loss on short solves under 5\%.

\textbf{Training cost.} Predictor: 8.8 H100-hours for stage 1 plus 2.5 hours for stage 2, including validation solves. Action costs for the shift policy: the 131,706-solve table took about 60 core-hours, of which the 7,699 training LPs use 69k solves; for the twenty actions of \oursF, 184k solves and 373 core-hours (Appendix~\ref{app:periter}). Labels of the harm gate: 7.7 core-hours. Published methods: Fan et al.\ 1.5 GPU-hours for stage 1 at 0.43 times our stage-1 samples; Lu \& Yan 2.5 GPU-hours plus 86 CPU-hours of late-path labels at 0.33 times; both second stages at our stage-2 budget, and Lu \& Yan additionally at 12 stage-1 epochs. Test generation: an estimated 58 CPU-hours for the 20k test sets and about 100 for the 50k and 80k sets.

\section{Published methods: reimplementation and the public benchmark}
\label{app:literature}

This appendix describes our reimplementations of the two published methods, checks that our decoder and training choices do not handicap them on our data (Table~\ref{tab:litfair}), and reports the public benchmark of \citet{lu2025learning} under its own protocol (Appendix~\ref{app:luyanbench}).

\textbf{Reimplementation.}
\begin{itemize}
\item Fan et al.\ (2023): the released GCN with the released balanced cross-entropy loss, 0.2M parameters, Adam at $10^{-3}$ with weight decay $10^{-4}$, learning rate $\times 0.1$ at one and two thirds of the budget, one LP per step. The paper's factorization-based basis repair is replaced by our crash on the basic log-odds, so that every learned method shares one decoder.
\item Lu \& Yan (2025): tripartite graph, 1.2M parameters, the basic-variable loss, the multi-level loss with $\theta = 0.9$ and ten levels from solver replays under iteration limits from the all-slack start, and label preprocessing (the paper's separate feasibility term is not implemented). Inputs are clipped at 10 with global normalization and a gradient-norm clip of 1, which the unclipped bound features (up to $7\times10^6$) require; LPs without late-path labels (15\%) use the terminal level only.
\item Both: stage 1 on our 114,874 LPs of \Trs (Fan et al.\ 5.1 epochs, Lu \& Yan 4.0, capped by label cost; a further Fan et al.\ continuation to 7.3 epochs gained nothing on its own loss), stage 2 on exactly the predictor's stage-2 data and budget (the full-budget Lu \& Yan model on the deduplicated set), and the checkpoint chosen by the method's own loss on the selection split of \Dtw. The full budget is 12 stage-1 epochs.
\end{itemize}
Differences from the papers that a reader should weigh: one model over all classes rather than one per dataset; presolved, stored LPs rather than raw files; our decoder rather than theirs. On their own benchmarks the papers report iteration ratios of roughly 0.3 to 0.6 against the default start (0.04 to 0.93 across individual datasets); our reimplementations reach 0.46 to 0.52 against the all-slack dual simplex on \Ttw and \Tfif (0.44 to 0.57 with \Tei).

\textbf{Decoder fairness and seed variance on our data.} Table~\ref{tab:litfair} compares decoders and training seeds on \Dtw (fully charged). Decoding the published models' scores with the papers' top-$m$ decoders instead of our crash makes them slower relative to the predictor, except for Fan et al.\ on W2, where the difference is not significant; with top-$m$ repair, Lu \& Yan takes $2.41\times$ the predictor's time on W1. Models trained across 26 classes select many dependent structural columns, HiGHS alters the installed basis on about 90\% of such starts, and the probe installs of the repair variant cost 0.13 to 1.2\,s; the triangular crash avoids both. A second training seed of each reimplementation, at full stage-1 budget and decoded by the crash, changes its charged time by $0.94$ to $1.08\times$ (paired; every interval includes 1), at most about 8\%, far below the $1.2$ to $2.9\times$ gap to the predictor; every interval against the predictor excludes 1 in the predictor's favor. Fan et al.'s failures are the known load-balancing wrong optima (4 to 5 per seed on W1, 1 on W2).

\begin{table}[h]
\centering\footnotesize
\caption{Decoders and training seeds of the reimplemented published methods on \Dtw: charged time relative to our predictor (class-balanced; W1: 392 LPs, W2: 69 LPs). Top-$m$ installs the $m$ highest-scoring columns and lets HiGHS repair the set. n.s.: not significantly different from the crash.}
\label{tab:litfair}
\setlength{\tabcolsep}{6pt}
\begin{tabular}{@{}llcc@{}}
\toprule
Method & Decoder or seed & W1 & W2\\
\midrule
\multicolumn{4}{@{}l}{Decoder (same model)}\\
Lu \& Yan & our crash & 1.41 & 1.21\\
 & top-$m$ & 1.82 & 1.72\\
Fan et al. & our crash & 1.50 & 2.84\\
 & top-$m$ & 1.69 & 2.27 (n.s.)\\
\midrule
\multicolumn{4}{@{}l}{Training seed (full budget, our crash)}\\
Lu \& Yan & first / second seed & 1.40 / 1.43 & 1.20 / 1.30\\
Fan et al. & first / second seed & 1.46 / 1.38 & 2.85 / 2.90\\
\bottomrule
\end{tabular}
\end{table}

The same holds on the W1 and W2 LPs of the larger test sets (check inline; charged time relative to the native solver on \Tfif / \Tei):
\begin{itemize}
\item Lu \& Yan at full budget: $1.02$ / $1.03$ with the crash, against $1.39$ / $1.46$ with the top-$m$ decoder;
\item Fan et al.: $0.99$ / $1.00$ against $1.12$ / $1.19$, with three top-$m$ solves returning a solver error after installation and counted as failures;
\item our predictor with the top-$m$ decoder: $0.73$ / $0.83$, indifferent to the decoder ($1.04\times$ \ci{0.88}{1.26} relative to the crash).
\end{itemize}
The crash-decoded published methods therefore do not understate the baselines in our setting, and the predictor's margin holds under every decoder. On the per-dataset benchmark below the opposite holds, which is why top-$m$ is the right decoder there.

\textbf{Results on \Ttw.} At full budget, Fan et al.\ removes 51\% of the iterations on \Ttw and is 10\% slower than the native solver; Lu \& Yan removes 53\% and is 6\% faster, and a second training seed of each gives $1.04\times$ and $0.96\times$. Their starts go through the same crash as ours with the same 40 to 60\,ms of overhead, so the difference lies in the scores; with the papers' own top-$m$ decoders they are slower still on our data ($1.16$ to $1.19\times$ on \Ttw). Our predictor takes $0.73\times$ \ci{0.64}{0.82} the time of Lu \& Yan, paired, with $0.67\times$ its iterations and $1.08\times$ its time per iteration (Table~\ref{tab:decomp}).

\subsection{The public benchmark of Lu \& Yan}
\label{app:luyanbench}

\textbf{Protocol.} The benchmark has four public datasets (Mirp1, Mirp2, Anonymous, Load balancing). We follow the protocol of \citet{lu2025learning}: one model per dataset, a 7:3 split (ours, since theirs is unpublished; 9 to 30 test LPs per dataset, 80 in total), presolved LPs, and fully charged time against the dual simplex with presolve off. Our predictor, pretrained on our synthetic data, is fine-tuned per dataset; the published models are trained from scratch per dataset as in the papers. The top-$m$ repair decoder fixes the number of basic variables at $m$ by rank and repairs singular selections as in \citet{fan2023smart}; top-$m$ without repair installs the top-$m$ set through \texttt{setBasis} and lets HiGHS repair it.

\textbf{Per-dataset decision policy.} The action set has 39 actions: our crash with the nine shifts and the top-$m$ decoder with seven shifts, each under default pricing or Devex; top-$m$ repair under either rule; the unshifted crash and top-$m$ starts with primal simplex; and the native solver with three settings. Every training LP is solved once under every action, with instruction counts as labels. For each dataset, 5-fold cross-validation over the training LPs selects among the best single action, a policy trained from scratch with the loss of Eq.~(\ref{eq:loss}), and a policy initialized from \oursF and fine-tuned; the winner is refit on all training LPs. No novelty gate is used, since every test LP comes from the dataset's own generator. With 19 to 70 training LPs per dataset, the selection keeps top-$m$ repair on Mirp1 and a constant action (top-$m$ with primal simplex) on Load balancing, and trained policies on the other two.

\textbf{Results.} Table~\ref{tab:luyanbench} in \S\ref{sec:luyan} gives the per-dataset results. The decision policy takes $0.58\times$ the iterations and $0.72\times$ \ci{0.68}{0.75} the time of the dual simplex, against $1.03\times$ for our reimplementation of Lu \& Yan and $0.92\times$ for Fan et al. Paired, it takes $0.69\times$ \ci{0.64}{0.75} the time of Lu \& Yan (Mirp1 $0.48$, Anonymous $0.64$, Load balancing $0.74$, Mirp2 $1.03$) and $0.78\times$ that of Fan et al. Against our fine-tuned predictor with top-$m$ repair it takes $0.84\times$ the time at the same iterations ($1.02\times$), so the gain comes from cheaper decoding and cheaper iterations. The best single action per dataset alone reaches $0.73\times$, so on datasets this small most of the gain comes from choosing among decoders, algorithms and pricing rules by measured cost, not from a per-LP choice; initializing the policy from \oursF, which was trained on our data, hurts ($0.86\times$).

\textbf{Raw (unpresolved) files.} On the unpresolved versions of the four datasets (same protocol, 79 test LPs), the decision policy takes $0.45\times$ the iterations and $0.585\times$ \ci{0.53}{0.65} the time of the dual simplex, against $0.58$ / $0.82$ for our reimplementation of Lu \& Yan, $0.58$ / $0.75$ for Fan et al.\ and $0.45$ / $0.65$ for our fine-tuned predictor with top-$m$ repair; paired, it takes $0.71\times$ \ci{0.63}{0.81} the time of Lu \& Yan and $0.78\times$ that of Fan et al. Lu \& Yan report $0.57$ / $0.73$ on raw files. On these files, however, the solver with its own presolve takes $0.29\times$ the time of the dual simplex without presolve, faster than every warm start, ours included; learned starts are useful here only when the LP is solved without presolve or presolve has already been applied, as on our data.

\textbf{Additional panel and zero-shot model.} A separate panel on the same split, without the decision policy, gives the following. Our predictor without fine-tuning (trained only on our synthetic data) reaches $0.91\times$ the iterations and $1.07\times$ the time with our crash and does not benefit from the top-$m$ decoders ($1.11$ to $1.15\times$ time). Mirp2 is at parity for every method on time, and it is the one dataset where the published top-$m$ repair beats our predictor on iterations ($1.24\times$ against it). Paired at each method's best decoder, our fine-tuned predictor takes:
\begin{itemize}
\item on the presolved files, $0.83$ \ci{0.78}{0.88} the iterations and $0.81$ \ci{0.76}{0.87} the time of Lu \& Yan, and $0.84$ / $0.89$ those of Fan et al.;
\item on the raw files, $0.78$ / $0.79$ those of Lu \& Yan and $0.79$ / $0.87$ those of Fan et al. There our predictor reaches $0.45\times$ the iterations and $0.65\times$ the time of the dual simplex, against the $0.57$ / $0.73$ that Lu \& Yan report, and the solver with presolve on takes $0.29\times$ the time of the dual simplex without presolve, faster than every warm start.
\end{itemize}

\textbf{Reproduction.} Our reimplementation reproduces the reported results of Lu \& Yan on the load-balancing and anonymous datasets and those of Fan et al.\ on Mirp1 (raw iterations $0.536$ against the paper's $0.539$), but not Lu \& Yan's Mirp1 figure. The paper reports $0.32\times$ the iterations; our reimplementation reaches $0.62$ to $0.70\times$ across five retrains that vary the input clipping and the preprocessing term (seed noise $\pm 0.04$), and only $0.63$ to $0.68\times$ on its own training LPs, while our own per-dataset predictor reaches $0.35\times$ there. We ruled out:
\begin{itemize}
\item the HiGHS version (1.7.2, 1.9.0 and 1.13.1 give identical iteration counts);
\item the algorithm choice and presolve;
\item the data sizes, which match the paper's table to within one;
\item the labels and the split.
\end{itemize}
What remains is something in the unreleased model or training beyond the paper's text, or an unreported evaluation detail. We report both figures. Full tables are released with the code.

\textbf{Reading of the public benchmark.} Three caveats apply: the split is ours, since theirs is unpublished, with 9 to 30 test LPs per dataset; our predictor is pretrained on our synthetic data and then fine-tuned per dataset, whereas theirs is trained from scratch per dataset; and zero-shot, without fine-tuning, no learned start beats the default there. On their smallest dataset, whose LPs solve in about a second, our time is dominated by decoding and basis installation. On these datasets the decoder matters more than the model: every per-dataset model gains from the top-$m$ repair, and our triangular crash is the wrong decoder there (the same fine-tuned predictor reaches only $0.83\times$ the iterations and $0.97\times$ the time with it). The shift policy adds at most 0.04 there, so the gain over the published methods comes from the predictor and the choice of decoder.

\section{Generalization: larger LPs, held-out variants and new classes}
\label{app:hc}

This appendix gives the full versions of Table~\ref{tab:general} (Table~\ref{tab:checked} for part a, Table~\ref{tab:transfer} for part b) and the details behind \S\ref{sec:generalization}: larger LPs, held-out variants and dual forms, unseen classes, and adaptation to new classes (Appendix~\ref{app:fewshot}).

\begin{table}[h]
\centering\footnotesize
\caption{Full version of Table~\ref{tab:general}b. Charged-time ratios, class-balanced, with 95\% intervals. \Tad: the shift policy adapted with $n$ labeled LPs per new class, against the unadapted policy, whose novelty gate falls back to the predictor's start on these classes. \Tva is a subset of the \Ttw panel.}
\label{tab:transfer}
\setlength{\tabcolsep}{5pt}
\begin{tabular}{@{}lcc@{}}
\toprule
Method & vs.\ native solver & vs.\ reference\\
\midrule
\multicolumn{3}{@{}l}{\Tnc: 234 LPs, 6 classes absent from training (reference: predictor)}\\
predictor & 1.42 \ci{0.95}{2.30} & --\\
\oursS, harm gate off (policy abstains on 234/234) & 1.42 \ci{0.96}{2.32} & 1.00\\
\oursS & 1.21 \ci{0.97}{1.62} & 0.85 \ci{0.62}{1.08}\\
Fan et al.\ (full budget) & 1.37 \ci{0.96}{2.09} & --\\
Lu \& Yan (full budget) & 0.94 \ci{0.49}{1.70} & --\\
classical crash & 0.67 \ci{0.30}{1.18} & --\\
\midrule
\multicolumn{3}{@{}l}{\Tad: 616 LPs, 5 classes absent from training (reference: unadapted policy = predictor's start)}\\
unadapted policy & 1.14 \ci{0.82}{1.68} & 1.00\\
\oursS (with harm gate) & 1.05 \ci{0.94}{1.18} & 0.93\\
adapted with $n = 10$ labeled LPs per class & 1.05 \ci{0.74}{1.49} & 0.93 \ci{0.85}{1.00}\\
adapted with $n = 25$ labeled LPs per class & 1.04 \ci{0.73}{1.48} & 0.91 \ci{0.84}{0.98}\\
online bandit over strategies, no labels & 0.93 \ci{0.79}{1.09} & 0.82 \ci{0.64}{0.97}\\
\midrule
\multicolumn{3}{@{}l}{\Tva: 296 LPs, 10 held-out formulations, 6 held-out dual forms (ref.: deployed predictor)}\\
predictor trained without them & 0.71 \ci{0.33}{1.41} & 1.33 \ci{1.13}{1.64}\\
\quad with its own shift policy, the deployed harm gate and check & 0.56 \ci{0.30}{0.86} & 1.05 \ci{0.67}{1.50}\\
deployed \oursF (trained on them) & 0.38 \ci{0.19}{0.63} & 0.71 \ci{0.49}{0.91}\\
\bottomrule
\end{tabular}
\end{table}

\begin{table}[h]
\centering\footnotesize
\caption{Full version of Table~\ref{tab:general}a: the W1 and W2 LPs of \Tfif (309 of 382 LPs) and \Tei (63 of 81). Charged time relative to the native solver, class-balanced, with 95\% intervals, failures in parentheses, and summed time. Every learned method runs the optimality check inline; gate evaluation is charged. The table omits W3, where Table~\ref{tab:main} shows the largest gains.}
\label{tab:checked}
\setlength{\tabcolsep}{4pt}
\begin{tabular}{@{}lcccc@{}}
\toprule
& \multicolumn{2}{c}{\Tfif (309)} & \multicolumn{2}{c}{\Tei (63)}\\
\cmidrule(lr){2-3}\cmidrule(l){4-5}
Method & time [95\% CI] & summed & time [95\% CI] & summed\\
\midrule
all-slack dual simplex & 0.94 \ci{0.72}{1.22} & 0.87 & 0.90 \ci{0.63}{1.28} & 1.13\\
native solver with Devex & 1.06 \ci{0.89}{1.34} (6) & 1.63 & 1.07 \ci{0.93}{1.26} & 1.28\\
predictor & 0.65 \ci{0.48}{0.87} (1) & 0.72 & 0.83 \ci{0.54}{1.27} (1) & 1.70\\
predictor, Devex on every LP & 0.66 \ci{0.46}{0.96} & 0.89 & 0.78 \ci{0.48}{1.30} (1) & 2.29\\
\oursS & 0.58 \ci{0.46}{0.72} & 0.50 & 0.67 \ci{0.48}{0.88} & 0.84\\
\quad with class-blind harm gate & 0.67 \ci{0.52}{0.83} & 0.74 & 0.75 \ci{0.55}{1.01} & 0.95\\
\oursA & 0.60 \ci{0.45}{0.78} (2) & 0.61 & 0.83 \ci{0.53}{1.28} (1) & 1.71\\
\quad + harm gate & 0.55 \ci{0.43}{0.69} & 0.47 & 0.66 \ci{0.48}{0.86} & 0.83\\
\quad + class-blind harm gate & 0.63 \ci{0.49}{0.79} & 0.69 & 0.74 \ci{0.53}{1.01} & 0.95\\
\oursF & 0.56 \ci{0.43}{0.74} (1) & 0.60 & 0.81 \ci{0.53}{1.25} (1) & 1.71\\
\quad + harm gate & 0.54 \ci{0.42}{0.69} & 0.49 & 0.65 \ci{0.47}{0.85} & 0.84\\
\quad + class-blind harm gate & 0.60 \ci{0.48}{0.76} (1) & 0.70 & 0.73 \ci{0.53}{0.99} & 0.96\\
\bottomrule
\end{tabular}
\end{table}

\textbf{Larger LPs.} The novelty gate sends most LPs of the larger sets to the default action (208 of 309 on \Tfif, 61 of 63 on \Tei), where the policy never gets to choose the native solver, which is why the harm gate carries most of the gain there. On \Dfif, without any training data at that size, \oursF and \oursA are within a few percent of each other, and the predictor with Devex on every LP loses on W2 (Table~\ref{tab:devreadouts}). On the test sets (Table~\ref{tab:checked}), \oursF with the harm gate takes $0.98\times$ the time of \oursA with the harm gate on \Tfif and $0.99\times$ on \Tei, while Devex on every LP gives no significant gain over the predictor. Table~\ref{tab:checked} covers only W1 and W2 and so omits the longest solves, on which Table~\ref{tab:main} shows the largest gains (\oursS at $0.407\times$ the native solver on W3 of \Tfif; Appendix~\ref{app:results}); its ratios are therefore higher than those of Table~\ref{tab:main}.

\textbf{Held-out variants and dual forms on development data.} Between a new instance of a training class and a new class lies a third case: a new formulation of a training class, or the dual form of one. On a development panel (357 LPs of 10 held-out formulations and 6 held-out dual forms, with 330 control LPs of matching formulations and forms; instruction counts, fully charged) we compare a full pipeline trained without these formulations and dual forms to the deployed one by a difference in differences: the ratio of the two pipelines' costs on the held-out LPs divided by the same ratio on the controls. It is $1.23$ \ci{0.99}{1.55} over all 16, $1.06$ \ci{0.86}{1.29} for the dual forms, which retain about 93\% of the predictor's gain over the native solver, and $1.34$ \ci{1.02}{1.93} for the variants, which retain about 28\%. The pipeline trained without the variants still beats the native solver on 7 of the 10 and fails badly on three (a load-balancing variant at $5.1\times$ the deployed pipeline's cost); the novelty gate abstains on 98.7\% of the held-out LPs, so the policy abstains there. On \Tva (Table~\ref{tab:transfer}), the predictor trained without the held-out formulations and dual forms takes $1.33\times$ the deployed predictor's time, and on the other 1,060 LPs of \Ttw $1.04\times$, so the loss is specific to the held-out formulations and dual forms. Its own shift policy adds nothing there ($0.72\times$ the native solver without the harm gate, against $0.71\times$ for its predictor alone); the gain to $0.56\times$ comes from the deployed harm gate, which was fit on \Trtw and so did see these formulations and dual forms.

\textbf{Unseen classes.} On the six classes of \Tnc, every learned start of ours is slower than the native solver (the predictor by $1.4\times$; job-shop scheduling $2.7\times$, diet problems $3\times$), while the classical crash is faster ($0.67\times$). Fan et al.\ at full budget is at $1.37\times$; Lu \& Yan is near parity only because on the diet problems its scores nearly reduce to the crash's ordering. The novelty gate flags every LP but falls back to the predictor's unshifted start, which leave-one-class-out validation on the training classes preferred over the native solver ($1.00$ against $1.45\times$); on new classes the order reverses ($1.42\times$ against $1.00$). The harm gate recovers part of the loss but leaves \oursS at $1.21\times$ the native solver, because it too was fit on the training classes.

\textbf{Per-class results on \Tnc.} The predictor against the native solver on W1 and W2: blending 0.87, container loading 0.88, process planning 1.05, resilient network design 1.22, job-shop scheduling 2.74, diet problem 3.43. W3 holds only diet problems (35 LPs), where the predictor takes $3.0\times$ the native solver, the classical crash $0.09\times$, Lu \& Yan's original-budget start $0.10\times$ with iteration counts identical to the crash's (its scores reduce to the crash ordering there), and its full-budget start $0.20\times$. The harm gate falls back on 45\% of \Tnc LPs (51\% on W1 and W2); it recovers job-shop scheduling ($2.74 \to 1.27$) but hurts container loading ($0.88 \to 1.07$) and worsens summed time on W1 and W2 ($1.01 \to 1.10$). On the 800-LP held-out-class panel from which \Tad was carved, the picture is the same (predictor 1.169, with harm gate 1.064).

\subsection{Adaptation to new classes}
\label{app:fewshot}

\begin{figure}[tbp]
\centering
\includegraphics[width=0.5\linewidth]{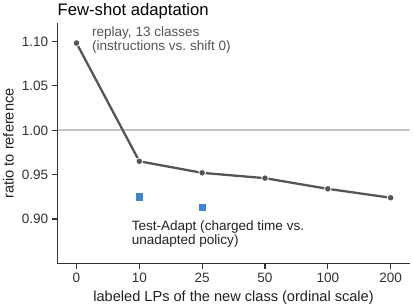}
\caption{Few-shot adaptation to unseen classes; the $n$ axis is ordinal. Line: replay on the training LPs with one class held out at a time, instructions relative to the unshifted start after adapting with $n$ labeled LPs of the held-out class (the 13 classes with long enough streams; Table~\ref{tab:fewshot-replay}). Squares: the charged panel on \Tad, time relative to the unadapted policy (the predictor's start).}
\label{fig:fewshot}
\end{figure}

We measured two ways to adapt to a new class from the solver's own cost: labeling a few LPs of the class, and an online bandit that needs no labels. Figure~\ref{fig:fewshot} summarizes the labeled route.

\textbf{Few-shot adaptation with labels.} With $n$ labeled LPs of a new class (every action solved once on each), the shift policy is retrained on the training data plus the $n$ new LPs, replicated to 10\% of the rows. In replay on the 7,699 training LPs, holding out one class at a time over three stream orders (Table~\ref{tab:fewshot-replay}), ten labeled LPs turn the harmful unadapted policy into a gain, and 50 bring it near the in-class reference of 0.929. Excluding LPs whose generator-instance cluster appears among the labeled LPs changes these numbers by less than 0.01. Charged for acquisition, a stream pays 1.02 at $n = 10$ and 1.28 at $n = 50$, so the option pays off only for streams of thousands of LPs per class.

\begin{table}[h]
\centering\footnotesize
\caption{Few-shot replay on the training LPs, one class held out at a time (three stream orders): instructions relative to shift 0 on the rest of the held-out class's stream after adapting with $n$ labeled LPs. At $n = 200$ only the 13 classes with at least 220 LPs qualify. The in-class reference is 0.929.}
\label{tab:fewshot-replay}
\setlength{\tabcolsep}{5pt}
\begin{tabular}{@{}lcccccc@{}}
\toprule
$n$ & 0 & 10 & 25 & 50 & 100 & 200\\
\midrule
13 classes with at least 220 LPs & 1.098 & 0.965 & 0.952 & 0.946 & 0.934 & 0.924\\
all classes with enough LPs & 1.076 & 0.944 & 0.932 & 0.929 & 0.912 & 0.924\\
\bottomrule
\end{tabular}
\end{table}

\textbf{Charged result on \Tad}. With 25 labeled LPs per class, the adapted policy with its novelty gate (whose reference set now includes the labeled LPs) takes $0.913\times$ \ci{0.842}{0.980} the time of the unadapted policy on the rest of the class; the unadapted policy's gate abstains on these classes, so this reference is the predictor's start (0.913 and 0.914 in two further timing repetitions), with instructions $0.920\times$ and iterations $0.998\times$, and is worse on none of the five classes (container loading 0.83, job-shop scheduling 0.84, process planning 0.93, blending 0.99, resilient network design 1.00). With $n = 10$ it takes $0.925\times$ \ci{0.847}{0.995}; without the novelty gate, the adapted policy takes $0.905\times$; against \oursS, which also has the harm gate, it takes $0.994\times$. Against the native solver it takes $1.038\times$ \ci{0.733}{1.482}, so it recovers most of the policy's in-distribution gain but does not win. Labeling costs 10.2 times the plain solves of the labeled LPs, so adaptation breaks even at about 4,900 LPs per class ($n = 25$) or 2,700 ($n = 10$).

\textbf{Labels that do not work.}
\begin{itemize}
\item Hindsight pseudo-labels from the solved LP (distance of each action's start to the final basis, structural recall, count gaps, bound flips) rank actions poorly: the median per-LP Spearman correlation with the true costs is 0.00 to 0.26. The best of them, distance to $B^\star$, picks actions worth $0.96\times$ shift 0 against the cost oracle's $0.79\times$; this is the overlap-oracle selector of \S\ref{sec:periter}. Few-shot adaptation on such labels gives 1.01 to 1.02.
\item A continual online variant on such labels gives 1.31, charged.
\end{itemize}

\textbf{Online bandit without labels.} The bandit keeps one learner per class stream with 13 arms: the nine shifts with dual simplex, the native solver with and without its choice of algorithm, the unshifted start with primal simplex, and the classical crash. Its context is 16 principal components of the penultimate features of a shift policy trained on \Trtw. A Bayesian linear model of each arm's log cost relative to the unshifted start, with a prior centered on that start and a prior covariance equal to the between-class covariance of the 25 training classes, is updated from the observed instruction count of the arm played. Arms are chosen by Thompson sampling after 20 warm-up LPs at the unshifted start, under a conservative budget that plays a nondefault arm only while an upper confidence bound on the summed relative cost of all such plays stays below 5\% of the stream so far; every stage is charged. Results:
\begin{itemize}
\item Development replay on the training LPs (instructions relative to shift 0): a noncontextual per-class bandit is limited by information (full-information follow-the-leader 0.961; realistic feedback 1.002); a contextual bandit on the trained shift policy's features reaches 0.975 \ci{0.942}{0.999} over full streams at zero acquisition cost, with gains arriving after about 200 LPs. The replay puts the break-even at roughly 125 LPs per class.
\item \Tad (the joint bandit over shifts and strategy arms, tuned on \Dtw, run on the 616 LPs as real charged streams in five orders): $0.885\times$ \ci{0.770}{1.013} the time of \oursS ($0.864\times$ in instructions), $0.820\times$ \ci{0.637}{0.968} the predictor's and $0.928\times$ \ci{0.793}{1.091} that of the native solver, over streams of 64 to 149 LPs per class; $1.03\times$ over the first 50 LPs of a stream and $0.85\times$ over the last 100; worst class process planning ($1.11\times$). A bandit over shifts alone gives $1.03\times$, so the held-out loss comes from the choice of strategy rather than of shift.
\item \Tnc (supplementary; streams of 26 to 91 LPs per class): negative. The joint bandit ends at $0.966\times$ the predictor's time and $1.129\times$ that of \oursS, because the streams end before it reaches the native solver's large wins on job-shop scheduling and diet problems (about $0.37\times$). For the diet problem, two of the three stream orders are included.
\end{itemize}

\section{Extended discussion: scope, limitations and deployment}
\label{app:discussion}

\textbf{When the system is meant to be used.} The setting in which learned warm starts pay off is a stream of new LPs from known problem classes that must be solved by simplex, for example because a basic solution is needed at the root of branch-and-bound, for cut generation or for sensitivity analysis, or because runs must be reproducible. In that setting the classes are known in advance, no earlier basis of the same LP is available, the training cost is paid once, and the relevant question is how much solve time is saved per LP. Our results support three practical recommendations. First, train on the classes that will be solved; the gains transfer to larger instances of those classes and to their dual forms, and partly to new formulations, but not to new classes. Second, deploy with the guardrails: the novelty gate keeps the policy from extrapolating, the harm gate is needed wherever the novelty gate takes decisions away from the policy (at sizes beyond the training range), and the optimality check is necessary on badly scaled classes, where it removed every wrong optimum we observed. Third, when a new class appears, expect no gain without in-class feedback; either label a few LPs of the class (about ten times their plain solve cost each) or let an online bandit adapt from its own solves over the first hundred or so LPs.

\textbf{Transfer to new problem classes.} None of the learned components transfers to problem classes absent from training, and neither did any selector or policy input we tried. On such classes the gated system takes $1.2\times$ the time of the native solver; the novelty gate recognizes the LPs as unfamiliar, but its fallback was chosen on training classes held out one at a time, which proved a poor stand-in for new classes. Adaptation recovers part of the loss, but needs labels or a stream of roughly 100 LPs or more per class.

\textbf{The pricing action.} The advantage of \oursF over \oursA is not significant in time on any test or development panel. What is consistent is how they reach their times: \oursF accepts more iterations that each cost less, and its gain concentrates on short solves. Only one alternative pricing rule was offered: Devex forced from the first iteration, against HiGHS' default, which starts with DSE and switches to Devex on its own in some cases. It was the only one of the 13 alternative solver settings in our probe that clearly helped.

\textbf{Scope of the theorem.} Theorem~\ref{thm:chain} is an existence result for solve-vector work under one pivot rule on a network-flow LP. The unbounded ratio holds in a cost model without per-iteration terms linear in the number of rows; an implementation with such a term, or with refactorization after a fixed number of updates, bounds the ratio of total work by a constant. HiGHS shows the separation in full with presolve and scaling off and Devex forced. Under its defaults, the exact initialization of the DSE weights, which a nonslack start pays unless it is already near optimal, masks most of it, and the defaults also change the paths. On our generated LPs that initialization is a median 0.4\% of the solve. The theorem motivates the design; the empirical case rests on the measurements of \S\ref{sec:empirical-gap}.

\textbf{Data, solver and hardware.} All main results use synthetic instances from one generator library and one solver. On the public benchmark of \citet{lu2025learning}, our decision policy beats the reimplemented baselines under their protocol, but it needs per-dataset training there, and on the benchmark's unpresolved files the solver's own presolve is faster than every warm start, and with 19 to 70 training LPs per dataset most of its gain comes from choosing a single action per dataset. Timing depends on hardware; we report paired ratios within panels and record instruction counts beside every time. Our timings charge a GPU forward pass; on CPU cores alone the forward pass is about ten times slower, which costs 29 to 30\% on solves under 10\,s (Appendix~\ref{app:cost}).

\textbf{Cost of the labels.} The cost table requires solving every training LP under every action, about 16 times the cost of solving each LP once for the twenty actions of \oursF. This is a one-time cost per set of classes, but it limits how many actions the policy can consider.

\section{Reproducibility statement}
\label{app:repro}

Code to reproduce all experiments and the trained models are available at \url{https://github.com/flixpar/learned-simplex-warm-start}.
All LPs are generated by an open-source generator library available at \url{https://github.com/flixpar/SyntheticLPs.jl}. The datasets and the timing protocol are described in Appendix~\ref{app:data}, and the training seeds in Appendix~\ref{app:seeds}. The solver is HiGHS 1.13.1 with the options in Table~\ref{tab:hyper}, which also lists every setting of the predictor, the policies and the gates. The decoder is given as Algorithm~\ref{alg:crash} and the deployment procedure as Algorithm~\ref{alg:deploy}. Timing follows the protocol of \S\ref{sec:setup} and Appendix~\ref{app:data}, and every time is recorded beside a hardware-independent instruction count. The published baselines were reimplemented as described in Appendix~\ref{app:literature}.

\end{document}